\documentclass[11 pt]{article}
\usepackage{amsmath,amsthm,amssymb,amsfonts}
\usepackage{graphicx}

\newcounter{Th}[section] \newcounter{Lm}[section] \newcounter{Ca}[section] \newcounter{Prop}[section]
\newcounter{ThA}
\newcounter{Problem}[section] \newcounter{Remark}[section] \newcounter{Example}[section]
\newcounter{Def}[section]\newcounter{Assum}[section]

\def\theTh{\arabic{section}.\arabic{Th}}
\def\theThA{\Alph{ThA}}
\def\theLm{\arabic{section}.\arabic{Lm}}
\def\theCa{\arabic{section}.\arabic{Ca}}

\def\theRemark{\arabic{section}.\arabic{Remark}}
\def\theExample{\arabic{section}.\arabic{Example}}
\def\theDef{\arabic{section}.\arabic{Def}}

\newenvironment{Th}[1][\relax]
    {\medspace\refstepcounter{Th}{\bf Theorem \theTh.}\ \it}
    {\rm\medspace}

\newenvironment{Lm}[1][\relax]
    {\medspace\refstepcounter{Lm}{\bf Lemma \theLm.}\ \it}
    {\rm\medspace}

\newenvironment{Ca}[1][\relax]
    {\medspace\refstepcounter{Ca}{\bf Corollary  \theCa.}\ \it}
    {\rm\medspace}

\newenvironment{Remark}[1][\relax]
    {\medspace\refstepcounter{Remark}{\bf Remark \theRemark.}\rm\ }
    {\medspace}

\newenvironment{Example}[1][\relax]
    {\medspace\refstepcounter{Example}{\bf Example \theExample.}\rm\ }
    {\medspace}

\newenvironment{Def}[1][\relax]
    {\medspace\refstepcounter{Def}{\bf Definition \theDef.}\rm\ }
    {\medspace}

\makeatother

\numberwithin{equation}{section}

\begin{document}

\author{\bfseries\large A.~I.~Tyulenev\thanks
{Moscow Institute of Physics and Technology
(State University). E-mail: tyulenev-math@yandex.ru. }}
\title{Some new function spaces of variable smoothness}
\maketitle

\section{Introduction}
The present paper is concerned with new modifications of Besov-type function spaces of variable smoothness, which are generalizations of the spaces $\widetilde{B}^{l}_{p,q}(\mathbb{R}^{n}, \{t_{k}\})$ of~\cite{Ty}.

The function spaces of variables smoothness (Besov-type and Lizorkin--\allowbreak Triebel-type spaces) and various generalizations thereof have been extensively studied.
We only mention the papers \cite{Be3}, \cite{Be4}, \cite{Be2}, \cite{KeVy}, \cite{KeDis}, \cite{KeWavelet}, \cite{Dachun}, \cite{Moura}, \cite{Ul} (and abundant references given therein).

It is interesting to note that the majority of studies on this subject have been concerned with spaces of variable smoothness consisting of distributions from the space $S'(\mathbb{R}^{n})$.
In this connection, the Littlewood--\allowbreak Paley theory and the machinery of Fourier analysis become basic research tools.

We  say that a weight sequence (defining the variable smoothness) $\{s_{k}\}=\{s_{k}(\cdot)\}_{k=0}^{\infty}$ lies in $ Y^{\alpha_{3}}_{\alpha_{1},\alpha_{2}}$ if,
for $\alpha_{3}\geq 0$, $\alpha_{1},\alpha_{2} \in \mathbb{R}$,

\smallskip
1) $\frac{1}{C_{1}} 2^{\alpha_{1}(k-l)} \le \frac{s_{k}(x)}{s_{l}(x)} \le C_{1} 2^{\alpha_{2}(k-l)}$, $l \le k \in \mathbb{N}_{0}$, $x \in \mathbb{R}^{n}$;
\vskip-30pt
\begin{gather}
\noalign{}\label{1.1}
\end{gather}
\vskip-15pt
2) $ s_{k}(x) \le C_{2} s_{k}(y)(1+2^{k}|x-y|)^{\alpha_{3}}$, $k \in \mathbb{N}_{0}$, $x,y \in \mathbb{R}^{n}$,

\smallskip
\noindent the constants $C_{1},C_{2} > 0$ in \eqref{1.1} are independent of both indexes $k,l$ and points~$x,y$.

In what follows we shall need the standard decomposition of unity. Let $B^{n}$ be the unit ball of $\mathbb{R}^{n}$, $\Psi_{0} \in S(\mathbb{R}^{n})$, $\Psi_{0}(x)=1$
for $x \in B^{n}$, $\operatorname{supp}\Psi_{0} \subset 2B^{n}$. For $j \in \mathbb{N}$, we set $\Psi_{j}(x):=\Psi_{0}(2^{-j}x)-\Psi_{0}(2^{-j+1}x)$, $x \in \mathbb{R}^{n}$.

In \cite{KeVy}, \cite{KeDis}, \cite{KeWavelet}, \cite{Moura} the Besov spaces of variable smoothness were defined as follows (here we indicate only the case of constant integration exponents).

\begin{Def}
\label{Def1.1}
Let $p,q \in (0,\infty]$, $\alpha_{1},\alpha_{2} \in \mathbb{R}$, $\alpha_{3} \geq 0$, $\{s_{k}\} \in Y^{\alpha_{3}}_{\alpha_{1},\alpha_{2}}$.
By $B^{\{s_{k}\}}_{p,q}(\mathbb{R}^{n})$ we shall denote the space of all distributions $f \in S'(\mathbb{R}^{n})$ with finite quasi-norm
\begin{equation}
\label{1.2}
\|f|B^{\{s_{k}\}}_{p,q}(\mathbb{R}^{n})\|:=\|s_{j}F^{-1}(\Psi_{j}F[f])|l_{q}(L_{p}(\mathbb{R}^{n}))\|.
\end{equation}
\end{Def}

In \eqref{1.2}, $F$ and $F^{-1}$ denote, respectively, the direct and inverse Fourier transform.
Formally replacing in Definition \ref{Def1.1} the weight sequence $\{s_{k}\}$ by the sequence $\{2^{ks}\gamma\}$ with  $s>0$, $\gamma \in A_{\infty}(\mathbb{R}^{n})$
we obtain the definition of the weighted Besov space (see \cite{HaSch}, \cite{Ry}) with Muckenhoupt weight.

Mention also should be made of the works \cite{HN}, \cite{Dachun}, \cite{Ul}, in which the axiomatic approach to function spaces (of both constant and
variable smoothness) was developed. Instead of the base space $L_{p}(\mathbb{R}^{n})$, a~study was made of a~more general function space equipped with norm \eqref{1.2} and satisfying
a~certain set of axioms. The spaces examined in \cite{Dachun}, \cite{Ul} include, as a~particular case, the scale of spaces of variable smoothness of \cite{KeDis},~\cite{KeWavelet}

In our opinion, it is also of interest to study the spaces of variable smoothness whose elements are not distributions, but rather  functions that are locally integrable in some power.
Such spaces were actively studied by O.\,V.~Besov. We only indicate the papers \cite{Be3}, \cite{Be4}, \cite{Be2} (and the references given therein).
It is worth noting that the aforementioned papers employed the classical methods of theory of functions and that the norm on a~space of functions of variable smoothness
was defined \textit{ab initio} with the help of classical differences.

Theorems on characterization of various function spaces of variable smoothness (and their generalizations) were put forward in \cite{KeVy}, \cite{Dachun},
where the ball means were used;
originally the norm on these spaces was defined using the Littlewood--\allowbreak Paley decomposition (\cite{KeVy}) or with the help of
Peetere maximal functions~\cite{Dachun}.
It was also assumed that a weight sequence $\{s_{k}\}$ lies in $Y^{\alpha_{3}}_{\alpha_{1},\alpha_{2}}$ with the additional assumption
\begin{equation}
\label{ogranicenie}
0 < \alpha_{1} \le \alpha_{2} < l.
\end{equation}

Besov \cite{Be4}, \cite{Be2}  studied the spaces of variable smoothness with $p,q \in (1,\infty)$.
It was also assumed that the weight sequence $\{s_{k}\} \in ^{\text{\rm loc}}Y^{\alpha_{3}}_{\alpha_{1},\alpha_{2}}$ under condition \eqref{ogranicenie}.
$Y^{\alpha_{3}}_{\alpha_{1},\alpha_{2}}$ in that condition~2) is replaced by the condition

\smallskip
\noindent 2$'$)\ \ $s_{k}(x) \le 2^{\alpha_{3}} s_{k}(y) , \quad k \in \mathbb{N}_{0}, \quad |x-y| \le 2^{-k}$.\hfill (1.4)
\smallskip

\setcounter{equation}{4}

Clearly, $Y^{\alpha_{3}}_{\alpha_{1},\alpha_{2}} \subset ^{\text{\rm loc}}Y^{\widetilde{\alpha}_{3}}_{\alpha_{1},\alpha_{2}}$, where $\widetilde{\alpha}_{3}$ depends only on~$\alpha_{3}$
and the constant $C_{2}$ of~\eqref{1.1}. The weighted class $^{\text{\rm loc}}Y^{\alpha_{3}}_{\alpha_{1},\alpha_{2}}$ is strictly larger than the class
$Y^{\alpha_{3}}_{\alpha_{1},\alpha_{2}}$, because the former contains functions of an exponential rate of growth at infinity.

It is also worth mentioning that the methods of \cite{Be4}, \cite{Be2},  \cite{KeVy}, \cite{KeDis}, \cite{KeWavelet} \cite{Dachun}, \cite{Moura}, \cite{Ul}
utilized to prove various assertions about the spaces $B^{\{s_{k}\}}_{p,q}(\mathbb{R}^{n})$ were based on pointwise estimates of the weight sequence $\{s_{k}\}$.
This machinery was used in \cite{HN}, \cite{KeWavelet}, \cite{Dachun}, \cite{Ul} to establish atomic decomposition theorems (as well as results on molecular and wavelet
expansions) provided that the atoms from such a~decomposition have zero high order moments. The number of zero moments for such atoms is governed by the exponents
$\alpha_{1},\alpha_{2},\alpha_{3}$. Unfortunately, it is not possible to check these conditions in specific problems. For example, if
the high-order moments of the atoms from the decomposition of a~function $f:\mathbb{R}^{n} \to \mathbb{R}$ are zero, then in general we may not assert that
the corresponding moments of traces of these atoms on the hyperplane $\mathbb{R}^{n-1}$ are zero.
In~\cite{Moura} the trace problem for Besov spaces of various smoothness was solved with the help of the atomic decomposition theorem under certain constraints on
the weight sequence $\{s_{k}\}$. These constraints allow one to avoid testing that the atoms from the trace decompositions have zero moments.
It will be demonstrated in \S\,6 that these conditions can be substantially relaxed.

The analysis of definitions of Besov spaces of variable smoothness used in \cite{Be4}, \cite{Be2}, \cite{KeVy}, \cite{Dachun} shows that
in this papers a~fairly restrictive condition \eqref{ogranicenie} (in the case when these spaces consisted of functions locally integrable in some power).
This constraint is natural in the case when $s_{k} \equiv C_{k}$ for all $k \in \mathbb{N}$ ($C_{k}$ are positive constants),
for otherwise one needs to have recourse to the theory of distributions. In the case of variable smoothness this condition is quite rough.
The constraint $l > \alpha_{2}$ was used in showing that the norm on a~Besov space is independent on the difference order.
This conditions is also fairly rough in the variable smoothness setting.
Indeed, in~\cite{Ty} with $p,q \in (1,\infty)$ the author has put forward new modifications of Besov spaces of variable smoothness
$\widetilde{B}^{l}_{p,q}(\mathbb{R}^{d},\{t_{k}\})$ (in~\cite{Ty} these spaces were denoted by $\overline{B}^{l}_{p,q}(\mathbb{R}^{d},\{\gamma_{k}\})$)
and showed the space $\widetilde{B}^{l}_{p,p}(\mathbb{R}^{d},\{\gamma_{k}\})$ is the trace of the weighted Sobolev space
$\widetilde{W}_{p}^{l}(\mathbb{R}^{n},\gamma)$ on the plane of dimension $1 \le d < n$, provided that a~weight $\gamma \in A^{\text{\rm loc}}_{p}(\mathbb{R}^{n})$ .
Note that the weight sequence $\{\gamma_{k}\}$ lies in the weighted class $^{\text{\rm loc}}Y^{\alpha_{3}}_{\alpha_{1},\alpha_{2}}$.
However condition \eqref{ogranicenie} may fail to hold for the sequence $\{\gamma_{k}\}$ if the weight is ``sufficiently bad'' (see Remarks \ref{R4.2}, \ref{R4.4} below).

Clearly this calls for a more sophisticated approach towards the very concept of variable smoothness.
The definition of the weighted class $^{\text{\rm loc}}Y^{\alpha_{3}}_{\alpha_{1},\alpha_{2}}$ requires correction.
There is also a~need in new methods that are capable, in particular, of dealing with the space $\widetilde{B}^{l}_{p,p}(\mathbb{R}^{d},\{\gamma_{k}\})$,
which is the trace of the weighted Sobolev spaces with weight $\gamma \in A^{\text{\rm loc}}_{p}(\mathbb{R}^{n})$.
In doing so new methods should not depend upon the pointwise behaviour of the weight sequence $\{t_{k}\}$.

In the present paper we introduce, for $p,q,r \in (0,\infty]$, the Besov space of variable smoothness $\widetilde{B}^{l}_{p,q,r}(\mathbb{R}^{n},\{t_{k}\})$,
which is a~subtle modification of the space $\widetilde{B}^{l}_{p,q}(\mathbb{R}^{n},\{t_{k}\})$ of~\cite{Ty}.
The norm on this space is defined in terms of the difference relations  $\delta^{l}_{r}g$ (see \eqref{2.7}).
Here, the weight sequence  $\{t_{k}\}$  lies in the new weighted class $X^{\alpha_{3}}_{\alpha,\sigma,p}$ (see Definition \ref{Def2.4}).
For our purposes the weighted class $X^{\alpha_{3}}_{\alpha,\sigma,p}$ proves to be more subtle (and altogether more natural!)
than the class $^{\text{\rm loc}}Y^{\alpha_{3}}_{\alpha_{1},\alpha_{2}}$.
Necessary and sufficient conditions for a~sequence $\{t_{k}\}$  to lie in the weighted class $X^{\alpha_{3}}_{\alpha,\sigma,p}$  are expressed in terms
(in a~sense) integral estimates, rather than pointwise estimates.

It is worth pointing out that the differences  $\delta^{l}_{r}(Q^{n})f$  were used in \cite{HN}
for the purpose of construction of
equivalent norms on Besov-type and Lizorkin--\allowbreak Trie\-bel-type spaces.
However, our spaces $\widetilde{B}^{l}_{p,q,r}(\mathbb{R}^{n},\{t_{k}\})$ do not fit the axiomatics of~\cite{HN} due to the less restrictive constraints which we place
on a~weight sequence $\{t_{k}\}$. Such constructions were recently used by O.\,V.~Besov \cite{Be5}, \cite{Be6}
for the study of spaces of functions of zero smoothness.

For the study of the space $\widetilde{B}^{l}_{p,q,r}(\mathbb{R}^{n},\{t_{k}\})$ we shall adjust the methods of nonlinear spline approximation,
which were developed in~\cite{DeVore} for the study of classical Besov spaces. It is worth observing that
methods of nonlinear spline approximation have not been used for the study of function spaces of variable smoothness and hence may be of independent interest.
Using these methods we will be able to put forward certain theorems on equivalent norms on the spaces
$\widetilde{B}^{l}_{p,q,r}(\mathbb{R}^{n},\{t_{k}\})$ and prove the atomic decomposition theorem for these spaces.
 We shall also characterize the trace of the space $\widetilde{B}^{l}_{p,q,r}(\mathbb{R}^{n},\{t_{k}\})$.
 This result extends, for a~constant $p$, the results of \cite{HaSch},~\cite{Moura}.

The paper is organized as follows. In \S\,1 we give auxiliary results to be used in the analysis that follows.
In~\S\,2 we put forward some fundamental properties of the new spaces $\widetilde{B}^{l}_{p,q,r}(\mathbb{R}^{n},\{t_{k}\})$ and compare them
with the space $B^{\{s_{k}\}}_{p,q}(\mathbb{R}^{n})$.
In~\S\,3 we extend results of the paper~\cite{Ty} (this section may be looked upon as a~rationale
for the further constructions that follow).
In~\S\,4 we present the central results of the paper and, in particular, put forward the atomic decomposition theorem, which will be used
in \S\,\S\,5 and~6 to derive a~few embedding and trace theorems.

\section{Definitions and auxiliary results}

Throughout the following convention will be adopted.
The symbol~$C$ will be used to denote (different) `insignificant' constants in various estimates.
Sometimes, if it is required for purposes of exposition, we shall indicate the parameters on which some or other constant depends.

By definition, a weight function (a~weight) is a~measurable function $\gamma:\mathbb{R}^{n} \to (0,+\infty)$.
Given a~measurable set $E \subset \mathbb{R}^{n}$, we define $\gamma(E):=\int\limits_{E}\gamma(x)\,dx$.

Next, by $L_{p}(E)$ we denote the space of all equivalence classes (consisting of functions vanishing almost everywhere) equipped with the norm
\begin{gather*}
\|f|L_{p}(E)\|:=\biggl(\int\limits_{E}|f(x)|^{p}\,dx\biggr)^{\frac{1}{p}}, \qquad 1\le p <\infty,\\
\|f|L_{\infty}(E)\|:=\mbox{ess}\sup|f(x)|.
\end{gather*}

Given a measurable function $g:\mathbb{R}^{n} \longrightarrow \mathbb{R}$, a~measurable set~$E$, and a~weight~$\gamma$, we denote by $L_{p}(E,\gamma)$
the space of all equivalence classes (consisting of functions that coincide almost everywhere) and equip it with norm $\|g|L_{p}(E,\gamma)\|:=\|\gamma
g|L_{p}(E)\|$.

In what follows, $Q^{n}$ will denote an open cube in the space $\mathbb{R}^{n}$ with sides parallel to coordinate axes, $r(Q^{n})$ will denote the
side length of a~cube $Q^{n}$, and $|Q^{n}|$ will denote its $n$-dimensional Lebesgue measure.
For $\delta >0$, by $\delta Q^{n}$ we shall mean the cube, concentric with a~cube $Q^{n}$, with side length $r(\delta Q^{n}):=\delta r(Q^{n})$.
For $m=(m_{1},\dots,m_{n}) \in \mathbb{Z}^{n}$, $k \in \mathbb{Z}$, we let
$Q_{k,m}^{n}:=\prod\limits_{i=1}^{n}(\frac{m_{i}}{2^{k}},\frac{m_{i}+1}{2^{k}})$ denote
the open dyadic cube of side $2^{-k}$, $\widetilde{Q}_{k,m}^{n}:=\prod\limits_{i=1}^{n}[\frac{m_{i}}{2^{k}},\frac{m_{i}+1}{2^{k}})$. We also define
$I^{n}:=\prod\limits_{i=1}^{n}(-1,1)$. Also, $\delta B^{n}$ ($\delta S^{n}$) is the $n$-dimensional ball (sphere) of radius~$\delta$ centred at the origin.

For $x \in \mathbb{R}^{n}$, $E \subset \mathbb{R}^{n}$, we define $x + E := \{y \in \mathbb{R}^{n}: y=x+z, z \in E\}$.

V.\,S.~Rychkov \cite{Ry} introduced the class of weights $A^{\text{\rm loc}}_{p}(\mathbb{R}^{n})$, which generalizes the well-known Muckenhoupt class
$A_{p}(\mathbb{R}^{n})$ (for $1 < p \le \infty$).

\begin{Def} (\cite{Ry})
\label{Def2.1} Let $p \in (1,\infty)$, $a > 0$. Given a~weight $\gamma$ we say that
$\gamma \in A^{\text{\rm loc}}_{p}(\mathbb{R}^{n})$ if
$$
C^{\text{\rm loc}}_{\gamma,p,a}:=\sup\limits_{Q^{n}:r(Q^{n}) \le
a}\frac{1}{|Q^{n}|}\int\limits_{Q^{n}}\gamma(x)\,dx\biggl[\frac{1}{|Q^{n}|}\int\limits_{Q^{n}}\gamma^{-\frac{p'}{p}}(x)\,dx\biggr]^{\frac{p}{p'}}<
+\infty.
$$
\end{Def}

\begin{Def}(\cite{Sawano})
\label{Def2.2} Let $a > 0$. We say that a weight $\gamma \in
A^{\text{\rm loc}}_{1}(\mathbb{R}^{n})$ if there exists
a~constant $C^{\text{\rm loc}}_{\gamma,1,a} > 0$ independent of $Q^{n}$ such that, for all cubes of side length $r(Q^{n}) \le a$,
$$
\frac{1}{|Q^{n}|}\int\limits_{Q^{n}}\gamma(\widetilde{x})\,d\widetilde{x} \le
A \gamma(x) \ \ \mbox{ for a.e.\ } \ x \in Q^{n}.
$$
By $C^{\text{\rm loc}}_{\gamma,1,a}$ we shall mean the smallest constant $A$ satisfying the above inequality.
\end{Def}

\begin{Def}(\cite{Ry})
\label{Def2.3}
Let $a > 0$. We say that a~weight $\gamma \in
A^{\text{\rm loc}}_{\infty}(\mathbb{R}^{n})$ if, for some $\alpha \in (0,1)$,
$$
\sup\limits_{r(Q^{n}) \le a}\biggl(\sup\limits_{F \subset Q^{n}, |F| \geq \alpha|Q^{n}|}\frac{\gamma(Q^{n})}{\gamma(F)}\biggr) < \infty.
$$
\end{Def}

\begin{Remark}(\cite{Ry})
\label{R2.1}
If a weight $\gamma \in A^{\text{\rm loc}}_{\infty}(\mathbb{R}^{n})$, then there exists a~number $p_{0} \in [1,\infty)$ such that $\gamma \in A^{\text{\rm loc}}_{p_{0}}(\mathbb{R}^{n})$.
\end{Remark}

\begin{Remark}
\label{R2.2}
For $p \in (1,+\infty]$ the definition of the class
$A^{\text{\rm loc}}_{p}(\mathbb{R}^{n})$ is independent of the choice of the parameter~$a$.
For various $a > 0$ the constants $C^{\text{\rm loc}}_{\gamma,p,a}$ are estimated by each other~\cite{Ry}.
 One may show that a~similar result also holds for $A^{\text{\rm loc}}_{1}(\mathbb{R}^{n})$.
\end{Remark}

Given $f \in L_{1}^{\text{\rm loc}}(\mathbb{R}^{n})$, $a > 0$, we let $M_{\le a}f$ denote the local version of the Hardy--Littlewood maximal function,
$$
M_{\le a}f(x):=\sup\limits_{x \in Q^{n}, r(Q^{n}) \le a}\frac{1}{|Q^{n}|}\int\limits_{Q^{n}}|f(y)|\,dy.
$$

The next theorem generalizes the classical result of Muckenhoupt \cite{St} (see, for example, \S\,5.3, Theorem~1).

\begin{Th} ({see \cite{Ry}})
\label{Th2.1}
Let $p \in (1,\infty)$, $\gamma \in A_{p}^{\text{\rm loc}}(\mathbb{R}^{n})$, $a > 0$. Then
there exists a~constant $C=C(n,p,a,\gamma)>0$ such that
$$
\int\limits_{\mathbb{R}^{n}}\gamma(x)\left\{M_{\le a}[f](x)\right\}^{p}\,dx \le C \int\limits_{\mathbb{R}^{n}}\gamma(x)|f(x)|^{p}\,dx
$$
for all $f \in L_{p}(\mathbb{R}^{n},\gamma^{\frac{1}{p}})$.
\end{Th}

\begin{Th} {\rm (Hardy's inequality for sequences)}
\label{Th2.2}
Let $0 < q \le \infty$, $\mu \le q$, $\beta \geq 0$, and let $\{a_{k}\}$ be a~sequence of real numbers. Then
\begin{equation}
\label{2.1}
\left(\sum\limits_{k=0}^{\infty}2^{qk\beta}|b_{k}|^{q}\right)^{\frac{1}{q}} \le C \left(\sum\limits_{k=0}^{\infty}2^{qk\beta}|a_{k}|^{q}\right)^{\frac{1}{q}}
\end{equation}
where
\begin{equation}
\label{2.2}
|b_{k}| \le C \left(\sum\limits_{j=k}^{\infty}|a_{k}|^{\mu}\right)^{\frac{1}{\mu}}, \quad \mbox{ provided that } \beta > 0 \quad \mbox{ or }
\end{equation}
\begin{equation}
\label{2.3}
|b_{k}| \le C 2^{-k\lambda} \left(\sum\limits_{j=0}^{k}2^{j\mu\lambda}|a_{k}|^{\mu}\right)^{\frac{1}{\mu}}, \quad \mbox{ provided that } \lambda>\beta,
\end{equation}
the constant $C > 0$ being independent of the sequence $\{a_{k}\}$.
\end{Th}

In what follows we shall also need the following elementary fact.

\begin{Lm}
\label{Lm2.1}
Let $r \in (0,\infty]$, $f_{j} \in L^{\text{\rm loc}}_{r}(\mathbb{R}^{n})$ for $j \in \mathbb{N}_{0}$. Then, for $\mu \le \min \{1,r\}$,
\begin{equation}
\label{2.4}
\|\sum\limits_{j=1}^{\infty}f_{j}|L_{r}(\mathbb{R}^{n})\| \le \left(\sum\limits_{j=1}^{\infty}\|f_{j}|L_{r}(\mathbb{R}^{n})\|^{\mu}\right)^{\frac{1}{\mu}}.
\end{equation}
\end{Lm}

The proof easily follows from the monotonicity in~$q$ of the $l_{q}$-norm.

Let $l \in \mathbb{N}$, $r \in (0,\infty]$, $\Omega$ be a domain in $\mathbb{R}^{n}$. For a~function
$g \in L_{r}^{\text{\rm loc}}(\Omega)$, $h \in \mathbb{R}^{n}$, $t > 0$ and a~cube $Q^{n}$, we define the differences of order~$l$ as follows:
\begin{equation}
\label{2.5}
\Delta^{l}(h,\Omega)g(x):=
\begin{cases} \sum\limits_{j=0}^{l}C_{l}^{j}(-1)^{l+j}g(x+jh), & [x,x+hl] \subset \Omega  ,\\
0, & \mbox{ otherwise};
\end{cases}
\end{equation}
\begin{equation}
\label{2.6}
\overline{\Delta}^{l}_{r}(t,\Omega)g(x):=\left(\frac{1}{t^{n}}\int\limits_{tI^{n}}|\Delta^{l}(h,\Omega)g(x)|^{r}\,dh\right)^{\frac{1}{r}}, \qquad  x \in \mathbb{R}^{n};
\end{equation}
\begin{equation}
\label{2.7}
 \delta^{l}_{r}(Q^{n},\Omega)g:=
\left(\frac{1}{[r(Q^{n})]^{2n}}\int\limits_{r(Q^{n})I^{n}}\int\limits_{Q^{n}}
|\Delta^{l}(h,\Omega)g(x)|^{r}\,dxdh\right)^{\frac{1}{r}}.
\end{equation}

We set  $\Delta^{l}(h)g:=\Delta^{l}(h,\mathbb{R}^{n})g$, $\overline{\Delta}^{l}_{r}(t)g:=\overline{\Delta}^{l}_{r}(t,\mathbb{R}^{n})g$, $\delta^{l}_{r}(Q^{n})g:=\delta^{l}_{r}(Q^{n},\mathbb{R}^{n})g$.

For a cube $Q^{n}$ with $l \in \mathbb{N}$, $r \in (0,\infty]$ we let $\omega_{l}(\varphi,Q^{n})_{r}$
denote the modulus of continuity of a~function
$\varphi \in L^{\text{\rm loc}}_{r}(\mathbb{R}^{n})$ on a~cube $Q^{n}$ in the  $L_{r}(Q^{n})$-metric; that is,
$$
\omega_{l}(\varphi,Q^{n})_{r}:=\sup\limits_{|h|>0}\|\Delta^{l}(h,Q^{n})\varphi|L_{r}(\mathbb{R}^{n})\|.
$$

The following two-sided estimate is well known (for $r \geq 1$ see \cite{Br}; for the general setting see~\cite{Oswald}).
\begin{equation}
\begin{split}
\label{2.8}
&C_{1}\delta^{l}_{r}(Q^{n},Q^{n})\varphi \le
|Q^{n}|^{-\frac{1}{r}} \omega_{l}(\varphi, Q^{n})_{r} \le
C_{2}\delta^{l}_{r}(Q^{n},Q^{n})\varphi,
\end{split}
\end{equation}
the constants $C_{1}$, $C_{2}$ in \eqref{2.8}
being independent of both the function~$\varphi$ and the cube~$Q^{n}$.

Let $l \in \mathbb{N}$, $r \in (0,\infty]$. For a~cube $Q^{n}$ we define the local best approximation to a~function $\varphi \in
L^{\text{\rm loc}}_{r}(\mathbb{R}^{n})$ in the $L_{r}(Q^{n})$-metric by the polynomials of degree $<l$ as follows:
$$
E_{l}(\varphi,Q^{n})_{r}:=\inf\limits_{deg(P)<l}\|\varphi-P|L_{r}(Q^{n})\|.
$$

Next, for a cube  $Q^{n}$  we define  the local best approximation to a~function $\varphi \in
L^{\text{\rm loc}}_{r}(\mathbb{R}^{n})$ in the  $L_{r}(Q^{n})$-metric by the polynomials
of \textit{coordinate degree}  $<l$ (the total degree of a~polynomial is, clearly, at most $n(l-1)$) as
$$
\widehat{E}_{l}(\varphi,Q^{n})_{r}:=\inf\limits_{deg_{i}(P)<l}\|\varphi-P|L_{r}(Q^{n})\|,
$$
the infimum being taken over all polynomials $P$ whose degree in the variable $x_{i}$ is smaller than~$l$ for each $i \in \{1,\dots,n\}$.

From the results of \cite{Oswald} it follows that, for  $l \in \mathbb{N}$, $r \in (0,\infty]$,
\begin{equation}
\label{2.9}
C_{3}\delta^{l}_{r}(Q^{n},Q^{n})\varphi \le
|Q^{n}|^{-\frac{1}{r}} E_{l}(\varphi, Q^{n})_{r} \le
C_{4}\delta^{l}_{r}(Q^{n},Q^{n})\varphi.
\end{equation}
the constants $C_{3},C_{4}$ in~\eqref{2.9} are independent of both the function $\varphi$ and the cube~$Q^{n}$.

Let $l \in \mathbb{N}$, $r \in (0,\infty]$, $Q^{n}$ be a cube.
A polynomial $P_{Q^{n}}$ will be said to be is a~polynomial of almost best approximation
to a~function $\varphi \in L^{\text{\rm loc}}_{r}(\mathbb{R}^{n})$
by polynomials of degree $<l$ in the $L_{r}(Q^{n})$-metric with constant $A \geq 1$ if
$$
\|\varphi-P_{Q^{n}}|L_{r}(Q^{n})\| \le A E_{l}(\varphi,Q^{n})_{r}.
$$

The definition a polynomial of almost $L_{r}(Q^{n})$-best approximation by polynomials of \textit{coordinate degree} $<l$
to a~function $\varphi \in L^{\text{\rm loc}}_{r}(\mathbb{R}^{n})$ with constant $A\geq 1$ is similar.

\begin{Def}
\label{Def2.4}
Let $p \in (0,\infty]$. A weight sequence $\{t_{k}\}$ is called $p$-admissible if $t_{k} \in L^{\text{\rm loc}}_{p}(\mathbb{R}^{n})$ for all $k \in \mathbb{N}_{0}$.
\end{Def}

\begin{Def}
\label{Def2.5}
Let $l \in \mathbb{N}$, $0 < p,q,r \le \infty$, and let $\{t_{k}\}$ be a~$p$-admissible weight sequence. We set
\begin{equation}
\begin{split}
\label{2.10}
&\overline{B}_{p,q,r}^{l}(\mathbb{R}^{n},\{t_{k}\}):=\{\varphi:
\varphi \in L_{r}^{\text{\rm loc}}(\mathbb{R}^{n}),
\|\varphi|\overline{B}_{p,q,r}^{l}(\mathbb{R}^{n},\{t_{k}\})\|<+\infty\}
\mbox{, where}\\
&\left\|\varphi|\overline{B}^{l}_{p,q,r}(\mathbb{R}^{n},\{t_{k}\})\right\|:=\left[\sum\limits_{k=1}^{\infty}
\|t_{k}\overline{\Delta}^{l}_{r}(2^{-k})\varphi|L_{p}(\mathbb{R}^{n})\|^{q}\right]^{\frac{1}{q}}+\left(\int\limits_{\mathbb{R}^{n}}t^{p}_{0}(x)\|\varphi|L_{r}(x+I^{n})\|^{p}\,dx\right)^{\frac{1}{p}}
;
\end{split}
\end{equation}
\begin{equation}
\begin{split}
\label{2.11}
&\widetilde{B}_{p,q,r}^{l}(\mathbb{R}^{n},\{t_{k}\}):=\{\varphi:
\varphi \in L_{r}^{\text{\rm loc}}(\mathbb{R}^{n}),
\|\varphi|\widetilde{B}_{p,q,r}^{l}(\mathbb{R}^{n},\{t_{k}\})\|<+\infty\}
\mbox{, where}\\
&\left\|\varphi|\widetilde{B}^{l}_{p,q,r}(\mathbb{R}^{n},\{t_{k}\})\right\|:=\left[\sum\limits_{k=1}^{\infty}
\|t_{k}\delta^{l}_{r}(\cdot+2^{-k}I^{n})\varphi|L_{p}(\mathbb{R}^{n})\|^{q}\right]^{\frac{1}{q}}+\left(\int\limits_{\mathbb{R}^{n}}t^{p}_{0}(x)\|\varphi|L_{r}(x+I^{n})\|^{p}\,dx\right)^{\frac{1}{p}};
\end{split}
\end{equation}
the modifications for $p=\infty$ or $q=\infty$ are clear.

Let $\gamma$ $-$ weight, $l>s>0$ . We set $\widetilde{B}^{s}_{p,q,r}(\mathbb{R}^{n},\gamma):=\widetilde{B}^{l}_{p,q,r}(\mathbb{R}^{n},\{2^{ks}\gamma\})$, $\overline{B}^{s}_{p,q,r}(\mathbb{R}^{n},\gamma):=\overline{B}^{l}_{p,q,r}(\mathbb{R}^{n},\{2^{ks}\gamma\})$.
The corresponding spaces will be called weighted Besov spaces with weight~$\gamma$.
\end{Def}

\begin{Remark}
\label{R2.3}
The space
$\widetilde{B}_{p,q,r}^{l}(\mathbb{R}^{n},\{t_{k}\})$ ($\overline{B}_{p,q,r}^{l}(\mathbb{R}^{n},\{t_{k}\})$)
may prove to be trivial, containing only the functions that vanish almost everywhere
We put forward a~condition on the parameters $l,p,q$ and
a~$p$-admissible  sequence $\{t_{k}\}$ that guarantees that the corresponding space be nontrivial. Let  $p,q \in (0,\infty]$ and any cube
$Q^{n}$ (with corresponding modifications in the case $p,q=\infty$)
$$
\left(\sum\limits_{k=0}^{\infty}\left(\int\limits_{Q^{n}}2^{-klp}t^{p}_{k}(x)\,dx\right)^{\frac{q}{p}}\right)^{\frac{1}{q}}
< \infty.
$$

Under this condition the set $C_{0}^{\infty}
\subset \widetilde{B}_{p,q,r}^{l}(\mathbb{R}^{n},\{t_{k}\})$ ($\overline{B}_{p,q,r}^{l}(\mathbb{R}^{n},\{t_{k}\})$). This easily follows
from the expansion of~$\varphi$ in a~Taylor formula with remainder in the Lagrange form.
\end{Remark}

\begin{Remark}
\label{R2.4}
The space $\widetilde{B}_{p,q,1}^{l}(\mathbb{R}^{n},\{t_{k}\})$ was introduced in $\cite{Ty}$ for $p,q \in (1,\infty)$ for weight sequences $\{t_{k}\} \in ^{\text{\rm loc}}Y^{\alpha_{3}}_{\alpha_{1},\alpha_{2}}$ without restrictions on the parameters $\alpha_{1},\alpha_{2}$.
The space $\overline{B}_{p,q,1}^{l}(\mathbb{R}^{n},\{t_{k}\})$ was studied by H.~Kempka and J.~Vybiral \cite{KeVy} for
weight sequences $\{t_{k}\} \in Y^{\alpha_{3}}_{\alpha_{1},\alpha_{2}}$ for $p,q \in (0,\infty]$ under condition~\eqref{ogranicenie}.
The space close to the space $\overline{B}_{p,q,1}^{l}(\mathbb{R}^{n},\{t_{k}\})$ (but different from it!) was studied by Besov  \cite{Be4},\cite{Be2} with
$p,q \in (1,\infty)$, $\{t_{k}\} \in ^{\text{\rm loc}}Y^{\alpha_{3}}_{\alpha_{1},\alpha_{2}}$ under condition~\eqref{ogranicenie}.
\end{Remark}

\begin{Def}
\label{Def2.6}
Let $p \in (0,\infty]$. For a $p$-admissible  weight sequence $\{t_{k}\}$ we set (in the case $p=\infty$ we assume that $\frac{kn}{p}=0$)
\begin{equation}
\label{2.12}
t_{k,m}:=\|t_{k}|L_{p}(Q^{n}_{k,m})\| \mbox{ for } k \in \mathbb{N}_{0}, m \in \mathbb{Z}^{n},
\end{equation}
\begin{equation}
\label{2.13}
\overline{t}_{k}(x):= 2^{\frac{kn}{p}}\sum\limits_{m \in \mathbb{Z}^{n}}t_{k,m}\chi_{\widetilde{Q}^{n}_{k,m}}(x)   \mbox{ for } k \in \mathbb{N}_{0}, x \in \mathbb{R}^{n}.
\end{equation}

In what follows, the multiple sequence $\{t_{k,m}\}$ (the weight sequence $\{\overline{t}_{k}\}$) defined by formula \eqref{2.12} (\eqref{2.13}) will be called
multiple sequence (weight sequence) $p$-associated with the weight sequence $\{t_{k}\}$.
\end{Def}

\begin{Def}
\label{Def2.7}
 Let $\alpha_{3} \geq 0$, $\alpha_{1},\alpha_{2} \in \mathbb{R}$, $\sigma_{1},\sigma_{2} \in (0,+\infty]$, $\alpha=(\alpha_{1},\alpha_{2})$, $\sigma=(\sigma_{1},\sigma_{2})$.
 By  $X^{\alpha_{3}}_{\alpha,\sigma,p}=X^{\alpha_{3}}_{\alpha,\sigma,p}(\mathbb{R}^{n})$ we will denote the set of $p$-admissible weight sequences $\{t_{k}\}$,
satisfying the following conditions:

1) There exist numbers $c_{1},c_{2} > 0$  such that
\begin{equation}
\begin{split}
\label{2.14}
\left(2^{kn}\int\limits_{Q^{n}_{k,m}}\overline{t}^{p}_{k}(x)\right)^{\frac{1}{p}}\left(2^{kn}\int\limits_{Q^{n}_{k,m}}(\overline{t}_{j})^{-\sigma_{1}}(x)\right)^{\frac{1}{\sigma_{1}}} \le C_{1}2^{\alpha_{1}(k-j)}, \qquad
 0 \le k \le j, m \in \mathbb{Z}^{n},
\end{split}
\end{equation}
\begin{equation}
\begin{split}
\label{2.15}
\left(2^{kn}\int\limits_{Q^{n}_{k,m}}\overline{t}^{p}_{k}(x)\right)^{-\frac{1}{p}}\left(2^{kn}\int\limits_{Q^{n}_{k,m}}\overline{t}^{\sigma_{2}}_{j}(x)\right)^{\frac{1}{\sigma_{2}}} \le C_{2}2^{\alpha_{2}(j-k)}, \qquad
 0 \le k \le j, m \in \mathbb{Z}^{n},
\end{split}
\end{equation}
(the modifications of \eqref{2.6} and  \eqref{2.7} for $\sigma_{1}=\infty$ and $\sigma_{2}=\infty$ are clear).

2) For all $k \in \mathbb{N}_{0}$
\begin{equation}
  \label{2.16}
  0<t_{k,m} \le 2^{\alpha_{3}} t_{k,\widetilde{m}}, \ \ \mbox{for } m,\widetilde{m} \in \mathbb{Z}^{n},  |m_{i}-\widetilde{m}_{i}| \le 1, i=1,..,n,
   \end{equation}
\end{Def}

\begin{Remark}
\label{R2.5}
We denote by $\widetilde{X}^{\alpha_{3}}_{\alpha,\sigma,p}$ subset of  $X^{\alpha_{3}}_{\alpha,\sigma,p}$ consisting of only $p$-admissible
weight sequences $\{t_{k}\}=\{\overline{t}_{k}\}$.  It is clear that
$\widetilde{X}^{\alpha_{3}}_{\alpha,\sigma,p}=^{\text{\rm loc}}Y^{\alpha_{3}}_{\alpha_{1},\alpha_{2}}$ for $p \in (0,\infty]$, $-\infty < \alpha_{1} \le \alpha_{2} < \infty$, $\alpha_{3} \geq 0$ and $\sigma_{1}=\sigma_{2}=\infty$.

Given fixed $-\infty < \alpha^{i}_{1} \le \alpha^{i}_{2} < \infty$, $\alpha_{3} \geq 0$, $\sigma^{i}_{1},\sigma^{i}_{2} \in (0,\infty]$, $i=1,2$,
 we set $\alpha^{i}:=(\alpha^{i}_{1},\alpha^{i}_{2})$, $\sigma^{i}=(\sigma^{i}_{1},\sigma^{i}_{2})$ for $i=1,2$.
Elementary arguments based on H\"older's inequality and the monotonicity of the $l_{q}$-norm  (in~$q$) prove the embedding
$\widetilde{X}^{\alpha_{3}}_{\alpha^{1},\sigma^{1},p} \subset \widetilde{X}^{\alpha_{3}}_{\alpha^{2},\sigma^{2},p}$, provided that $\alpha^{2}_{1}=\alpha^{1}_{1}+n\min\{\frac{1}{\sigma^{2}_{1}}-\frac{1}{\sigma^{1}_{1}},0\}$, $\alpha^{2}_{2}=\alpha^{1}_{2}+n\max\{\frac{1}{\sigma^{1}_{2}}-\frac{1}{\sigma^{2}_{2}},0\}$.
\end{Remark}

\begin{Remark}
\label{R2.6}
Clearly, it may happen that a multiple sequence $\{t_{k,m}\}$ is $p$-associated with several weight sequences. However, this will not be an impediment
for further constructions if  $\{t_{k}\} \in X^{\alpha_{3}}_{\alpha,\sigma,p}$ . Indeed, using  \eqref{2.16} and elementary arguments we have
\begin{equation}
\begin{split}
\label{2.19}
&\left\|\varphi|\widetilde{B}^{l}_{p,q,r}(\mathbb{R}^{n},\{t_{k}\})\right\| \sim \left\|\varphi|\widetilde{B}^{l}_{p,q,r}(\mathbb{R}^{n},\{\overline{t}_{k}\})\right\| \sim \left\|\varphi|\widetilde{B}^{l}_{p,q,r}(\mathbb{R}^{n},\{t_{k,m}\})\right\|^{(1)}:=\\
&=\left\|\left(\sum\limits_{m
\in \mathbb{Z}^{n}}
t^{p}_{k,m}\left[\delta^{l}_{r}(Q_{k,m}^{n})\varphi\right]^{p}\right)^{\frac{1}{p}}\left|\right.l_{q}\right\|+\left(\sum\limits_{m \in
\mathbb{Z}^{n}}t^{p}_{0,m}\|\varphi|L_{r}(Q_{0,m}^{n})\|^{p}\right)^{\frac{1}{p}}
\end{split}
\end{equation}
(with corresponding modifications in the case $p=\infty$).

Here, the constant through which one norm is estimated in terms of the other one in \eqref{2.19} will depend only on  $\alpha_{3},l,p,n$.

For a fixed $p \in (0,\infty]$, it is clear that there exists a bijection between the multiple sequences $\{t_{k,m}\}$ and the sets of weight sequences $\{t_{k}\} \in X^{\alpha_{3}}_{\alpha,\sigma,p}$,
for which the multiple sequence $\{t_{k,m}\}$ is $p$-associated with the weight sequence $\{t_{k}\}$.

Considering \eqref{2.19}, in what follows the space $\widetilde{B}^{l}_{p,q,r}(\mathbb{R}^{n},\{t_{k}\})$ will also be denoted by the symbol
$\widetilde{B}^{l}_{p,q,r}(\mathbb{R}^{n},\{t_{k,m}\})$.
\end{Remark}

\begin{Def}
\label{Def2.8}
Let $p \in (0,\infty)$, $d \in \mathbb{N}_{0}$, and let a weight $\gamma^{p} \in A^{\text{\rm loc}}_{\infty}(\mathbb{R}^{n+d})$.
We set $\Xi^{d,n}_{k,m}:=Q^{n}_{k,m}\times \left(2^{-k}B^{d} \setminus 2^{-k-1}B^{d}\right)$ for $k \in \mathbb{N}_{0}$, $m \in Z^{n}$.
The multiple sequence $\widehat{\gamma}_{k,m}$ defined by
$$
\widehat{\gamma}_{k,m}:=\|\gamma|L_{p}(\Xi^{d,n}_{k,m})\| \mbox{ for } k \in \mathbb{N}_{0}, m \in \mathbb{Z}^{n}
$$
will be called the multiple sequence generated by the weight $\gamma$.
\end{Def}

The following important properties of the sequence $\{\widehat{\gamma}_{k,m}\}$ will be required in what follows.

\begin{Lm}
\label{Lm2.2}
Let $p \in (0,\infty)$, $d \in \mathbb{N}_{0}$, a weight $\gamma^{p} \in A^{\text{\rm loc}}_{\infty}(\mathbb{R}^{n+d})$,
and let the multiple sequence $\widehat{\gamma}_{k,m}$ be generated by the weight~$\gamma$.
Also let ${m} \in \mathbb{Z}^{n}$, $k,j \in \mathbb{N}_{0}$, $j \geq k$, $G_{j,k,m}$ be an arbitrary set of cubes $Q^{n}_{j,\widetilde{m}} \subset Q^{n}_{k,m}$. Then

{\rm 1)} the inequality holds
\begin{equation}
\label{2.20}
\sum\limits_{Q^{n}_{j,\widetilde{m}} \subset Q^{n}_{k,m}}\widehat{\gamma}^{p}_{j,\widetilde{m}} \le C 2^{(k-j)d \delta}\widehat{\gamma}^{p}_{k,m},
\end{equation}
in which the constant $C>0$ and the number $\delta(\gamma)>0$ is independent of $k,j,\widetilde{m}$;

{\rm 2)} the inequality holds
\begin{equation}
\label{2.21}
\sum\limits_{Q^{n}_{j,\widetilde{m}} \in G_{j,k,m}}\widehat{\gamma}^{p}_{j,\widetilde{m}} \le C \left(\frac{|\bigcup\limits_{Q^{n}_{j,\widetilde{m}} \in G_{j,k,m}}Q^{n}_{j,\widetilde{m}}|}{|Q^{n}_{k,m}|}\right)^{\delta'}\sum\limits_{Q^{n}_{j,\widetilde{m}} \subset Q^{n}_{k,m}}\widehat{\gamma}^{p}_{j,\widetilde{m}},
\end{equation}
in which the constants $C>0$ and $\delta'>0$ dependent only on $\gamma,n,d$;

{\rm 3)} for $a \geq 1$, $|m-\widetilde{m}| \le a$, $k \geq 0$,
$$
 2^{-\delta_{3}}\widehat{\gamma}_{k,m} \le \widehat{\gamma}_{k,\widetilde{m}} \le 2^{\delta_{3}} \widehat{\gamma}_{k,m},
$$
where the number $\delta_{3}(\gamma) \geq 0$ depends  only on $\gamma,n,p,d$;

{\rm 4)} for any cube $Q^{n}_{k,m}$ and any cube
$Q^{n}_{k+1,\widetilde{m}} \subset Q^{n}_{k,m}$,
$$
\widehat{\gamma}_{k,m} \le C \widehat{\gamma}_{k+1,\widetilde{m}},
$$
 for $k \geq 0$, $m \in \mathbb{Z}^{n}$, where the constant $C>0$ depends only on $\gamma,n,d,p$.
\end{Lm}

\textbf{Proof.}
To prove 3) it suffices to take some cube $Q^{n+d}$ containing both sets $\Xi^{d,n}_{k,m}$ and $\Xi^{d,n}_{k,\widetilde{m}}$ and
use the fact that $\gamma^{p}$ satisfies the doubling condition on the cube $Q^{n+d}$ with doubling constant depending only on the constant
$C^{\text{\rm loc}}_{\gamma,p,r(Q^{n+d})}$ (the proof of the last fact is similar to the proof of the corresponding result in \cite{St}, Ch.~5).
The proof of property~4) is similar to that of property~3).

Let us prove property 1); property 2) is dealt with similarly. It is easily seen that
\begin{equation}
\label{2.22}
\frac{|\bigcup\limits_{Q^{n}_{j,\widetilde{m}}\subset Q^{n}_{k,m}}\Xi^{d,n}_{j,\widetilde{m}}|}{| \Xi^{d,n}_{k,m}|} \le C 2^{(k-j)d}.
\end{equation}

Using Definition \ref{Def2.3} and Remark~\ref{R2.1} one may easily prove that for some $\delta(\gamma) > 0$, for any cube $Q^{n}$, $r(Q^{n}) \le a$, and any
measurable set $F \subset Q^{n}$,
\begin{equation}
\label{2.23}
\frac{\gamma^{p}(F)}{\gamma^{p}(Q^{n})} \le C \left(\frac{|F|}{|Q^{n}|}\right)^{\delta(\gamma)},
\end{equation}
the constant $C > 0$ being independent of both the cube $Q^{n}$ and the set~$F$.

From \eqref{2.22}, \eqref{2.23} we get \eqref{2.20}, completing the proof of the lemma.

We let $\delta_{1}(\gamma):=\delta_{1}(\gamma,n,d)$ denote the supremum over all $\delta$ for which \eqref{2.20} holds.
Similarly, $\delta_{2}(\gamma):=\delta_{2}(\gamma,n,d)$ will denote the supremum over all $\delta'$ satisfying~\eqref{2.11}.

Note that in general $\delta_{1} \neq \delta_{2}$. Indeed, let $\gamma^{p}(x_{1},x_{2}):=x_{1}^{\beta}$ with $(x_{1},x_{2}) \in \mathbb{R}^{2}$, $\beta > 0$.
Clearly, $\gamma^{p} \in A_{\infty}(\mathbb{R}^{2})$. Also, $\delta_{1}(\gamma)=\frac{1}{2}$ for any $\beta > 0$, whereas $\delta_{2}(\gamma)$ depends on
$\beta > 0$.

\begin{Example}
\label{Ex2.1}
For future purposes we give an important example of a~weight sequence $\{t_{k}\} \in \widetilde{X}^{\alpha_{3}}_{\alpha,\sigma,p}$.
We note that this example is the main impetus for practical applications of the classes $X^{\alpha_{3}}_{\alpha,\sigma,p}$.

Let $d \in \mathbb{N}_{0}$, $p \in  (0,\infty)$, a weight $\gamma^{p} \in A^{\text{\rm loc}}_{\infty}(\mathbb{R}^{n+d})$, and let a multiple sequence $\{\widehat{\gamma}_{k,m}\}$ be generated
 by the weight $\gamma$. By Remark~\ref{R2.1} we have $\gamma^{p} \in A^{\text{\rm loc}}_{p_{0}}(\mathbb{R}^{n})$ for some $p_{0} \in [1,\infty)$.
 Assume that a~weight sequence $\{s_{k}\} \in ^{\text{\rm loc}}Y^{\alpha'_{3}}_{\alpha'_{1},\alpha'_{2}}$. We set  $\{s_{k,m}\} = \|s_{k}|L_{p}(Q^{n}_{k,m})\| $, $t_{k,m}:=\widehat{\gamma}_{k,m}(2^{\frac{kn}{p}}s_{k,m})$ for $k \in \mathbb{N}_{0}$, $m \in \mathbb{Z}^{n}$. Then the weight sequence $\{t_{k}\} = \{\overline{t}_{k}\} \in \widetilde{X}^{\alpha_{3}}_{\alpha,\sigma,p}$ for $\alpha_{3}=\alpha'_{3}+\delta_{3}(\gamma)$, $\alpha_{2}=\alpha'_{2}-\frac{d(\delta_{1}(\gamma)-\varepsilon)}{p}$, $\alpha_{1}=\alpha'_{1}+\frac{n}{\sigma_{1}}+\frac{n}{p}-\frac{(n+d)p_{0}}{p}+\frac{dp_{0}}{pp'_{0}}\left(\delta_{1}(\gamma^{-\frac{pp'_{0}}{p_{0}}},n,d)-\varepsilon\right)$, $\sigma_{2}=p$,
 and $\sigma_{1} = p\frac{p_{0}'}{p_{0}}$ for any $\varepsilon > 0$. Indeed,   \eqref{2.15} and \eqref{2.16} easily follow from assertions 1) and 3) of Lemma~\ref{Lm2.1}.
 Let us verify \eqref{2.14} with $p_{0} > 1$, the case $p_{0}=1$ is dealt with similarly. By Definition~\ref{Def2.2},
$$
\left(2^{kn}\int\limits_{Q^{n}_{k,m}}\overline{t}^{p}_{k}(x)\right)^{\frac{1}{p}}\left(2^{kn}\int\limits_{Q^{n}_{k,m}}(\overline{t}_{j})^{-\sigma_{1}}(x)\right)^{\frac{1}{\sigma_{1}}} \le
 C_{1} 2^{(k-j)(\alpha'_{1}+\frac{n}{\sigma_{1}}+\frac{n}{p})}\widehat{\gamma}_{k,m}\left(\sum_{\substack{\widetilde{m} \in \mathbb{Z}^{n}\\
Q^{n}_{j,\widetilde{m}} \subset Q^{n}_{k,m}}}\frac{1}{(\widehat{\gamma}_{j,\widetilde{m}})^{p\frac{p_{0}'}{p_{0}}}}\right)^{\frac{p_{0}}{pp'_{0}}} \le
$$
$$
\le C_{2} 2^{(k-j)(\alpha'_{1}+\frac{n}{\sigma_{1}}+\frac{n}{p})}\widehat{\gamma}_{k,m}2^{j(n+d)\frac{p_{0}}{p}}\left(\sum_{\substack{\widetilde{m} \in \mathbb{Z}^{n}\\
Q^{n}_{j,\widetilde{m}} \subset Q^{n}_{k,m}}}\int\limits_{\Xi^{d,n}_{j,\widetilde{m}}}\gamma^{-p\frac{p'_{0}}{p_{0}}}(x)\,dx\right)^{\frac{p_{0}}{pp'_{0}}} \le
$$
$$
 \le C_{3} 2^{(k-j)(\alpha'_{1}+\frac{n}{\sigma_{1}}+\frac{n}{p}+\frac{dp_{0}}{pp'_{0}}(\delta_{1}(\gamma^{-\frac{pp'_{0}}{p_{0}}},n,d)-\varepsilon))}\widehat{\gamma}_{k,m}2^{j(n+d)\frac{p_{0}}{p}}\left(\int\limits_{\Xi^{d,n}_{k,m}}\gamma^{-p\frac{p'_{0}}{p_{0}}}(x)\,dx\right)^{\frac{p_{0}}{pp'_{0}}} \le C_{4} 2^{(k-j)\alpha_{1}}.
$$

It is worth pointing out that $\alpha_{i}=\alpha'_{i}$ ($i=1,2$) in the case $d=0$.
\end{Example}

Let $p,q\in (0,\infty]$, $r \in (0,p]$, $\alpha_{1},\alpha_{2} \in \mathbb{R}$, $\alpha_{3} \geq 0$, $\sigma_{1},\sigma_{2} \in (0,\infty]$, and let $\{t_{k,m}\}$ be
the $p$-associated multiple sequence with a~$p$-admissible  weight sequence $\{t_{k}\} \in X^{\alpha_{3}}_{\alpha,\sigma,p}$, $c > 1$.
In the space $\widetilde{B}^{l}_{p,q,r}(\mathbb{R}^{n},\{t_{k}\})$, we consider
the quasi-norms generated by the multiple sequence $\{t_{k,m}\}$:
\begin{equation}
\begin{split}
\label{2.24}
&\|\varphi|\widetilde{B}^{l}_{p,q,r}(\mathbb{R}^{n},\{t_{k,m}\},c)\|^{(2)}:=\left\|\left(\sum\limits_{m
\in
\mathbb{Z}^{n}}t^{p}_{k,m}\left[\delta^{l}_{r}(c
Q^{n}_{k,m},c Q^{n}_{k,m})\varphi\right]^{p}\right)^{\frac{1}{p}}|l_{q}\right\|+\left(\sum\limits_{m \in
\mathbb{Z}^{n}}t^{p}_{0,m}\|\varphi|L_{r}(Q_{0,m}^{n})\|^{p}\right)^{\frac{1}{p}},\\
\end{split}
\end{equation}
\begin{equation}
\begin{split}
\label{2.25}
&\|\varphi|\widetilde{B}^{l}_{p,q,r}(\mathbb{R}^{n},\{t_{k,m}\},c)\|^{(3)}:=\left\|\left(\sum\limits_{m
\in \mathbb{Z}^{n}}t^{p}_{k,m}\left[2^{\frac{kn}{r}} E_{l}(\varphi,c
Q^{n}_{k,m})_{r}\right]^{p}\right)^{\frac{1}{p}}|l_{q}\right\|+\left(\sum\limits_{m \in
\mathbb{Z}^{n}}t^{p}_{0,m}\|\varphi|L_{r}(Q_{0,m}^{n})\|^{p}\right)^{\frac{1}{p}},
\end{split}
\end{equation}
\begin{equation}
\begin{split}
\label{2.26}
&\|\varphi|\widetilde{B}^{l}_{p,q,r}(\mathbb{R}^{n},\{t_{k,m}\},c)\|^{(4)}:=\left\|\left(\sum\limits_{m
\in \mathbb{Z}^{d}}t^{p}_{k,m}\left[2^{\frac{kn}{r}}\omega_{l}(\varphi,c
Q^{n}_{k,m})_{r}\right]^{p}\right)^{\frac{1}{p}}|l_{q}\right\|+\left(\sum\limits_{m \in
\mathbb{Z}^{n}}t^{p}_{0,m}\|\varphi|L_{r}(Q_{0,m}^{n})\|^{p}\right)^{\frac{1}{p}}.
\end{split}
\end{equation}

\begin{Th}
\label{Th2.3}
Let $p,q,r \in (0,\infty]$,  $\alpha_{3} \geq 0$, $\alpha_{1},\alpha_{2} \in \mathbb{R}$, $\sigma_{1},\sigma_{2} \in (0,+\infty]$,
 $\{t_{k}\} \in X^{\alpha_{3}}_{\alpha,\sigma,p}$ be a~$p$-admissible sequence, and $\{t_{k,m}\}$ be the associated multiple sequence. Then, for $i=1,2,3,4$, the quasi-norms
$\|\cdot|\widetilde{B}^{l}_{p,q,r}(\mathbb{R}^{n},\{t_{k,m}\},c)\|^{(i)}$
are equivalent on the space $\widetilde{B}^{l}_{p,q,r}(\mathbb{R}^{n},\{t_{k}\})$.
\end{Th}

\textbf{Proof}. This theorem was proved in \cite{Ty} with $r=1$, $p,q \in (1,\infty)$, $\{t_{k}\} \in ^{\text{\rm loc}}Y^{\alpha_{3}}_{\alpha_{1},\alpha_{2}}$.
In the general setting the proof is similar if we take into account \eqref{2.8}, \eqref{2.9}, \eqref{2.19}, Remark \eqref{R2.5} and use the estimate
\begin{equation}
\label{2.1'}
\delta^{l}_{r}(c Q^{n}_{k,m},c Q^{n}_{k,m})\varphi \le C \sum_{\substack{\widetilde{m} \in \mathbb{Z}^{n}\\
Q^{n}_{k-j(c),\widetilde{m}}\bigcap cQ^{n}_{k,m} \neq
\emptyset}}\delta^{l}_{r}(Q^{n}_{k-j(c),\widetilde{m}})\varphi, \hbox{ where }
\end{equation}
where $c > 1$ is from the hypotheses of the theorem, and  $j(c) \geq 1$ is the smallest natural number such that $2^{k-j(c)} > c2^{k}$. The constant $C$ in \eqref{2.1'} depends only on $n,r,l,c$.

We set  $\widetilde{B}^{l}_{p,q}(\mathbb{R}^{n},\{t_{k}\}):=\widetilde{B}^{l}_{p,q,1}(\mathbb{R}^{n},\{t_{k}\})$, $\widetilde{B}^{l}_{p}(\mathbb{R}^{n},\{t_{k}\}):=\widetilde{B}^{l}_{p,p,1}(\mathbb{R}^{n},\{t_{k}\})$.

\begin{Th}
\label{Th2.4}
Let $p,q,r \in (0,\infty]$,  $l \in \mathbb{N}$, $\alpha_{3} \geq 0$, $\alpha_{1},\alpha_{2} \in \mathbb{R}$, $\sigma_{1},\sigma_{2} \in (0,+\infty]$,
$\{t_{k}\} \in X^{\alpha_{3}}_{\alpha,\sigma,p}$ be a~$p$-admissible  weight sequence. Then the space $\widetilde{B}_{p,q,r}^{l}(\mathbb{R}^{n},\{t_{k}\})$ is complete.
\end{Th}

The proof, which is close in spirit to that of the completeness of the classical Besov space \cite{Be1}, depends on the completeness of the space
$L_{r}(Q^{n})$ for any cube $Q^{n}$, uses Remark~\ref{R2.5} and the equivalent norm  \eqref{2.19}. We suppress the details, which are quite standard.

\begin{Lm}
\label{Lm2.5}
Let $p,q \in (0,\infty]$, $r \in (0,p]$, $l \in \mathbb{N}$, $\alpha_{3} \geq 0$, $\alpha_{1},\alpha_{2} \in \mathbb{R}$, $\sigma_{1},\sigma_{2} \in (0,+\infty]$,
$\{t_{k}\} \in \widetilde{X}^{\alpha_{3}}_{\alpha,\sigma,p}$ be a~$p$-admissible  weight sequence. Then
$\overline{B}_{p,q,r}^{l}(\mathbb{R}^{n},\{t_{k}\}) \subset \widetilde{B}_{p,q,r}^{l}(\mathbb{R}^{n},\{t_{k}\})$.
\end{Lm}

\textbf{Proof}. We shall consider only the case $p,q \in (0,\infty)$, the arguments in the case $p=\infty$ or $q=\infty$ are similar.
We compare the first terms in norms \eqref{2.14} and \eqref{2.15}. Applying H\"older's inequality to the integral in~$y$ and using~\eqref{2.20},
we have, for $k \in \mathbb{N}$,
\begin{equation}
\begin{split}
\label{2.28}
&\sum\limits_{m \in
\mathbb{Z}^{n}}t^{p}_{k,m}\left[\delta^{l}_{r}(Q_{k,m}^{n})\varphi\right]^{p}
\le \sum\limits_{m \in
\mathbb{Z}^{n}}t^{p}_{k,m}2^{\frac{npk}{r}+nk}\int\limits_{Q^{n}_{k,m}}\left[\int\limits_{\frac{1}{2^{k}}I^{n}}|\Delta^{l}(h)\varphi(y)|^{r}\,dh\right]^{\frac{p}{r}}\,dy
=\\
&=\int\limits_{\mathbb{R}^{n}}t^{p}_{k}(y)[\overline{\Delta}^{l}_{r}(2^{-k})\varphi (y)]^{p}\,dy.
\end{split}
\end{equation}

The proof of the embedding will be completed if we raise the both parts of \eqref{2.28} to the power $\frac{q}{p}$ and sum over all~$k$.

\begin{Th}
\label{Th2.5}
Let $p,q \in (0,\infty]$,  $p \neq \infty$ , $p_{0} \in [1,\infty)$, $0 < r_{1} \le r_{2} \le \frac{p}{p_{0}}$,  $l \in \mathbb{N}$, $\gamma^{p}
\in A^{\text{\rm loc}}_{p_{0}}(\mathbb{R}^{n})$. Let $\alpha_{3} \geq 0$, $0< \alpha_{1} \le \alpha_{2} < l$ and let a weight sequence
$\{s_{k}\} \in ^{\text{\rm loc}}Y^{\alpha_{3}}_{\alpha_{1},\alpha_{2}}$. Also let $t_{k}(x)=\gamma(x)s_{k}(x)$ for $k \in \mathbb{N}_{0}$, $x \in \mathbb{R}^{n}$. Then
$\overline{B}_{p,q,r_{1}}^{l}(\mathbb{R}^{n},\{t_{k}\}) = \widetilde{B}_{p,q,r_{2}}^{l}(\mathbb{R}^{n},\{t_{k}\})$, the corresponding norms being equivalent.
\end{Th}

The proof of Theorem~\ref{Th2.5} depends on the atomic decomposition theorem for the space $\widetilde{B}_{p,q,r_{2}}^{l}(\mathbb{R}^{n},\{t_{k}\})$, and so we defer it to the end of~\S\,4.

\begin{Remark}
\label{R2.7}
For $\gamma \equiv 1$ the conclusion of Theorem~\ref{Th2.5} may be extended also to the case $p=\infty$.
\end{Remark}

The following result, which was proved in~\cite{KeVy}, will be given in a~simplified form with constant $p$ and~$q$.

\begin{Th}
\label{Th2.6}
Let $p,q \in (0,\infty]$, $\alpha_{1}> n\left(\frac{1}{\min\{p,1\}}-1\right)\left[1+\frac{\alpha_{3}}{n}p\right]$, $l > \alpha_{2}$, $\{s_{k}\} \in Y^{\alpha_{3}}_{\alpha_{1},\alpha_{2}}$. Then $B^{\{s_{k}\}}_{p,q}(\mathbb{R}^{n})=\overline{B}_{p,q,1}^{l}(\mathbb{R}^{n},\{s_{k}\})$, the corresponding quasi-norms being equivalent.
\end{Th}

Combining Theorems~\ref{Th2.5}, \ref{2.6} and Remark \ref{R2.7} we obtain

\begin{Ca}
\label{Th2.7}
Let $p,q \in [1,\infty]$, $r_{1},r_{2} \in [1,p]$, $\alpha_{3} \geq 0$, $\alpha_{1} > 0$, $l > \alpha_{2}$, $\{s_{k}\} \in Y^{\alpha_{3}}_{\alpha_{1},\alpha_{2}}$. Then $B^{\{s_{k}\}}_{p,q}(\mathbb{R}^{n})=\overline{B}_{p,q,r_{1}}^{l}(\mathbb{R}^{n},\{s_{k}\}) = \widetilde{B}_{p,q,r_{2}}^{l}(\mathbb{R}^{n},\{s_{k}\})$,
the corresponding norms being equivalent.
\end{Ca}

\begin{Remark}
\label{R2.8}
The question of the coincidence (or noncoincidence) of the spaces $\widetilde{B}_{p,q,r}^{l}(\mathbb{R}^{n},\{s_{k}\})$, $\overline{B}_{p,q,r}^{l}(\mathbb{R}^{n},\{s_{k}\})$ and $B^{\{s_{k}\}}_{p,q}(\mathbb{R}^{n})$ for weaker (in comparison with (Theorems  \ref{Th2.5}, \ref{Th2.6} or~\ref{Th2.7})) constraints on the variable smoothness $\{s_{k}\}$
is a~matter for the future.
\end{Remark}

\begin{Remark}
\label{R2.9}
Combining Theorem \ref{Th2.5} with Theorem 3.14 of \cite{HN} with $p \in (0,\infty)$, $0 < r \le p$, $q \in (0,\infty]$, $s>0$, $l > s$,  $\gamma^{p} \in A_{\frac{p}{r}}(\mathbb{R}^{n})$
we obtain $\overline{B}_{p,q,r}^{s}(\mathbb{R}^{n},\gamma) = \widetilde{B}_{p,q,r}^{s}(\mathbb{R}^{n},\gamma)=B^{s}_{p,q}(\mathbb{R}^{n}, \gamma)$, the corresponding norms being equivalent
\end{Remark}

\section{Trace space of weighted Sobolev space}

As was pointed out in the introduction, the main impetus for the study of the spaces $\widetilde{B}_{p,q,r}^{l}(\mathbb{R}^{n},\{t_{k}\})$ stems from their application
in the problem of traces of weighted Sobolev spaces.

For a brief overview of the available literature on traces of weighted Sobolev spaces we refer to~\cite{Ty}; we do not dwell on this here.

For a brief overview of the available literature on traces of weighted Sobolev spaces we refer to~\cite{Ty}; we do not dwell on this here.

Let $p \in [1,\infty]$, $l \in \mathbb{N}$, $\gamma$ be a~weight. We fix $n,d \in \mathbb{N}$. A~point of the $(n+d)$-dimensional Euclidean space
$\mathbb{R}^{n+d}:=\mathbb{R}^{n} \times \mathbb{R}^{d}$ will be written as a~pair $(x,y)$.
The plane given in~$\mathbb{R}^{n+d}$ by the equation $y=0$ will be identified with the space $\mathbb{R}^{n}$. For $a > 0$,
we set $^{n}R^{d}_{a}:=\mathbb{R}^{n+d} \setminus (\mathbb{R}^{n} \times aB^{d})$ and put
$\Xi^{d,n}_{k,m}:=Q^{n}_{k,m}\times(\frac{B^{d}}{2^{k}}\setminus \frac{B^{d}}{2^{k+1}})$ for $k \in \mathbb{N}_{0}$, $m \in \mathbb{Z}^{n}$.

By $W^{l}_{p}(\mathbb{R}^{n+d},\gamma)$ we shall denote the linear space of classes of equivalent functions having on~$\mathbb{R}^{n}$ all
(Sobolev) generalized derivatives up to order~$l$ inclusively. We equip this space with the norm
$$
\|f|W^{l}_{p}(\mathbb{R}^{n+d},\gamma)\|=\sum\limits_{|\alpha| \le l}\|D^{\alpha}f|L_{p}(\mathbb{R}^{n+d}, \gamma)\|.
$$

In \cite{Ty} the trace problem was solved for the spaces $\widetilde{W}_{p}^{l}(\mathbb{R}^{n},\gamma)$,
which slightly differ from the spaces $W^{l}_{p}(\mathbb{R}^{n},\gamma)$ in terms of the norm form. More precisely,
the norm of the space $\widetilde{W}^{l}_{p}(\mathbb{R}^{n},\gamma)$ does not include some mixed derivatives.

In what follows we shall require a~certain averaging operator, which was constructed in~\cite{Ty}.
We shall not give the details of the construction of this operator. For a~function  $\varphi \in L^{\text{\rm loc}}_{1}(\mathbb{R}^{n})$ we set
\begin{equation}
\label{3.1}
E_{\varepsilon}[\varphi](x):=\frac{1}{\varepsilon^{2n}}\sum\limits_{j=1}^{l}\mu_{j}\int\limits_{\mathbb{R}^{n}}\Theta\left(\frac{y-x}{\varepsilon}\right)\int\limits_{\mathbb{R}^{n}}\Theta\left(\frac{z-y}{j\varepsilon}\right)\varphi(z)\,dzdy, \mbox{ for } x \in \mathbb{R}^{n}.
\end{equation}

In \eqref{3.1} function $\Theta \in C^{\infty}_{0}$ is chosen appropriately,
$\int\limits_{\mathbb{R}^{n}}\Theta(x)\,dx=1$, and $\mu_{j}$ are specially chosen constants (see~\cite{Ty}).
Given $k \in \mathbb{N}_{0}$ we define $E_{k}[g]:=E_{2^{-k}}[g]$.

\begin{Lm}
\label{Lm3.1}
Let a function $\varphi \in L^{\text{\rm loc}}_{1}(\mathbb{R}^{n})$. Then, for any number $\varepsilon>0$ and a~multi-index
$\alpha, |\alpha|=l$ for $x \in \mathbb{R}^{n}$
\begin{equation}
\label{3.2} \left|D^{\alpha}_{x}E_{\varepsilon}[\varphi](x)\right| \le
\frac{1}{\varepsilon^{l}}\delta^{l}(x+\varepsilon I^{n})\varphi.
\end{equation}

Moreover, for any numbers $0 < \varepsilon_{1} < \varepsilon_{2}$, a multi-index ~$\beta$, and $x \in \mathbb{R}^{n}$,
\begin{equation}
\label{3.3}
\left|D^{\beta}E_{\varepsilon_{1}}[\varphi](x)-D^{\beta}E_{\varepsilon_{2}}[\varphi](x)\right|
\le C
\int\limits_{\varepsilon_{1}}^{\varepsilon_{2}}\frac{1}{t^{1+|\beta|}}\delta^{l}(
x+t I^{n})\varphi\,dt.
\end{equation}
\end{Lm}

\textbf{Proof}. For $\beta=0$ the proof is given in~\cite{Ty}. The general case is dealt with similarly.

\smallskip

In this section we are not going to give a~precise definition of the trace of a~function
$f \in L^{\text{\rm loc}}_{1}(\mathbb{R}^{n+d})$ on the plane $y=0$. This is a~standard definition and may be found, for example,
in Chapter~5 of the book~\cite{10}.

Assume that a multiple sequence $\{\widehat{\gamma}_{k,m}\}$ is generated by a weight $\gamma \in A_{\infty}^{\text{\rm loc}}(\mathbb{R}^{n+d})$.
Next, assume that parameters $l \in \mathbb{N}$ and $p \in (1,\infty)$ are fixed.
We set $\gamma^{l}_{k}(x):=\gamma_{k}(x):=2^{k(l+\frac{n}{p})}\sum\limits_{m \in \mathbb{Z}^{n}}\chi_{\widetilde{Q}^{n}_{k,m}}(x)\widehat{\gamma}_{k,m}$ for $k \in \mathbb{N}_{0}$, $x \in \mathbb{R}^{n}$, $\gamma_{k,m}:=2^{kl}\widehat{\gamma}_{k,m}$ for $k \in \mathbb{N}_{0}, m \in \mathbb{Z}^{n}$.

The main result of the present section is the following

\begin{Th}
\label{Th3.1}
Let $p \in (1,\infty)$, $r \in [1,p) $, $\gamma^{p} \in
A_{\frac{p}{r}}^{\text{\rm loc}}(\mathbb{R}^{n+d})$, $f \in
W_{p}^{l}(\mathbb{R}^{n+d},\gamma)$, $l > \frac{d}{r}$. Then there exists the trace $\varphi \in
\widetilde{B}^{l}_{p,p,r}(\mathbb{R}^{n},\{\gamma_{k}\})$ of the function~$f$, and moreover,
\begin{equation}
\label{3.4}
\left\|\varphi|\widetilde{B}^{l}_{p,p,r}(\mathbb{R}^{n},\{\gamma_{k}\})\right\|
\le C_{1}\|f|W_{p}^{l}(\mathbb{R}^{n+d},\gamma)\|.
\end{equation}
The constant $C_{1}$ in  \eqref{3.4}  is independent of the function~$f$.

Conversely, if a function $\varphi \in
\widetilde{B}^{l}_{p,p,r}(\mathbb{R}^{n},\{\gamma_{k}\})$, then there exists a~function $f \in
W_{p}^{l}(\mathbb{R}^{n+d},\gamma)$ such that $\varphi$ is the trace of~$f$ on $\mathbb{R}^{n}$, and moreover,
\begin{equation}
\label{3.5}
\|f|W_{p}^{l}(\mathbb{R}^{n+d},\gamma)\|
\le C_{2} \left\|\varphi|\widetilde{B}^{l}_{p,p,r}(\mathbb{R}^{n},\{\gamma_{k}\})\right\|,
\end{equation}
the constant $C_{2}$ in \eqref{3.5} being independent of the function~$\varphi$.
\end{Th}

The following result in an important step in the proof of Theorem~\ref{Th3.1}.

\begin{Lm}
\label{Lm3.2}
Let $p \in (1,\infty)$, $r \in [1,p)$, $\gamma^{p} \in A_{\frac{p}{r}}^{\text{\rm loc}}(\mathbb
R^{n+d})$, $f \in W_{p}^{l}(\mathbb{R}^{n+d},\gamma)$, $l
> \frac{d}{r}$. Then there exists the trace~$\varphi$ of the function~$f$ on~$\mathbb{R}^{n}$. Moreover, for an arbitrary cube $Q^{n}_{k,m}$,
\begin{equation}
\begin{split}
\label{3.6}
&\delta^{l}_{r}(Q_{k,m}^{n})\varphi \le
\frac{C_{3}}{2^{k(l-\frac{n+d}{r})}}\left\{\left
\|D^{\alpha}_{x}f|L_{r}\left(C_{4}Q_{k,m}^{n}\times
\frac{C_{5}}{2^{k}}B^{d}\right)\right\|+\left
\|D^{\beta}_{y}f|L_{r}\left(C_{4}Q_{k,m}^{n}\times
\frac{C_{5}}{2^{k}}B^{d}\right)\right\|\right\}.
\end{split}
\end{equation}
The constants $C_{3},C_{4},C_{5}$ in~\eqref{3.6} are independent of both the function~$f$ and the cube $Q^{n}_{k,m}$.
\end{Lm}

\textbf{The proof} of the lemma is a straightforward modification of that of Lemma 3.1 of~\cite{Ty}.

\textbf{Proof of Theorem 3.1}.

\textit{Step 1}. Let $f \in W^{l}_{p}(\mathbb{R}^{n+d}, \gamma)$. Then by Lemma  \ref{Lm3.2} there exists the trace~$\varphi$ of the function~$f$
on the plane $\mathbb{R}^{n}$, and so it suffices
to prove estimate~\eqref{3.4}. In turn, this estimate follows by an obvious modification of the argument made in the proof of Theorem~3.1 of~\cite{Ty}
(one needs only to recourse to estimate \eqref{3.6} in an appropriate place).

\textit{Step 2}. We shall carry out this part of the proof in detail. As distinct from~\cite{Ty}, we shall need to estimate all mixed derivatives.
So let $\{\psi_{k}\}_{k=0}^{\infty}$ be a~partition of unity for the ball $B^{d}$.
Note that $\psi_{0} \in C^{\infty}(B^{d} \setminus \frac{1}{2}B^{d})$, $\psi_{k} \in
C_{0}^{\infty}(\frac{1}{2^{k-1}}B^{d} \setminus \frac{1}{2^{k+1}}B^{d})$ for $k \in \mathbb{N}$ and $|D^{\beta}\psi_{k}(y)| \le
\frac{C_{1}}{(\delta_{k})^{|\beta|}}$ for $y \in B^{d}$, $k \in
\mathbb{N}_{0}$. Assume that, for any $k \in \mathbb{N}_{0}$, only two functions $\psi_{k}$ and $\psi_{k+1}$ do not vanish on the set
$2^{-k}B^{d} \setminus 2^{-k-1}B^{d}$. Hence,
$D^{\beta}\psi_{k}(y)=-D^{\beta}\psi_{k+1}(y)$ for $y \in 2^{-k}B^{d} \setminus 2^{-k-1}B^{d}$.

The existence of a~sequence $\{\psi_{k}\}_{k=0}^{\infty}$ with the above properties may be proved as it was done, for example,
in \S\,4.5 of the book~\cite{10} in the proof of the trace theorem for unweighted Sobolev spaces.

We set
$$
f(x,y):=\sum\limits_{k=1}^{\infty}\psi_{k}(y)E_{2^{-k}}[\varphi](x) , \qquad  (x,y) \in \mathbb{R}^{n} \times B^{d},
$$
where, the operator $E_{\varepsilon}$ (with  $\varepsilon > 0$) is defined in \eqref{2.28}. We extend the function $f$ by zero on the set
$^{n}\mathbb{R}^{d}_{1}$.

A multi-index $\alpha$ will be written as $(\alpha^{1},\alpha^{2})=(\alpha^{1}_{1},\dots,\alpha^{1}_{n},\alpha^{2}_{1},\dots,\alpha^{2}_{d})$.

Clearly,
\[
\iint\limits_{\mathbb{R}^{n}\times\frac{1}{2}B^{d}}\gamma^{p}(x,y)\{\sum\limits_{|\alpha|=l,\alpha^{2}=0}|D^{\alpha}f(x,y)|^{p}+\sum\limits_{|\alpha|=l,|\alpha^{2}|>0}|D^{\alpha}f(x,y)|^{p}\}\,dxdy=
\]

\[
=\sum\limits_{k=1}^{\infty}\sum\limits_{m \in
\mathbb{Z}^{n}}\iint\limits_{\Xi_{k,m}^{d,n}}\gamma^{p}(x,y)\{\sum\limits_{|\alpha|=l,\alpha^{2}=0}|D^{\alpha}f(x,y)|^{p}+\sum\limits_{|\alpha|=l,|\alpha^{2}|>0}|D^{\alpha}f(x,y)|^{p}\}\,dxdy.
\]

Taking into account properties of the functions $\psi_{k}$ and applying estimate \eqref{2.30}, we see that$$
\sum\limits_{|\alpha|=l,\alpha^{2}>0}\iint\limits_{\Xi_{k,m}^{d,n}}\gamma^{p}(x,y)|D^{\alpha}f(x,y)|^{p}\,dxdy=\\
$$
$$
=\sum\limits_{|\alpha|=l,\alpha^{2}>0}\iint\limits_{\Xi_{k,m}^{d,n}}\gamma^{p}(x,y)|D^{\alpha^{2}}_{y}\psi_{k}(y)D^{\alpha^{1}}_{x}E_{2^{-k}}\varphi(x)+D^{\alpha^{2}}_{y}\psi_{k+1}(y)D^{\alpha^{1}}_{x}E_{2^{-(k+1)}}\varphi(x)|^{p}\,dxdy
\le
$$
$$
\le \sum\limits_{|\alpha^{1}|=l-|\alpha^{2}|}2^{k|\alpha^{2}|p}\iint\limits_{\Xi_{k,m}^{d,n}}\gamma^{p}(x,y)|D^{\alpha^{1}}_{x}E_{2^{-k}}\varphi(x)-D^{\alpha^{1}}_{x}E_{2^{-(k+1)}}\varphi(x)|^{p}\,dxdy
\le
$$
\begin{equation}
\begin{split}
\label{3.7}
&\le C
2^{2nkp}\gamma^{p}_{k,m}\biggl[\int\limits_{\widetilde{C}Q_{k,m}^{n}}\int\limits_{\frac{1}{2^{k}}I^{n}}|\Delta^{l}(h)\varphi(z)|\,dhdz\biggr]^{p}\,dxdy \mbox{ for } k \in \mathbb{N},m \in \mathbb{Z}^{n}.
\end{split}
\end{equation}

The constant $\widetilde{C} \geq 1$, which is the dilation coefficients of the cubes  $Q^{n}_{k,m}$, depends only on the diameter of the support of the function~$\Theta$ from~\eqref{3.1}.

Similarly, it follows from \eqref{3.2} that
\begin{equation}\label{3.8}
\begin{gathered}
\sum\limits_{|\alpha^{1}|=l}\iint\limits_{\Xi_{k,m}^{d,n}}\gamma^{p}(x,y)|D^{\alpha^{1}}_{x}f(x,y)|^{p}\,dxdy
\le  \\
\le C \sum\limits_{|\alpha^{1}|=l} \iint\limits_{\Xi_{k,m}^{d,n}}\gamma^{p}(x,y)\max
\{|D^{\alpha^{1}}_{x}E_{2^{-k}}\varphi(x)|^{p},
|D^{\alpha^{1}}_{x}E_{2^{-k-1}}\varphi(x)|^{p} \}\,dxdy \le \\
 \le C
2^{2npk}\gamma^{p}_{k,m}\biggl[\int\limits_{\widetilde{C} Q_{k,m}^{n}}\int\limits_{\frac{1}{2^{k}}I^{n}}|\Delta^{l}(h)\varphi(z)|\,dhdz\biggr]^{p}, \qquad
 k \in \mathbb{N},m \in \mathbb{Z}^{n}.
\end{gathered}
\end{equation}

Using the definition of the function $f$ and employing H\"older's inequality with exponents $r$, $r'$, we have, for  $|\alpha|=l$
\begin{gather}
\sum\limits_{|\alpha|=l}\iint\limits_{^{n}\mathbb{R}^{d}_{\frac{1}{2}}}\gamma^{p}(x,y)|D^{\alpha}f(x,y)|^{p}\,dxdy
\le C \sum\limits_{m \in \mathbb{Z}^{n}}\gamma^{p}_{0,m}\|\varphi|L_{1}(\widetilde{C} Q^{n}_{0,m})\|^{p} \le C \sum\limits_{m \in \mathbb{Z}^{n}}\gamma^{p}_{0,m}\|\varphi|L_{r}(\widetilde{C} Q^{n}_{0,m})\|^{p} \le \notag \\
\le C \sum\limits_{m \in \mathbb{Z}^{n}}\gamma^{p}_{0,m}\|\varphi|L_{r}( Q^{n}_{0,m})\|^{p},
\label{3.9}
\end{gather}
since the cubes $\widetilde{C} Q^{n}_{k,m}$ have finite overlapping multiplicity  (the constant  $\widetilde{C}$ is the same as in \eqref{3.7}).

Hence, summing up estimates \eqref{3.7}, \eqref{3.8} in $k$ and~$m$, taking into that the
cubes $nQ_{k,m}^{n}$ have finite overlapping multiplicity (for $n \in \mathbb{N}$), and employing estimate \eqref{3.9}, this gives
\begin{equation}
\label{3.10}
\sum\limits_{|\alpha|=l}\|D^{\alpha}f|L_{p}(\mathbb{R}^{n+d},\gamma)\| \le C \|\varphi|\widetilde{B}^{l}_{p,p,r}(\mathbb{R}^{n},\{\gamma_{k,m}\})\|.
\end{equation}

To estimate the generalized derivatives $D^{\alpha}f$ for $|\alpha| < l$ we write, for each $(x,y) \in \mathbb{R}^{n} \times B^{d}$, the integral representation
of the function $D^{\alpha}f$ in a~cone (see \S\,3.4, \cite{10}),
$V(x,y)=\{(x,y)(1-t)+t(x',y')|t \in [0,1], (x',y') \in \frac{1}{2}B^{n+d}(x,y+3)\}$ (here, $\frac{1}{2} B^{n+d}(x,y+3)$ is the ball of radius $\frac{1}{2}$
centred at~$(x,y+3)$), and use Remark~16 of \S\,3.5 in~\cite{10}.

Let $|\alpha| < l$. Since $f(x,y)=0$ for $|y| > 1$, we have
$$
|D^{\alpha}f(x,y)| \le C \sum\limits_{|\beta|=l} \iint\limits_{(x,0)+(I^{n} \times B^{d})}|D^{\beta}f(\widetilde{x},\widetilde{y})|\,d\widetilde{x}d\widetilde{y} ,
\quad  (x,y) \in \mathbb{R}^{n} \times B^{d}.
$$

Hence, using the obvious inclusion $A^{\text{\rm loc}}_{\frac{p}{r}}(\mathbb{R}^{n+d}) \subset A^{\text{\rm loc}}_{p}(\mathbb{R}^{n+d})$,
and employing H\"older's inequality, we obtain, for $m \in \mathbb{Z}^{n}$, $|\alpha| < l$,
\begin{equation}
\begin{split}
\label{3.11}
&\iint\limits_{Q^{n}_{0,m} \times B^{d}}\gamma^{p}(x,y)|D^{\alpha}f(x,y)|^{p}\,dxdy \le \\
& \le C \sum\limits_{|\beta|=l}\left[\gamma^{p}(\widetilde{C} Q^{n}_{0,m} \times B^{d})\right]\left[\gamma^{-p'}(\widetilde{C} Q^{n}_{0,m} \times B^{d})\right]^{\frac{p}{p'}}\iint\limits_{\widetilde{C} Q^{n}_{0,m} \times B^{d}}\gamma^{p}(x,y)|D^{\beta}f(x,y)|^{p}\,dxdy \\
& \le C \sum\limits_{|\beta|=l}\iint\limits_{\widetilde{C} Q^{n}_{0,m} \times B^{d}}\gamma^{p}(x,y)|D^{\beta}f(x,y)|^{p}\,dxdy.
\end{split}
\end{equation}

Summing up estimate \eqref{3.11} over $m \in \mathbb{Z}^{n}$ and taking into account the finite multiplicity
of the cubes $\widetilde{C} Q_{k,m}^{n}$  we obtain \eqref{3.5} in view of~\eqref{3.10}.

It remains to show that $\varphi=\operatorname{tr}\left|_{y=0}\right.f$.
We fix an arbitrary cube $Q^{n}$. Almost every point $x \in \mathbb{R}^{n}$ is a~Lebesgue point of the function $\varphi$, because $\varphi \in
L^{\text{\rm loc}}_{1}(\mathbb{R}^{n})$. Hence $g_{\delta}(x):=\lim\limits_{\delta \to 0}\frac{1}{\delta^{n}}\int\limits_{x+\delta I^{n}}|\varphi(x')-\varphi(x)|\,dx'=0$
for almost all $x \in \mathbb{R}^{n}$. Consequently, by the Lebesgue convergence theorem,
\begin{equation}
\label{3.12}
\int\limits_{Q^{n}}|\varphi(x)-E_{\delta}[\varphi](x)|\,dx \le
C\int\limits_{Q^{n}}g_{\delta}(x) dx \to 0 \mbox{as} \delta \to 0.
\end{equation}

From \eqref{3.12} and the definition of the function~$f$ it easily follows that $\varphi$ is the trace of the function $f$ on the plane $y=0$.

The proof of the theorem is complete.

\section{Atomic decomposition of functions from the spaces $\widetilde{B}^{l}_{p,q,r}(\mathbb{R}^{n},\{t_{k,m}\})$}

Our aim in this section is to prove the atomic decomposition theorem for functions~$\varphi$ from the space
$\widetilde{B}^{l}_{p,q,r}(\mathbb{R}^{n},\{t_{k}\})$.
This theorem is one of the principal tools in establishing various embedding
and trace theorems (see \S\S\,5 and~6 below).

We shall also put forward conditions on a~$p$-admissible weight sequence $\{t_{k,m}\}$ securing the coincidence of
the spaces $\widetilde{B}^{l}_{p,q,r}(\mathbb{R}^{n},\{t_{k}\})$ for various $l \in \mathbb{N}$ with equivalence of
the corresponding norms.

Our arguments will depend to a large extent on the methods of the paper~\cite{DeVore}.

We fix numbers $n, d \in \mathbb{N}$ and define
$\Xi^{d,n}_{k,m}:=Q^{n}_{k,m}\times(\frac{B^{d}}{2^{k}}\setminus \frac{B^{d}}{2^{k+1}})$ for $k \in \mathbb{N}_{0}, m \in \mathbb{Z}^{n}$.
A~point of an $(n+d)$-dimensional Euclidean space will be written as a~pair $(x,y):=(x_{1},\dots,x_{n},y_{1},\dots,y_{d})$.

Let $l \in \mathbb{N}$, $N^{l-1}$ be a~$B$-spline of degree $l-1$ with knots at the points  $t_{i}=i$, $i \in \{0,1,..,l\}$. More precisely,
$$
N^{l-1}(t):=[0,1,..,l](t-\cdot)_{+}^{l-1}.
$$
Here, we use the standard notation for the divided difference (see [7, Ch.~1]).

For $k \in \mathbb{N}_{0}$, $m \in \mathbb{Z}^{n}$, we set
$$
N^{l-1}_{k,m}(x):=\prod\limits_{i=1}^{n}N^{l-1}(2^{k}(x_{i}-\frac{m_{i}}{2^{k}})), \mbox{ for } x \in \mathbb{R}^{n}.
$$

The functions $N^{l-1}_{k,m}$ were first introduced by Curry and Schoenberg~\cite{Carry}.

We list some properties of the $B$-splines $N^{l-1}_{k,m}$ that will be required in what follows.
The corresponding proofs may be found, for example, in~\cite{Carry},~\cite{Fix}.

1) The $B$-splines $N^{l-1}_{k,m}$ form a~partition of unity on~$\mathbb{R}^{n}$ for each fixed
$k \in \mathbb{N}_{0}$. That is,
\begin{equation}
\label{4.1}
\sum\limits_{m \in \mathbb{Z}^{n}}N^{l-1}_{k,m}(x)=1, \mbox{ for } x \in
\mathbb{R}^{n}.
\end{equation}
Here, the overlapping multiplicity of the supports of splines $N^{l-1}_{k,m}$ is finite
and is independent of both~$k$ and~$m$. We also note that $\operatorname{supp} N^{l-1}_{k,m} \subset
\frac{m}{2^{k}}+[0,\frac{l}{2^{k}}]^{n}$ and $N^{l-1}_{k,m}(x) \in (0,1]$ for $x
\in \frac{m}{2^{k}}+(0,\frac{l}{2^{k}})^{n}$.

2) On each cube $Q^{n}_{k,m}$ the function $N^{l-1}_{k,m}$ is a polynomial of degree $\le l-1$ in each variable.

3) The spline $N^{l-1}$ has continuous derivative of order $l-2$. At knots $t_{i}=i$, $i \in \{0,1,\dots,l\}$, the spline $N^{l-1}$
has finite one-sided derivatives of order $l-1$. Hence,
\begin{equation}
\label{4.2}
\Delta^{l}(h)N^{l-1}_{k,m}(x) \le C (2^{k}|h|)^{l-1}, \qquad  x,h \in
\mathbb{R}^{n}.
\end{equation}

4) Any spline $S=\sum\limits_{m \in Z^{n}}\beta_{k,m}N^{l-1}_{k,m}$ may be expanded into a~series in splines $N^{l-1}_{j,m}$ for $j \geq k$; that is,
$S=\sum\limits_{m \in Z^{n}}\widehat{\beta}_{k,m}(S)N^{l-1}_{j,m}$.

\smallskip

We let $\Sigma^{l-1}_{k}$ denote the set of all splines $S$ of the form
$$
S(x):=\sum\limits_{m \in \mathbb{Z}^{n}} \beta_{k,m} N^{l-1}_{k,m}(x) \mbox{ for } x\in \mathbb{R}^{n}.
$$

For future purposes we shall require the concept of a~quasi-interpolant, which was first introduced in~\cite{Fix}.
Quasi-interpolants were also used in the papers \cite{DeVore}, \cite{Ir} for constructing equivalent norms on unweighted dyadic Besov-type spaces and on classical Besov spaces.

\begin{Def}(\cite{Fix})
\label{Def4.1}
Given $k \in \mathbb{N}_{0}$, $m \in \mathbb{Z}^{n}$, we let $\xi_{k,m}=(\xi_{k,m_{1}},\dots,\xi_{k,m_{n}})$ denote the centre of the cube $Q^{n}_{k,m}$. Assume that
all partial derivatives $D^{\nu}f$, $\nu_{j} \le l-1$, $j \in \{1,\dots,n\}$ of~$f$ are continuous at each point $\xi_{k,m}$.
For $k \in \mathbb{N}_{0}$, $m \in \mathbb{Z}^{n}$, we set
$$
Q^{l-1}_{k}(f):=\sum\limits_{m \in \mathbb{Z}^{n}}\alpha_{k,m}(f)N^{l-1}_{k,m}, \qquad \mbox{ where }
$$
$$
\alpha_{k,m}(f):=\sum\limits_{\substack{0 \le \nu_{j} \le l-1\\
j \in \{1,..,n\}}}a_{k,\nu,m}D^{\nu}f(\xi_{k,m}),
$$
$$
a_{k,m,\nu}:=\prod_{i=1}^{n}a_{k,m_{i},\nu_{i}}, a_{k,m_{i},\nu_{i}}:=\frac{(-1)^{l-1-\nu_{i}}}{(l-1)!}D^{l-1-\nu_{i}}\psi_{m_{i}}(\xi_{k,m_{i}}),
$$
$$
\psi_{m_{i}}(t):=\prod\limits_{j=1}^{l-1}(\frac{m_{i}+j}{2^{k}}-t) \mbox{ for } t \in \mathbb{R}.
$$
The operator $Q^{l-1}_{k}$ is called a quasi-interpolant.
\end{Def}

Throughout this section we fix a~constant $A \geq 1$.

Let $P_{Q^{n}_{k,m}}$ be a~polynomial of almost best approximation in the $L_{r}(Q^{n}_{k,m})$-metric to a~function $\varphi \in L_{r}^{\text{\rm loc}}(\mathbb{R}^{n})$ (for $r \in (0,\infty]$) by
 polynomials of \textit{coordinate degree} $<l$ on the cube $Q^{n}_{k,m}$ with constant $A$.

We set $g_{k}(x):=\sum\limits_{m \in \mathbb{Z}^{n}}P_{Q^{n}_{k,m}}(x)\chi_{Q^{n}_{k,m}}(x)$ for $x \in \mathbb{R}^{n}$, $k \in \mathbb{N}_{0}$.
Finally, following \cite{DeVore} we define, for $r \in (0,\infty]$,
\begin{equation}
\begin{split}
\label{4.3}
&T^{l-1}_{k}(\varphi,r)(x):=Q^{l-1}_{k}(g_{k})(x), \mbox{ for } x \in \mathbb{R}^{n}, k \in \mathbb{N}_{0},\\
&T^{l-1}_{-1}(\varphi,r)(x):=\varphi(x),\mbox{ for } x \in \mathbb{R}^{n}.
\end{split}
\end{equation}

We note that in \cite{DeVore} the operator $T^{l-1}_{k}$ acts on functions~$\varphi$ defined on the unit cube.

\begin{Remark}
\label{R4.1}
The operator $Q^{l-1}_{k}$ is a projection operator from the space of piecewise-polynomial functions to the space $\Sigma^{l-1}_{k}$ (for the proof we refer to~\cite{Fix}),
and hence
$$
T^{l-1}_{k}(\varphi,r)(x)=\sum\limits_{m \in \mathbb{Z}^{n}}\alpha_{k,m}(T^{l-1}_{k}(\varphi,r))N^{l-1}_{k,m}(x) \mbox{ for } x \in \mathbb{R}^{n}.
$$
\end{Remark}

\begin{Lm}
\label{Lm4.1}
Let $r \in (0,\infty]$. Then for any function $\varphi \in L_{r}^{\text{\rm loc}}(\mathbb{R}^{n})$  and $k \in \mathbb{N}_{0}$ the following estimate holds:
\begin{equation}
\label{4.4}
\|\varphi-T^{l-1}_{k}(\varphi,r)|L_{r}(Q^{n}_{k,m})\| \le C \widehat{E}_{l}(\varphi,(1+l)Q^{n}_{k,m})_{r} \le C E_{l}(\varphi,(1+l)Q^{n}_{k,m})_{r}.
\end{equation}
The constant $C$ in \eqref{4.4}  depends only on  $l,n,r,A$.
\end{Lm}

\textbf{Proof}. The first inequality in \eqref{4.4} follows from estimate (4.25) of~\cite{DeVore}, the second inequality  is clear.

Let $p,r \in (0,\infty]$,  $\varphi \in L^{\text{\rm loc}}_{r}(\mathbb{R}^{n})$, $\alpha_{3} \geq 0$, $\alpha_{i} \in \mathbb{R}$, $\sigma_{i} \in (0,\infty]$ ($i=1,2$).
For a~multiple sequence  $\{t_{k,m}\}$ which is $p$-associated with $p$-admissible  weight sequence $\{t_{k}\} \in X^{\alpha_{3}}_{\alpha,\sigma,p}$, we set
\begin{equation}
\begin{split}
\label{4.5}
&s^{l}_{k}:=s^{l}_{k}(\varphi,\{t_{k,m}\})_{r,p}:=\inf\limits_{S \in
\Sigma^{l-1}_{k}}\left(\sum\limits_{m \in
\mathbb{Z}^{n}}t^{p}_{k,m}\|\varphi-S|L_{r}(Q^{n}_{k,m})\|^{p}
\right)^{\frac{1}{p}} \mbox{ for } k \in \mathbb{N}_{0},\\
&s^{l}_{-1}:=s^{l}_{-1}(\varphi,\{t_{0,m}\})_{r,p}=\left(\sum\limits_{m \in
\mathbb{Z}^{n}}t^{p}_{0,m}\|\varphi|L_{r}(Q^{n}_{k,m})\|^{p}
\right)^{\frac{1}{p}}.
\end{split}
\end{equation}

We note that $s^{l}_{k}(\varphi,\{t_{k,m}\})_{r,p}< \infty$ for $\varphi \in \widetilde{B}^{l}_{p,q,r}(\mathbb{R}^{n},\{t_{k,m}\})$. Indeed,
in view of \eqref{4.4} and Theorem~\ref{Th2.3}, we have, for $k \in \mathbb{N}_{0}$,
\begin{equation}
\begin{split}
\label{4.6}
&\inf\limits_{S \in
\Sigma^{l-1}_{k}}\left(\sum\limits_{m \in
\mathbb{Z}^{n}}t^{p}_{k,m}\|\varphi-S|L_{r}(Q^{n}_{k,m})\|^{p}
\right)^{\frac{1}{p}} \le \left(\sum\limits_{m \in
\mathbb{Z}^{n}}t^{p}_{k,m}\|\varphi-T_{k}^{l-1}(\varphi)|L_{r}(Q^{n}_{k,m})\|^{p}
\right)^{\frac{1}{p}} \le \\
&\le C\|\varphi|\widetilde{B}^{l}_{p,q,r}(\mathbb{R}^{n},\{t_{k}\})\| < \infty.
\end{split}
\end{equation}

\begin{Def}
\label{Def4.2}
Let $p,r \in (0,+\infty]$,  $s^{l}_{k}(\varphi,\{t_{k,m}\})_{r,p} < \infty$. We say that $U^{l-1}_{k}:=U^{l-1}_{k}(\varphi,\{t_{k,m}\},p) \in
\Sigma^{l-1}_{k}$
is a~spline of almost best approximation with constant~$A$ to a~function $\varphi \in L^{\text{\rm loc}}_{r}(\mathbb{R}^{n})$  if
\begin{equation}
\label{4.7}
\left(\sum\limits_{m \in
\mathbb{Z}^{n}}t^{p}_{k,m}\|\varphi-U^{l-1}_{k}|L_{r}(Q^{n}_{k,m})\|^{p}
\right)^{\frac{1}{p}} \le A s^{l}_{k}(\varphi,\{t_{k,m}\})_{r,p}.
\end{equation}
\end{Def}

\begin{Lm}(\cite{DeVore})
\label{Lm4.2}
Assume that a spline $S \in \Sigma^{l-1}_{k}$. Then, for any $r \in (0,+\infty]$ and any cube $Q^{n}_{k,m}$,
\begin{equation}
\label{4.8}
C_{1}\|S|L_{r}(Q^{n}_{k,m})\| \le
\left(\sum\limits_{\substack{\widetilde{m} \in \mathbb{Z}^{n}\\
Q^{n}_{k,m} \bigcap \operatorname{supp} N^{l-1}_{k,\widetilde{m}} \neq
\emptyset}}|\alpha_{k,\widetilde{m}}(S)|^{r}2^{-kn}\right)^{\frac{1}{r}} \le
C_{2}\|S|L_{r}(C_{3}Q^{n}_{k,m})\|,
\end{equation}
the constants $C_{1},C_{2},C_{3} > 0$ being independent of both the cube  $Q^{n}_{k,m}$ and the spline~$S$.
The corresponding modifications in the case $r=\infty$ are clear.
\end{Lm}

The next theorem is an extension of Theorem~4.8 of \cite{DeVore} (which was concerned with classical Besov spaces) to the case of Besov spaces of variable smoothness.

We recall that the symbols $\delta_{i}(\gamma,n,d)$ ($i=1,2$) were introduced right after the completion of the proof of Lemma~2.2.

\begin{Th}
\label{Th4.1}
Let $p,r \in (0,\infty)$,  $\varphi \in L_{r}^{\text{\rm loc}}(\mathbb{R}^{n})$, $\alpha_{3} \geq 0$, $0 < \alpha_{1} \le \alpha_{2}$, a~weight sequence $\{s_{k}\} \in ^{\text{\rm loc}}Y^{\alpha_{3}}_{\alpha_{1},\alpha_{2}}$. Assume that a~weight $\gamma^{p} \in A^{\text{\rm loc}}_{\infty}(\mathbb{R}^{n+d})$ and
a~multiple  sequence $\{\widehat{\gamma}_{k,m}\}$ is generated by the weight~$\gamma$.
Next, let  $t_{k,m}:=2^{\frac{kn}{p}}s_{k,m}\widehat{\gamma}_{k,m}$ for $k \in \mathbb{N}_{0}$, $m \in \mathbb{Z}^{n}$.

Then, for any sufficiently small $\varepsilon > 0$,
\begin{equation}
\begin{split}
\label{4.9}
&\left(\sum\limits_{m \in
\mathbb{Z}^{n}}t^{p}_{k,m}[\delta^{l}_{r}(Q^{n}_{k,m})\varphi]^{p}\right)^{\frac{1}{p}} \le\\
 &\le C 2^{-k(\widetilde{\lambda}+\frac{d(\delta_{1}(\gamma,n,d)-\varepsilon)}{p}-\alpha_{2})}\left(\sum\limits_{j=-1}^{k}2^{j\mu(\widetilde{\lambda}+\frac{d(\delta_{1}(\gamma,n,d)-\varepsilon)}{p}-\alpha_{2})}\left(s_{j}^{l}(\varphi,\{t_{k,m}\})_{r,p}\right)^{\mu}\right)^{\frac{1}{\mu}},
\end{split}
\end{equation}
where $\widetilde{\lambda}:= \min
\{l,l-1+\frac{\delta_{2}(\gamma,n,d)-\varepsilon}{p}\}$, $\mu \le \min\{1,r,p\}$,
and the constant $C > 0$ depends on $\alpha_{1},\alpha_{2},\alpha_{3}$, $r$, $l$ and the weight~$\gamma$,
but is independent of the function~$\varphi$.
\end{Th}

\textbf{Proof.} The main idea of the proof of Theorem~\ref{Th3.1} follows that of Theorem~4.8 in~\cite{DeVore}.
However, certain modifications of the proof of~\cite{DeVore} are required to account for the properties of the multiple sequence $\{\widehat{\gamma}_{k,m}\}$,
which were indicated in Lemma~\ref{Lm2.1}.

We fix $\varepsilon \in [0,\min\{\delta_{1}(\gamma,n,d),\delta_{2}(\gamma,n,d\})$ and define  $\widetilde{\delta_{1}}:=\widetilde{\delta_{1}}(\gamma,n,d):=\delta_{1}(\gamma,n,d)-\varepsilon$, $\widetilde{\delta_{2}}:=\widetilde{\delta_{2}}(\gamma,n,d):=\delta_{2}(\gamma,n,d)-\varepsilon$.

Let $U^{l-1}_{j}:=U^{l-1}_{j}(\varphi,\{t_{k,m}\},p)$ be a~spline of almost best approximation with constant $A\geq1$ to a~function
$\varphi \in L^{\text{\rm loc}}_{r}(\mathbb{R}^{n})$. Given
$j \in \mathbb{N}_{0}$, we set $u^{l-1}_{j}:=U^{l-1}_{j}-U^{l-1}_{j-1}$ ($U^{l-1}_{-1}
\equiv 0$). Then, clearly,

\begin{equation}
\label{4.10}
\Delta^{l}(h)\varphi(x)=\Delta^{l}(h)[\varphi-U^{l-1}_{k}](x)+\sum\limits_{j=0}^{k}\Delta^{l}(h)u^{l-1}_{j}(x),
\qquad  x,h \in \mathbb{R}^{n}.
\end{equation}

The inequality
$$
\delta^{l}_{r}(Q^{n}_{k,m})\varphi \le \left([\delta^{l}_{r}(Q^{n}_{k,m})(\varphi-U^{l-1}_{k})]^{\mu}+\sum\limits_{j=0}^{k}[\delta^{l}_{r}(Q^{n}_{k,m})u^{l-1}_{j}]^{\mu}\right)^{\frac{1}{\mu}}
$$
is an easy consequence of \eqref{4.10} and Lemma \ref{2.1} with $\mu \le \min\{1,r,p\}$.
Hence, since $\frac{p}{\mu} \geq 1$ and using Minkowski's inequality,
\begin{equation}
\begin{split}
\label{4.11}
&\left(\sum\limits_{m \in
\mathbb{Z}^{n}}t^{p}_{k,m}[\delta^{l}_{r}(Q^{n}_{k,m})\varphi]^{p}\right)^{\frac{1}{p}}
\le \\
&\le   \left(\sum\limits_{m \in
\mathbb{Z}^{n}}t^{p}_{k,m}\left([\delta^{l}_{r}(Q^{n}_{k,m})(\varphi-U^{l-1}_{k})]^{\mu}+\sum\limits_{j=0}^{k}[\delta^{l}_{r}(Q^{n}_{k,m})u^{l-1}_{j}]^{\mu}\right)^{\frac{p}{\mu}}\right)^{\frac{\mu}{p}{\frac{1}{\mu}}} \\
&\le  \left(\left(\sum\limits_{m \in
\mathbb{Z}^{n}}t^{p}_{k,m}[\delta^{l}_{r}(Q^{n}_{k,m})(\varphi-U^{l-1}_{k})]^{p}\right)^{\frac{\mu}{p}}+\sum\limits_{j=0}^{k}\left(\sum\limits_{m \in
\mathbb{Z}^{n}}t^{p}_{k,m}[\delta^{l}_{r}(Q^{n}_{k,m})u^{l-1}_{j}]^{p}\right)^{\frac{\mu}{p}}\right)^{\frac{1}{\mu}}\\
&=\left((R^{1})^{\frac{\mu}{p}}+\sum\limits_{j=0}^{k} (R^{2}_{j})^{\frac{\mu}{p}}\right)^{\frac{1}{\mu}}.
\end{split}
\end{equation}

Since the cubes $(1+l)Q^{n}_{k,m}$ have finite overlapping multiplicity (which is independent of $k$ and~$m$) and using Remark~\ref{R2.2}, it is easily seen that
\begin{equation}
\label{4.12}
(R^{1})^{\frac{\mu}{p}} \le C \left(\sum\limits_{m \in
\mathbb{Z}^{n}}t^{p}_{k,m}2^{\frac{knp}{r}}\|\varphi-U^{l-1}_{k}|L_{r}((1+l)Q^{n}_{k,m})\|^{p}\right)^{\frac{\mu}{p}} \le C \left(s_{k}^{l}(\varphi,\{t_{k,m}\})_{r,p}\right)^{\mu}.
\end{equation}

It is worth noting that the differences $\delta^{l}_{r}$ (rather than $\overline{\Delta}^{l}_{r}$)
were crucial in obtaining estimate~\eqref{4.12} .

Next, for each $j \in \mathbb{N}_{0}$ the function $u^{l-1}_{j}$ may be expanded into a~series in $B$-splines $N^{l-1}_{j,m}$ (by property~4), see above); that is,
\begin{equation}
\label{4.13}
u^{l-1}_{j}(x)=\sum\limits_{m \in Z^{n}}\alpha_{j,m}(u^{l-1}_{j})N^{l-1}_{j,m}(x), \mbox{ for } x \in \mathbb{R}^{n}.
\end{equation}

Since, for any point $x \in \mathbb{R}^{n}$, only a~finite number (independent of $j$ and $m$) of splines $N^{l-1}_{j,m}$ are nonzero, we have
\begin{equation}
\label{4.14}
|\Delta^{l}(h)u^{l-1}_{j}(x)|^{r} \le C \sum\limits_{x \in \mbox{ supp } N^{l-1}_{j,m}}|\alpha_{j,m}|^{r}|\Delta^{l}(h)N^{l-1}_{j,m}(x)|^{r}.
\end{equation}

For $k \geq j \in \mathbb{N}_{0}$, $\widetilde{m} \in \mathbb{Z}^{n}$, we let
$\Gamma_{j,\widetilde{m}}$ denote the set of all cubes $Q^{n}_{k,m} \subset Q^{n}_{j,\widetilde{m}}$.
Next, let $\Gamma^{1}_{j,\widetilde{m}}$ denote the set of all cubes $Q^{n}_{k,m} \subset Q^{n}_{j,\widetilde{m}}$ for which
$(1 +l)Q^{n}_{k,m} \subset Q^{n}_{j,\widetilde{m}}$. We also define $\Gamma^{2}_{j,\widetilde{m}}:=\Gamma_{j,\widetilde{m}} \setminus \Gamma^{1}_{j,\widetilde{m}}$.

For further purposes, we shall require the following estimate of the measure of the set
$F_{j,\widetilde{m}}:=\cup_{Q^{n}_{k,m} \in \Gamma^{2}_{j,\widetilde{m}}}Q^{n}_{k,m}$ (the proof is similar to that of the corresponding estimate in~\cite{DeVore}).
\begin{equation}
\label{4.15}
\frac{|F_{_{j,\widetilde{m}}}|}{|Q^{n}_{j,\widetilde{m}}|} \le C 2^{j-k}.
\end{equation}

Using \eqref{4.15} and \eqref{2.11}, we obtain
\begin{equation}
\label{4.16}
\sum\limits_{Q^{n}_{k,m} \in \Gamma^{2}_{j,\widetilde{m}}}\widehat{\gamma}^{p}_{k,m} \le C 2^{(j-k)\widetilde{\delta}_{2}(\gamma,n,d)}\sum\limits_{Q^{n}_{k,m} \in \Gamma_{j,\widetilde{m}}}\widehat{\gamma}^{p}_{k,m}.
\end{equation}

For the cubes $Q^{n}_{k,m} \in \Gamma^{1}_{j,\widetilde{m}}$ we have
\begin{equation}
\label{4.17}
\delta^{l}_{r}(Q^{n}_{k,m})N^{l-1}_{j,\widetilde{m}} \le C 2^{(j-k)l},
\end{equation}
inasmuch as $N^{l-1}_{j,\widetilde{m}}$ is a polynomial on the cube $Q^{n}_{j,\widetilde{m}}$.

For the cubes $Q^{n}_{k,m} \in \Gamma^{2}_{j,\widetilde{m}}$, we have by~\eqref{4.2}
\begin{equation}
\label{4.18}
\delta^{l}_{r}(Q^{n}_{k,m})N^{l-1}_{j,\widetilde{m}} \le C 2^{(j-k)(l-1)}.
\end{equation}

Next, since $\{s_{k}\} \in ^{\text{\rm loc}}Y^{\alpha_{3}}_{\alpha_{1},\alpha_{2}}$ it clearly follows that
\begin{equation}
\label{4.19}
s_{k,m} \le C 2^{(\alpha_{2}-\frac{n}{p})(k-j)}s_{j,\widetilde{m}},
\end{equation}
provided that
$Q^{n}_{k,m} \subset Q^{n}_{j,\widetilde{m}}$, $k \geq j \in \mathbb{N}_{0}$, $m,\widetilde{m} \in \mathbb{Z}^{n}$ (the constant $C>0$ depends only on
the weight sequence $\{s_{k}\}$).

Combining estimates \eqref{4.8}, \eqref{4.14}, \eqref{4.17}, \eqref{4.18}, and using properties \eqref{2.10}, \eqref{4.16} and \eqref{4.19}
 of the multiple sequences $\{\widehat{\gamma}_{k,m}\}$ and $\{s_{k,m}\}$, this establishes
$$
R^{2}_{j} \le C \sum\limits_{\widetilde{m} \in \mathbb{Z}^{n}} \sum\limits_{Q^{n}_{k,m} \in \Gamma^{1}_{j,\widetilde{m}}}2^{p(k-j)\alpha_{2}}2^{\frac{jn}{p}}s^{p}_{j,\widetilde{m}}\widehat{\gamma}^{p}_{k,m}2^{(j-k)lp}\biggl[\sum\limits_{\substack{m \in \mathbb{Z}^{n}\\
Q^{n}_{j,\widetilde{m}} \bigcap \operatorname{supp}N^{l-1}_{j,m} \neq
\emptyset}}[\alpha_{j,m}]^{r}\biggr]^{\frac{p}{r}} +
$$

$$
+ C \sum\limits_{\widetilde{m} \in \mathbb{Z}^{n}} \sum\limits_{Q^{n}_{k,m} \in \Gamma^{2}_{j,\widetilde{m}}}2^{p(k-j)\alpha_{2}}2^{\frac{jn}{p}}s^{p}_{j,\widetilde{m}}\widehat{\gamma}^{p}_{k,m}2^{(j-k)(l-1)p}\biggl[\sum\limits_{\substack{m \in \mathbb{Z}^{n}\\
Q^{n}_{j,\widetilde{m}} \bigcap \operatorname{supp}N^{l-1}_{j,m} \neq
\emptyset}}[\alpha_{j,m}]^{r}\biggr]^{\frac{p}{r}} \le
$$

$$
\le C  \sum\limits_{m \in \mathbb{Z}^{n}}2^{p(k-j)\alpha_{2}}2^{\frac{jn}{p}}s^{p}_{j,m}\widehat{\gamma}^{p}_{j,m}2^{(j-k)d \widetilde{\delta}_{1}}2^{(j-k)pl}2^{\frac{jnp}{r}}\|u^{l-1}_{j}|L_{r}(C Q^{n}_{j,m})\|^{p} +
$$
$$
+ C \sum\limits_{m \in \mathbb{Z}^{n}}2^{p(k-j)\alpha_{2}}2^{\frac{jn}{p}}s^{p}_{j,m}\widehat{\gamma}^{p}_{j,m}2^{(j-k)(\widetilde{\delta}_{2}+d\widetilde{\delta}_{1})}2^{(j-k)p(l-1)}2^{\frac{jnp}{r}}\|u^{l-1}_{j}|L_{r}(C Q^{n}_{j,m})\|^{p} \le
$$

\begin{equation}
\begin{split}
\label{4.20}
&\le C 2^{(\widetilde{\lambda}+\frac{d \widetilde{\delta}_{1}}{p}-\alpha_{2})p(j-k)} \sum\limits_{m \in
\mathbb{Z}^{n}}
t^{p}_{j,m}2^{\frac{jnp}{r}}\|u^{l-1}_{j}|L_{r}(Q^{n}_{j,m})\|^{p} \le \\
&\le C 2^{(\widetilde{\lambda}+\frac{d \widetilde{\delta}_{1}}{p}-\alpha_{2})p(j-k)} (s_{j}^{l}(\varphi,\{t_{k,m}\})_{r,p})^{p}+(s_{j-1}^{l}(\varphi,\{t_{k,m}\})_{r,p})^{p}.
\end{split}
\end{equation}
Substituting estimates \eqref{4.12} and \eqref{4.20} into \eqref{4.11} completes the proof of the theorem.

We now prove a theorem similar to Theorem~\ref{Th4.1}, but under weaker constraints on a~multiple sequence
$\{t_{k,m}\}$. However, in doing so we increase the order of splines approximating a~given function.

\begin{Th}
\label{Th4.2}
Let $p,r \in (0,\infty]$, $\alpha_{1},\alpha_{2} \in \mathbb{R}$, $\alpha_{3} \geq 0$, $\sigma_{1} \in (0,\infty]$, $\sigma_{2}=p$, and let a~multiple sequence
$\{t_{k,m}\}$ be $p$-associated with $p$-admissible weight sequence $\{t_{k}\} \in X^{\alpha_{3}}_{\alpha,\sigma,p}$. Then, for $\mu \le \min\{1,r,p\}$,
\begin{equation}
\begin{split}
\label{4.21}
\left(\sum\limits_{m \in
\mathbb{Z}^{n}}t^{p}_{k,m}[\delta^{l}_{r}(Q^{n}_{k,m})\varphi]^{p}\right)^{\frac{1}{p}} \le
C 2^{-k(l-\alpha_{2})}\left(\sum\limits_{j=-1}^{k}2^{j\mu(l-\alpha_{2})}\left(s_{j}^{l+1}(\varphi,\{t_{k,m}\})_{r,p}\right)^{\mu}\right)^{\frac{1}{\mu}}
\end{split}
\end{equation}
{\rm (}with obvious modifications for $p=\infty$ or $\mu=\infty)$.
Here, the constant $C > 0$ depends on $\alpha_{1},\alpha,\sigma$, $r$, $l$,  but is independent of the function~$\varphi$.
\end{Th}

\textbf{Proof}. To a large extent we shall follow the proof of Theorem~\ref{Th4.1}. We only indicate the differences.

Clearly, for all $j \in \mathbb{N}_{0}$, $\widetilde{m} \in \mathbb{Z}^{n}$ and for all cubes $Q^{n}_{k,m} \subset Q^{n}_{j,\widetilde{m}}$, we have
\begin{equation}
\label{4.22}
\delta^{l}_{r}(Q^{n}_{k,m})N^{l}_{j,\widetilde{m}} \le C 2^{(j-k)l}
\end{equation}
with constant $C >0$ independent of $k,j,m,\widetilde{m}$.

Using this fact, as well as \eqref{2.15} instead of \eqref{4.16}, \eqref{4.19}, \eqref{4.19}, we proceed step by step in the same manner as in the proof of
Theorem~\ref{Th4.1}, replacing all the splines $U_{k}^{l-1}$ in the proof of Theorem~\ref{Th4.1} by the splines $U_{k}^{l}$.
In view of \eqref{4.22}, we clearly need not resort to deal with the sets $\Gamma^{1}_{j,\widetilde{m}}$ and $\Gamma^{2}_{j,\widetilde{m}}$,
and so estimate \eqref{4.20} is substantially simplified. Eventually, we obtain the conclusion of Theorem~\ref{Th4.1}.

\smallskip

The following result is a corollary to Theorem~\ref{Th4.1}.

\begin{Ca}
\label{Ca4.1}
Let $p,q,r \in (0,\infty]$, $p \neq \infty$,  $\varphi \in L_{r}^{\text{\rm loc}}(\mathbb{R}^{n})$, a~weight sequence $\{s_{k}\} \in ^{\text{\rm loc}}Y^{\alpha_{3}}_{\alpha_{1},\alpha_{2}}$,
and a~multiple sequence $\{s_{k,m}\}$
be $p$-associated with the weight sequence $\{s_{k}\}$. Next, let
$d \in \mathbb{N}_{0}$, a~weight $\gamma^{p} \in A^{\text{\rm loc}}_{\infty}(\mathbb{R}^{n+d})$ and
a~multiple sequence $\{\widehat{\gamma}_{k,m}\}$ be generated by the weight~$\gamma$. Let $\alpha_{2} < \lambda + \frac{d \delta_{1}(\gamma)}{p}$ for $\lambda:= \min
\{l,l-1+\frac{\delta_{2}(\gamma)}{p}\}$. We set $t_{k,m}:=s_{k,m}\widehat{\gamma}_{k,m}$ for $k \in \mathbb{N}_{0}$, $m \in \mathbb{Z}^{n}$.

Then a necessary a sufficient condition that a function $\varphi $ be in $\widetilde{B}_{p,q,r}^{l}(\mathbb{R}^{n},\{t_{k,m}\})$ is that
\begin{equation}
\label{4.23}
N_{1}(\varphi,l):=\left(\sum\limits_{j=-1}^{\infty}(s^{l}_{k}(\varphi,\{t_{k}\})_{r,p})^{q}\right)^{\frac{1}{q}}  < \infty.
\end{equation}
Moreover,
\begin{equation}
\begin{split}
\label{4.24}
&N_{1}(\varphi,l)\sim N_{2}(\varphi,l) \sim \|\varphi|\widetilde{B}_{p,q,r}^{l}(\mathbb{R}^{n},\{t_{k,m}\})\|, \qquad \mbox{ where }\\
&N_{2}(\varphi,l):=\left(\sum\limits_{k=0}^{\infty} \left(\sum\limits_{m \in \mathbb{Z}^{n}}t^{p}_{k,m}\|\varphi-T^{l-1}_{k}(\varphi,r)|L_{r}(Q^{n}_{k,m})\|^{p}\right)^{\frac{q}{p}}\right)^{\frac{1}{q}} + \left(\sum\limits_{m \in \mathbb{Z}^{n}}t^{p}_{0,m}\|\varphi|L_{r}(Q_{0,m}^{n})\|^{p}\right)^{\frac{1}{p}}.
\end{split}
\end{equation}
\end{Ca}

\textbf{Proof}. We consider the case $q < \infty$, because the case $q=\infty$ is dealt with similarly.
Let $\varphi \in \widetilde{B}_{p,q,r}^{l}(\mathbb{R}^{n},\{t_{k,m}\})$. The estimate
$N_{1}(\varphi,l) \le N_{2}(\varphi,l) \le C \|\varphi|\widetilde{B}_{p,q,r}^{l}(\mathbb{R}^{n},\{t_{k,m}\})\|$ follows from~\eqref{4.6}.

Assume now that $N_{1}(\varphi,l) < \infty$. We choose $\varepsilon \geq 0 $ so small that $\alpha_{2}-n < \widetilde{\lambda} + \frac{d \widetilde{\delta}_{1}}{p}$.
Applying Theorem \ref{Th4.1} with $\mu \le \min\{1,q,r\}$ and next using Theorem~\ref{Th2.2}, this gives
\begin{equation}
\begin{split}
\label{4.25}
&\|\varphi|\widetilde{B}_{p,q,r}^{l}(\mathbb{R}^{n},\{t_{k,m}\})\|^{q} \le\\
 &\le \sum\limits_{k=0}^{\infty}2^{-kq(\widetilde{\lambda}  + \frac{d \widetilde{\delta}_{1}}{p}-(\alpha_{2}-\frac{n}{p}))}\left(\sum\limits_{j=-1}^{k}2^{j\mu(\widetilde{\lambda}+\frac{d \widetilde{\delta}_{1}}{p}-(\alpha_{2}-\frac{n}{p}))}\left(s^{l}_{j}(\varphi,\{t_{k}\})_{r,p}\right)^{\mu}\right)^{\frac{q}{\mu}} \le C [N_{1}(\varphi,l)]^{q},
\end{split}
\end{equation}
proving the corollary.

In a similar manner Theorem~\ref{Th4.2} applies to obtain the following result.

\begin{Ca}
\label{Ca4.2}
Let $p,q,r \in (0,\infty]$,   $\alpha_{1},\alpha_{2} \in \mathbb{R}$, $\alpha_{3} \geq 0$, $\sigma_{1} \in (0,\infty]$,$\sigma_{2}=p$, a~multiple sequence $\{t_{k,m}\}$ be $p$-associated
 with a~$p$-admissible weight sequence $ \{t_{k}\} \in X^{\alpha_{3}}_{\alpha,\sigma,p}$. If $l > \alpha_{2}$, then a~function $\varphi \in \widetilde{B}_{p,q,r}^{l}(\mathbb{R}^{n},\{t_{k,m}\})$
if and only if
\begin{equation}
\label{4.26}
N_{1}(\varphi, l+1)  < \infty.
\end{equation}
Moreover,
\begin{equation}
\begin{split}
\label{4.27}
N_{1}(\varphi,l+1)\sim N_{2}(\varphi,l+1) \sim \|\varphi|\widetilde{B}_{p,q,r}^{l}(\mathbb{R}^{n},\{t_{k,m}\})\|.
\end{split}
\end{equation}
\end{Ca}

Theorem~\ref{Th4.2} can be used to obtain a~result on equivalent norms in the space $\widetilde{B}_{p,q,r}^{l}(\mathbb{R}^{n},\{t_{k}\})$  for different (sufficiently large)~$l$.

\begin{Ca}
\label{Ca4.3}
Let $p,q,r \in (0,\infty]$, $\alpha_{1},\alpha_{2} \in \mathbb{R}$, $\alpha_{3} \geq 0$, $\sigma_{1} \in (0,\infty]$, $\sigma_{2}=p$, $\{t_{k}\} \in X_{\alpha,\sigma,p}^{\alpha_{3}}$ and $l > \alpha_{2}$.
Then, for $l' > l$,
$$
\|\varphi|\widetilde{B}_{p,q,r}^{l'}(\mathbb{R}^{n},\{t_{k}\})\| \sim \|\varphi|\widetilde{B}_{p,q,r}^{l}(\mathbb{R}^{n},\{t_{k}\})\|.
$$
\end{Ca}

\textbf{Proof}. To obtain the the estimate
$$
\|\varphi|\widetilde{B}_{p,q,r}^{l}(\mathbb{R}^{n},\{t_{k}\})\| \le C N_{2}(\varphi,l') \le C \|\varphi|\widetilde{B}_{p,q,r}^{l'}(\mathbb{R}^{n},\{t_{k}\})\|
$$
it suffices to employ Corollary \ref{Ca4.2} and Lemma \ref{Lm4.1}. The reverse estimate is clear.

\begin{Remark}
\label{R4.2}
As was pointed out in the introduction, the methods of \cite{Be2} (\cite{KeVy}) enable us to obtain a~result similar to Corollary~\ref{Ca4.3}
for the space $\overline{B}^{l}_{p,q,1}(\mathbb{R}^{n},\{t_{k}\})$, provided that $\{t_{k}\} \in ^{\text{\rm loc}}Y^{\alpha_{3}}_{\alpha_{1},\alpha_{2}}$ ($\{t_{k}\} \in Y^{\alpha_{3}}_{\alpha_{1},\alpha_{2}}$) and $l > \alpha_{2}$.

The following simple example shows that Corollary \ref{Ca4.3} applies for weaker assumptions on the variable smoothness.

 Let $p \in (1,\infty)$. We claim that there exists a~weight sequence $\{\gamma^{1}_{k}\}$ such that
 $\{\gamma^{1}_{k}\} \in \widetilde{X}^{\alpha_{3}}_{\alpha,\sigma,p}$, $\alpha_{2} < l$, $\sigma_{2}=p$ and
 $\{\gamma^{1}_{k}\} \in ^{\text{\rm loc}}Y^{\alpha_{3}}_{\alpha'_{1},\alpha'_{2}}$ for $l < \alpha'_{2}$, but
 $\{\gamma^{1}_{k}\} \notin ^{\text{\rm loc}}Y^{\alpha_{3}}_{\alpha''_{1},\alpha''_{2}}$ for any  $\alpha''_{2} \le l$. Indeed, let
 $\varepsilon \in (0,n)$ be a~sufficiently small number that will be chosen later.
We define  $(\gamma^{1})^{p}(x,x_{n+1}):=\prod\limits_{i=1}^{n+1}\frac{1}{|x_{i}|^{1-\varepsilon}}$ for $(x,x_{n+1}) \in \mathbb{R}^{n+1} \setminus \{0\}$. Note that
$(\gamma^{1})^{p} \in A_{1}(\mathbb{R}^{n+1})$. Consider the multiple sequence
$(\gamma^{1}_{k,m})^{p}:=2^{klp}\iint\limits_{\Xi^{1,n}_{k,m}}(\gamma^{1})^{p}(x,x_{n+1})\,dx dx_{n+1}$ for $k \in \mathbb{N}_{0}$, $m \in \mathbb{Z}^{n}$. We set
$(\gamma^{1})^{p}(x):=\sum\limits_{m \in \mathbb{Z}^{n}}\chi_{\widetilde{Q}^{n}_{k,m}}(x)2^{kn}(\gamma^{1}_{k,m})^{p}$ for $x \in \mathbb{R}^{n}$, $k \in \mathbb{N}_{0}$.

Clearly, the inequality
$$
\frac{(\gamma^{1}_{k+1})^{p}(0)}{(\gamma^{1}_{k})^{p}(0)} = 2^{pl+n}\frac{\iint\limits_{\Xi^{1,n}_{k+1,0}}(\gamma^{1})^{p}(x,x_{n+1})\,dx dx_{n+1}}{\iint\limits_{\Xi^{1,n}_{k,0}}(\gamma^{1})^{p}(x,x_{n+1})\,dx dx_{n+1}} \geq 2^{pl+\frac{n}{2}}
$$
is satisfied for sufficiently small $\varepsilon \in (0,1)$. Hence, $\{\gamma^{1}_{k}\}\in Y^{\alpha_{3}}_{\alpha_{1},\alpha_{2}}$ only if $\alpha_{2} \geq l+\frac{n}{2p}>l$.

However, by Example \ref{Ex2.1} we have $\{\gamma^{1}_{k}\} \in \widetilde{X}^{\alpha_{3}}_{\alpha,\sigma,p}$ for $\sigma_{2}=p$, $\alpha_{2}=l-\frac{\delta_{1}(\gamma)}{p}<l$, and hence
all the hypotheses of Corollary \ref{Ca4.3} are satisfied.

We note that by Remark \ref{R2.3} the space $\widetilde{B}_{p,p,r}^{l}(\mathbb{R}^{n},\{\gamma^{1}_{k}\})$ is nontrivial for any $p \in (1,\infty)$, $r \in [1,p]$.
Moreover, the space $\widetilde{B}_{p,p,r}^{l}(\mathbb{R}^{n},\{\gamma^{1}_{k}\})$ (by Theorem~\ref{Th3.1}) is the trace of the Sobolev space
$W^{l}_{p}(\mathbb{R}^{n+1},\gamma)$.
\end{Remark}

The following theorem is an important step in the proof of the atomic decomposition theorem.
However, the estimate obtained here may be of independent interest.

Given $p \in (0,\infty]$, $\theta \in (0,p]$, we set $p_{\theta}:=\frac{p}{\theta}$, provided that $p$ and $\theta$ are not simultaneously infinite. For $\theta=p=\infty$ we set
$p_{\theta}:=1$. For $0<\theta<p=\infty$ we assume $p_{\theta}=\infty$.

\begin{Th}
\label{Th4.3}
Let $p,q,r \in (0,\infty]$, $\theta \in (0,\min\{p,r\}]$, be a $p$-admissible  weight sequence $\{t_{k}\} \in X^{\alpha_{3}}_{\alpha,\sigma,p}$ with $\sigma_{1}=\theta p'_{\theta}$,
$\alpha_{1}> n(\frac{1}{\theta}-\frac{1}{r})$, $\sigma_{2} \in (0,\infty]$, $\alpha_{3} \geq 0$, a~multiple sequence $\{t_{k,m}\}$ be a~$p$-associated with weight sequence
$\{t_{k}\}$. Assume that functions $V_{k} \in L^{\text{\rm loc}}_{r}(\mathbb{R}^{n})$ (with $k \in \mathbb{N}_{0}$) and
$$
\biggl(\sum\limits_{k=0}^{\infty}\biggl(\sum\limits_{m \in \mathbb{Z}^{n}}t^{p}_{k,m}2^{\frac{knp}{r}}\|v_{k}|L_{r}(Q^{n}_{k,m})\|^{p}\biggr)^{\frac{q}{p}}\biggr)^{\frac{1}{q}} < \infty,
$$
where $v_{k}:=V_{k}-V_{k-1}$($V_{-1} \equiv 0$) for $k \in \mathbb{N}_{0}$.

Then the series $\sum\limits_{k=0}^{\infty}v_{k}$ converges in $L^{\text{\rm loc}}_{r}(\mathbb{R}^{n})$ to some function  $\varphi \in L_{r}^{\text{\rm loc}}(\mathbb{R}^{n})$,
and moreover,
\begin{equation}
\begin{split}
\label{4.28}
&\biggl(\sum\limits_{k=0}^{\infty}\biggl(\sum\limits_{m \in \mathbb{Z}^{n}}t^{p}_{k,m}2^{\frac{knp}{r}}\|\varphi-V_{k}|L_{r}(Q^{n}_{k,m})\|^{p}\biggr)^{\frac{q}{p}}\biggr)^{\frac{1}{q}} +
\biggl(\sum\limits_{m \in \mathbb{Z}^{n}}t^{p}_{0,m}\|\varphi|L_{r}(Q^{n}_{0,m})\|^{p}\biggr)^{\frac{1}{p}}\le\\
&\le C \biggl(\sum\limits_{k=0}^{\infty}\biggl(\sum\limits_{m \in \mathbb{Z}^{n}}t^{p}_{k,m}2^{\frac{knp}{\theta}}\|v_{k}|L_{r}(Q^{n}_{k,m})\|^{p}\biggr)^{\frac{q}{p}}\biggr)^{\frac{1}{q}},
\end{split}
\end{equation}
in which the constant $C > 0$ depends on $\alpha_{1},\alpha,\sigma$, $r$,  but is independent of the function sequence $\{V_{k}\}$.
\end{Th}

\textbf{Proof}.We consider only the case $p,q \neq \infty$.

We claim that the series $\sum\limits_{k=0}^{\infty}v_{k}$ converges in $L_{r}^{\text{\rm loc}}(\mathbb{R}^{n})$ to some
function $\varphi \in L_{r}^{\text{\rm loc}}(\mathbb{R}^{n})$. Indeed,
it suffices to show that the series $\sum\limits_{k=0}^{\infty}v_{k}$ converges in $L_{r}(Q^{n}_{0,m})$ for any cube $Q^{n}_{0,m}$.

For any $j_{1} \le  j_{2} \in \mathbb{N}_{0}$ and $\mu \le \min\{1,\theta\}$ it follows from Lemma~\ref{Lm2.1} that
\begin{equation}
\begin{split}
\label{4.29}
&\|V_{j_{1}}-V_{j_{2}}|L_{r}(Q^{n}_{0,m})\| \le \biggl(\sum\limits_{j=j_{1}}^{\infty}\|v_{j}|L_{r}(Q^{n}_{0,m})\|^{\mu}\biggr)^{\frac{1}{\mu}} =\\
&= \biggl(\sum\limits_{j=j_{1}}^{\infty}\biggl(\sum_{\substack{\widetilde{m}
\in \mathbb{Z}^{n}\\Q^{n}_{j,\widetilde{m}}\subset
Q^{n}_{0,m}}}\|v_{j}|L_{r}(Q^{n}_{j,\widetilde{m}})\|^{\theta}\biggr)^{\frac{\mu}{\theta}}\biggr)^{\frac{1}{\mu}}=K_{j_{1},m}.
\end{split}
\end{equation}

Applying H\"older's inequality to the inner sum (in $\widetilde{m}$) with exponents $\tau$ and $\tau'$ and using \eqref{2.6}, we have, for any $\mu \le \min\{1,\theta\}$
\begin{equation}
\begin{split}
\label{4.30}
&\biggl(K_{j_{1},m}\biggr)^{\mu}=\sum\limits_{j=j_{1}}^{\infty}\biggl(\sum_{\substack{\widetilde{m}
\in \mathbb{Z}^{n}\\Q^{n}_{j,\widetilde{m}}\subset
Q^{n}_{0,m}}}\frac{t^{\theta}_{j,\widetilde{m}}}{t^{\theta}_{j,\widetilde{m}}}\|v_{j}|L_{r}(Q^{n}_{j,\widetilde{m}})\|^{\theta}\biggr)^{\frac{\mu}{\theta}} \le \\
& \le C \sum\limits_{j=j_{1}}^{\infty}\frac{1}{2^{\frac{jn\mu}{r}}}\frac{t^{\mu}_{0,m}}{t^{\mu}_{0,m}}\biggl(\sum_{\substack{\widetilde{m}
\in \mathbb{Z}^{n}\\Q^{n}_{j,\widetilde{m}}\subset
Q^{n}_{0,m}}}\biggl[\frac{1}{t_{j,\widetilde{m}}}\biggr]^{\theta p'_{\theta}}\biggr)^{\frac{\mu}{\theta p'_{\theta}}}\biggl(\sum_{\substack{\widetilde{m}
\in \mathbb{Z}^{n}\\Q^{n}_{j,\widetilde{m}}\subset
Q^{n}_{0,m}}}2^{\frac{jnp}{r}}t^{p}_{j,\widetilde{m}}\|v_{j}|L_{r}(Q^{n}_{j,\widetilde{m}})\|^{p}\biggr)^{\frac{\mu}{p}} \le \\
&\le C \sum\limits_{j=j_{1}}^{\infty}\frac{2^{j\mu n(\frac{1}{\theta}-\frac{1}{r})}}{2^{j\mu \alpha_{1}}}\biggl(\sum_{\substack{\widetilde{m}
\in \mathbb{Z}^{n}\\Q^{n}_{j,\widetilde{m}}\subset
Q^{n}_{0,m}}}2^{\frac{jnp}{r}}t^{p}_{j,\widetilde{m}}\|v_{j}|L_{r}(Q^{n}_{j,\widetilde{m}})\|^{p}\biggr)^{\frac{\mu}{p}}.
\end{split}
\end{equation}

We choose $\mu \le \min\{1,q,r\}$ and take $q_{\mu}:=\frac{q}{\mu} \geq 1$.
An application of H\"older's inequality with exponents $q_{\mu}$ and $q'_{\mu}$ to the right-hand side of~\eqref{4.30} shows that
\begin{equation}
\begin{split}
\label{4.31}
&K_{j_{1},m} \le C \biggl(\sum\limits_{j=j_{1}}^{\infty}\frac{1}{2^{j\mu q'_{\mu}(\alpha_{1}-(\frac{n}{\theta}-\frac{n}{r}))}}\biggr)^{\frac{1}{\mu q'_{\mu}}}\biggl(\sum\limits_{j=j_{1}}^{\infty}
\biggl(\sum_{\substack{\widetilde{m}
\in \mathbb{Z}^{n}\\Q^{n}_{j,\widetilde{m}}\subset
Q^{n}_{0,m}}}2^{\frac{jnp}{r}}t^{p}_{j,\widetilde{m}}\|v_{j}|L_{r}(Q^{n}_{j,\widetilde{m}})\|^{p}\biggr)^{\frac{q}{p}}\biggr)^{\frac{1}{q}}
\end{split}
\end{equation}

Note that for a fixed $j_{1}$ the right-hand side of inequalities \eqref{4.31} tends to infinity as $\alpha_{1}$ tends to  $\frac{n}{\theta}-\frac{n}{r}$.

From \eqref{4.29}, \eqref{4.31} and since the space
$L_{r}(Q^{n}_{0,m})$ is complete, we obtain the required convergence of the series
$\sum\limits_{k=0}^{\infty}v_{k}$ in $L_{r}(Q^{n}_{0,m})$ to some function $\varphi_{m} \in L_{r}(Q^{n}_{0,m})$.

Let us prove \eqref{4.28}. Applying Lemma \ref{2.1} with $f_{j}=0$ with  $j<k$ and $f_{j}:=v_{j}\chi_{Q^{n}_{k,m}}$ with $j\geq k$,
and then using Minkowski's inequality (because $\frac{p}{\mu} \geq 1$), this gives
\begin{equation}
\begin{split}
\label{4.32}
&\biggl(\sum\limits_{m \in
\mathbb{Z}^{n}}2^{\frac{nkp}{r}}t^{p}_{k,m}\|\varphi-V_{k}|L_{r}(Q^{n}_{k,m})\|^{p}\biggr)^{\frac{1}{p}}
\le 2^{\frac{nk}{r}}
\biggl(\sum\limits_{m \in
\mathbb{Z}^{n}}t^{p}_{k,m}\biggl[\sum\limits_{j=k}^{\infty}\|v_{j}|L_{r}(Q^{n}_{k,m})\|^{\mu}\biggr]^{\frac{p}{\mu}}\biggr)^{\frac{1}{p}}\le\\
&\le 2^{\frac{nk}{r}}
\biggl(\sum\limits_{j=k}^{\infty}\biggl(\sum\limits_{m \in
\mathbb{Z}^{n}}t^{p}_{k,m}\|v_{j}|L_{r}(Q^{n}_{k,m})\|^{p}\biggr)^{\frac{\mu}{p}}\biggr)^{\frac{1}{\mu}}\le\\
&\le 2^{\frac{nk}{r}} \biggl(\sum\limits_{j=k}^{\infty}\biggl(\sum\limits_{m
\in \mathbb{Z}^{n}}t^{p}_{k,m}\biggl[\sum_{\substack{\widetilde{m} \in
\mathbb{Z}^{n}\\Q^{n}_{j,\widetilde{m}}\subset
Q^{n}_{k,m}}}\|v_{j}|L_{r}(Q^{n}_{j,\widetilde{m}})\|^{\theta}\biggr]^{\frac{p}{\theta}}\biggr)^{\frac{\mu}{p}}\biggr)^{\frac{1}{\mu}}=:R_{k}.
\end{split}
\end{equation}

Arguing as in the proof of estimate \eqref{4.30}, we obtain for $j \geq k+1$
\begin{equation}
\begin{split}
\label{4.33}
&2^{\frac{nkp}{r}}\sum\limits_{m \in \mathbb{Z}^{n}}t^{p}_{k,m}\biggl[\sum_{\substack{\widetilde{m} \in
\mathbb{Z}^{n}\\Q^{n}_{j,\widetilde{m}}\subset
Q^{n}_{k,m}}}\|v_{j}|L_{r}(Q^{n}_{j,\widetilde{m}})\|^{\theta}\biggr]^{\frac{p}{\theta}}
\le\\
&\le 2^{\frac{nkp}{r}} \sum\limits_{m \in \mathbb{Z}^{n}}t^{p}_{k,m}\biggl(\sum_{\substack{\widetilde{m}
\in \mathbb{Z}^{n}\\Q^{n}_{j,\widetilde{m}}\subset
Q^{n}_{k,m}}}\biggl[\frac{1}{t_{j,\widetilde{m}}}\biggr]^{\theta p'_{\theta}}\biggr)^{\frac{p}{\theta p'_{\theta}}}\biggl[\sum_{\substack{\widetilde{m} \in
\mathbb{Z}^{n}\\Q^{n}_{j,\widetilde{m}}\subset
Q^{n}_{k,m}}}t^{p}_{j,\widetilde{m}}\|v_{j}|L_{r}(Q^{n}_{j,\widetilde{m}})\|^{p}\biggr] \le \\
&\le C 2^{(k-j)p(\alpha_{1}-(\frac{n}{\theta}-\frac{n}{r}))}\sum\limits_{m \in \mathbb{Z}^{n}}t^{p}_{j,m}2^{\frac{njp}{r}}\|v_{j}|L_{r}(Q^{n}_{j,\widetilde{m}})\|^{p}
\end{split}
\end{equation}

We take $\mu \le \min\{1,\theta,q\}$. Since $\alpha_{1} > \frac{n}{\theta}-\frac{n}{r}$ by
the hypothesis of the lemma, it follows by Hardy's inequality and \eqref{4.32}, \eqref{4.33} that
\begin{equation}
\begin{split}
\label{4.34}
&\biggl(\sum\limits_{k=0}^{\infty}R^{q}_{k}\biggr)^{\frac{1}{q}} \le C \biggl(\sum\limits_{k=0}^{\infty}\biggl(\sum\limits_{m \in
\mathbb{Z}^{n}}t^{p}_{j,m}2^{\frac{jnp}{r}}\|v_{j}|L_{r}(Q^{n}_{j,m})\|^{p}\biggr)^{\frac{q}{p}}\biggr)^{\frac{1}{q}}.
\end{split}
\end{equation}

Now the required estimate for the first term on the left of \eqref{4.27} follows from \eqref{4.30} and~\eqref{4.33}.

Let us estimate the second term in the left-hand side of~\eqref{4.28}. Similarly to \eqref{4.29}, $\|\varphi|L_{r}(Q^{n}_{0,m})\| \le C K_{1,m}$. Hence,
using \eqref{4.31} we easily obtain the estimate
\begin{equation}
\begin{split}
\label{4.36}
&\biggl(\sum\limits_{m \in \mathbb{Z}^{n}}t^{p}_{0,m}\|\varphi|L_{r}(Q^{n}_{0,m})\|^{p}\biggr)^{\frac{1}{p}}
\le C \biggl[ \sum\limits_{j=0}^{\infty}\biggl(\sum\limits_{m \in \mathbb{Z}^{n}}t^{p}_{j,m}2^{\frac{jnp}{r}}\|v_{j}|L_{r}(Q^{n}_{j,m})\|^{p}\biggr)^{\frac{q}{p}} \biggr]^{\frac{1}{q}}.
\end{split}
\end{equation}
The proof of Theorem~\ref{Th4.3} is complete.

For $\varphi \in L_{r}^{\text{\rm loc}}(\mathbb{R}^{n})$ we assume that
$$
\varphi=\sum\limits_{k=0}^{\infty}v^{l}_{k}  \mbox{ in } L^{\text{\rm loc}}_{r}(\mathbb{R}^{n}), \mbox{ where }
v^{l}_{k}(x):=\sum\limits_{m \in \mathbb{Z}^{n}}\beta_{k,m} N^{l}_{k,m}(x), \quad k \in \mathbb{N}_{0}, x \in \mathbb{R}^{n}.
$$

Given $p,q,r \in (0,\infty]$, we define (with corresponding  modifications for $p=\infty$ or $q=\infty$)
\begin{equation}
\label{4.37}
N_{3}(\varphi,l+1):=\inf \biggl(\sum\limits_{k=0}^{\infty}\biggl(\sum\limits_{m \in \mathbb{Z}^{n}}t^{p}_{k,m}|\beta_{k,m}|^{p}\biggr)^{\frac{q}{p}}\biggr)^{\frac{1}{q}},
\end{equation}
where the infimum in \eqref{4.37} is taken over all series $\sum\limits_{k=0}^{\infty}v^{l}_{k}$ convening in $L_{r}^{\text{\rm loc}}(\mathbb{R}^{n})$ to the function~$\varphi$.

For $\varphi \in L_{r}^{\text{\rm loc}}(\mathbb{R}^{n})$ we also set  (with corresponding modifications in the case $p=\infty$ or $q=\infty$)
$$
N_{4}(\varphi,l+1):=\biggl(\sum\limits_{k=0}^{\infty}\biggl(\sum\limits_{m \in \mathbb{Z}^{n}}t^{p}_{k,m}|\alpha_{k,m}(T^{l}_{k}(\varphi,r))|^{p}\biggr)^{\frac{q}{p}}\biggr)^{\frac{1}{q}}.
$$

The next result extends Theorem 5.1 of \cite{DeVore} to the case of Besov spaces of variable smoothness $\widetilde{B}_{p,q,r}^{l}(\mathbb{R}^{n},\{t_{k}\})$.

\begin{Ca} (\textbf{the atomic decomposition})
\label{Ca4.4}
Let $p,q,r \in  (0,\infty]$, $\theta \in (0,\min\{r,p\}]$. Next, $\{t_{k}\} \in X^{\alpha_{3}}_{\alpha,\sigma,p}$
be a~$p$-admissible  weight sequence with
$\alpha_{1}>n(\frac{1}{\theta}-\frac{1}{r})$, $l > \alpha_{2}$, $\sigma_{1}=\theta p'_{\theta}$, $\sigma_{2}=p$, a~multiple sequence $\{t_{k,m}\}$ be $p$-associated with the weight sequence $\{t_{k}\}$.  Then

{\rm 1)} each function $\varphi \in \widetilde{B}_{p,q,r}^{l}(\mathbb{R}^{n},\{t_{k}\})$
may be expanded into an $L_{r}^{\text{\rm loc}}(\mathbb{R}^{n})$-convergent series of splines $N^{l}_{k,m}$; that is,
\begin{equation}
\begin{split}
\label{4.38}
&\varphi=\sum\limits_{k=0}^{\infty}v^{l}_{k}(\varphi) \mbox{ in the sense of  } L_{r}^{\text{\rm loc}}({R}^{n}), \mbox{ where }\\
&v^{l}_{k}(\varphi)(x)=\sum\limits_{m \in \mathbb{Z}^{n}}\beta_{k,m}(\varphi)N^{l}_{k,m}(x) \mbox{ for } x \in \mathbb{R}^{n}.
\end{split}
\end{equation}
Moreover, for some constant $C > 0$
$$
N_{3}(\varphi,l+1) \le  N_{4}(\varphi,l+1) \le C\|\varphi|\widetilde{B}_{p,q,r}^{l}(\mathbb{R}^{n},\{t_{k}\})\|;
$$

{\rm 2)} if, for some multiple sequence $\{\beta_{k,m}\}$,
$$
\biggl(\sum\limits_{k=0}^{\infty}\biggl(\sum\limits_{m \in \mathbb{Z}^{n}}t^{p}_{k,m}|\beta_{k,m}|^{p}\biggr)^{\frac{q}{p}}\biggr)^{\frac{1}{q}} < \infty,
$$
then the series $\sum\limits_{k=0}^{\infty}\sum\limits_{m \in \mathbb{Z}^{n}}\beta_{k,m}N^{l}_{k,m}$ converges in $L_{r}^{\text{\rm loc}}(\mathbb{R}^{n})$ to some function
$\varphi \in \widetilde{B}_{p,q,r}^{l}(\mathbb{R}^{n},\{t_{k}\})$ and there exist constants $C_{1}, C_{2} > 0$ such that
$$
\|\varphi|\widetilde{B}_{p,q,r}^{l}(\mathbb{R}^{n},\{t_{k}\})\| \le C_{1} N_{3}(\varphi,l+1) \le C_{2} N_{4}(\varphi,l+1).
$$
\end{Ca}

The proof is similar to those of Theorem 5.1 and Corollary 5.3 of \cite{DeVore}, one only needs to appropriately use Theorem~\ref{Th4.3}, Corollary \ref{Ca4.2}, Lemma~\ref{Lm4.2} and
the clear estimate
\begin{equation}
\label{4.39}
\biggl(\sum\limits_{\substack{\widetilde{m} \in \mathbb{Z}^{n}\\
Q^{n}_{k,m} \bigcap \operatorname{supp} N^{l-1}_{k,\widetilde{m}} \neq
\emptyset}}|\beta_{k,\widetilde{m}}|^{p} \biggr)^{\frac{1}{p}} \le \biggl(\sum\limits_{\substack{\widetilde{m} \in \mathbb{Z}^{n}\\
Q^{n}_{k,m} \bigcap \operatorname{supp} N^{l-1}_{k,\widetilde{m}} \neq
\emptyset}}|\beta_{k,\widetilde{m}}|^{r}\biggr)^{\frac{1}{r}} \le C \biggl(\sum\limits_{\substack{\widetilde{m} \in \mathbb{Z}^{n}\\
Q^{n}_{k,m} \bigcap \operatorname{supp} N^{l-1}_{k,\widetilde{m}} \neq
\emptyset}}|\beta_{k,\widetilde{m}}|^{p}\biggr)^{\frac{1}{p}},
\end{equation}
the modifications in \eqref{4.39} for $ p=\infty$ and $r =\infty$ are straightforward.

The constant $C$ in \eqref{4.39} depends only on $n,l,p,r$.

\begin{Remark}
\label{R4.3}
Under the hypotheses of Corollary \ref{Ca4.4}  the set $\Sigma^{l}$ is dense in the space $\widetilde{B}_{p,q,r}^{l}(\mathbb{R}^{n},\{t_{k}\})$, $p,q \in (0,\infty)$.
Indeed, let $\varphi \in \widetilde{B}_{p,q,r}^{l}(\mathbb{R}^{n},\{t_{k}\})$. Then, by Corollary \ref{Ca4.4},
$$
\varphi=\sum\limits_{k=0}^{\infty}\sum\limits_{m \in \mathbb{Z}^{n}}\alpha_{k,m}(\varphi)N^{l}_{k,m} \hbox{ in the sense of } L^{\text{\rm loc}}_{r}(\mathbb{R}^{n})
$$
and, for any $\varepsilon > 0$,
$$
\biggl(\sum\limits_{k=0}^{\infty}\biggl(\sum\limits_{m \in \mathbb{Z}^{n}}t^{p}_{k,m}|\alpha_{k,m}(\varphi)|^{p}\biggr)^{\frac{q}{p}}\biggr)^{\frac{1}{q}} \le (1+\varepsilon)\inf \biggl(\sum\limits_{k=0}^{\infty}\biggl(\sum\limits_{m \in \mathbb{Z}^{n}}t^{p}_{k,m}|\beta_{k,m}|^{p}\biggr)^{\frac{q}{p}}\biggr)^{\frac{1}{q}} \le
$$
$$
\le C \|\varphi|\widetilde{B}_{p,q,r}^{l}(\mathbb{R}^{n},\{t_{k}\})\|.
$$
Hence, taking $\varphi_{n}:=\sum\limits_{k=0}^{n}\sum\limits_{m \in \mathbb{Z}^{n}}\alpha_{k,m}(\varphi)N^{l}_{k,m}$, it follows from Corollary \ref{Ca4.4} that
$$
\|\varphi_{n}-\varphi|\widetilde{B}_{p,q,r}^{l}(\mathbb{R}^{n},\{t_{k}\})\| \le C \biggl(\sum\limits_{k=n+1}^{\infty}\biggl(\sum\limits_{m \in \mathbb{Z}^{n}}t^{p}_{k,m}|\beta_{k,m}|^{p}\biggr)^{\frac{q}{p}}\biggr)^{\frac{1}{q}} \to 0, n \to \infty.
$$
\end{Remark}

\begin{Remark}
\label{R4.4}
 Let $p \in (1,\infty)$. We claim that one may choose parameters $r \in (1,p)$, $\alpha_{3}$ , $\alpha$, $\sigma$ and a~weight sequence
 $\{\gamma^{2}_{k}\} \in \widetilde{X}^{\alpha_{3}}_{\alpha,\sigma,p}$ so as to satisfy all the hypotheses of Corollary \ref{Ca4.4}. Besides,
 $\{\gamma^{2}_{k}\} \in Y^{\alpha_{3}}_{\alpha'_{1},\alpha'_{2}}$ for some $\alpha'_{1}<0$, but $\{\gamma^{2}_{k}\} \notin Y^{\alpha_{3}}_{\alpha''_{1},\alpha'_{2}}$
 for any $\alpha''_{1} \geq 0$.

Indeed, let $\varepsilon \in (0,p-1)$ be a sufficiently small number, which will be specified later. We set $(\gamma^{2})^{p}(x_{1},\dots,x_{n+1}):=\prod\limits_{i=1}^{n+1}|x_{i}|^{p-1-\varepsilon}$.
Note that $(\gamma^{2})^{p} \in A_{\frac{p}{\theta}}(\mathbb{R}^{n+1})$ for some $\theta \in (1,p)$. For $l \in \mathbb{N}$ consider the multiple sequence $(\gamma^{2}_{k,m})^{p}:=2^{klp}\iint\limits_{\Xi^{1,n}_{k,m}}(\gamma^{2})^{p}(x,x_{n+1})\,dxdx_{n+1}$ for $k \in \mathbb{N}_{0}$, $m \in \mathbb{Z}^{n}$. Let $(\gamma^{2}_{k})^{p}(x)=2^{nk}\sum\limits_{m \in \mathbb{Z}^{n}}\chi_{\widetilde{Q}^{n}_{k,m}}(x)(\gamma^{2}_{k,m})^{p}$ for $k \in \mathbb{N}_{0}$, $x \in \mathbb{R}^{n}$.

The inequality
$$
\frac{(\gamma^{2}_{k+1})^{p}(0)}{(\gamma^{2}_{k})^{p}(0)} = 2^{pl+n}\frac{\iint\limits_{\Xi^{1,n}_{k+1,0}}(\gamma^{2})^{p}(x,x_{n+1})\,dxdx_{n+1} dx_{n+1}}{\iint\limits_{\Xi^{1,n}_{k,0}}(\gamma^{2})^{p}(x,x_{n+1})\,dx dx_{n+1}}= C(p) 2^{pl+n - (p-\varepsilon)(n+1)}
$$
is clear.

If $p>n$, $l \le n$, then one easily checks that $\operatorname{sup}\{\alpha_{1}|\{\gamma^{2}_{k}\} \in Y^{\alpha_{3}}_{\alpha_{1},\alpha_{2}}\} < 0$ for
 sufficiently small $\varepsilon > 0$. On the other hand, from \ref{Ex2.1} it easily follows that the weight sequence
 $\{\gamma^{2}_{k}\} \in \widetilde{X}^{\alpha_{3}}_{\alpha,\sigma,p}$ with  $\sigma_{1}=p \frac{p'_{0}}{p_{0}}$ ($p_{0}=\frac{p}{\theta}$), $\sigma_{2}=p$, $\alpha_{1}=l+\frac{n}{p}+\frac{n}{\sigma_{1}}-\frac{(n+1)p_{0}}{p}=l-\frac{1}{\theta}> 0$ (because in our setting  $d=1$, $p_{0}=\frac{p}{\theta}$,  $\theta>1$) $\alpha_{2}=l-\frac{\delta_{1}(\gamma^{2},n,1)}{p} < l$. Hence, all the hypotheses of Corollary~\ref{Ca4.4} are satisfied.
\end{Remark}

\begin{Remark}
\label{R4.5}
Let $p,q \in (0,\infty]$, $p \neq \infty$, $r \in (0,p]$, $\gamma^{p} \in A^{\text{\rm loc}}_{\frac{p}{r}}(\mathbb{R}^{n})$, $s > 0$.
We set $t_{k,m}=2^{ks}\int\limits_{Q^{n}_{k,m}}\gamma^{p}(x)\,dx$ for $k \in \mathbb{N}_{0}, m \in \mathbb{Z}^{n}$.
Arguing as in Example 2.1, we conclude that the sequence $\{t_{k,m}\}$ satisfies the hypotheses of
Corollary~\ref{Ca4.4}. Hence, an application of Theorem \ref{Th2.5} gives the atomic decomposition theorem for the weighted Besov space
$\widetilde{B}^{s}_{p,q,r}(\mathbb{R}^{n},\gamma)=\overline{B}^{s}_{p,q,r}(\mathbb{R}^{n},\gamma)$ as a~particular case of Corollary~\ref{Ca4.4}.
Problems on atomic decomposition of the spaces
$B^{s}_{p,q}(\mathbb{R}^{n},\gamma)$  (and their generalizations) were studied in the papers \cite{HN}, \cite{Sawano}
by different methods (see also the references given therein).
 \end{Remark}

\textbf{Proof of Theorem~\ref{Th2.5}}.

\textit{Step 1}. We prove the embedding
$\overline{B}^{l}_{p,q,r}(\mathbb{R}^{n},\{t_{k}\}) \subset \widetilde{B}^{l}_{p,q,r}(\mathbb{R}^{n},\{t_{k}\})$ for any $0<r \le \frac{p}{p_{0}}$.
By the definition of the class $A_{p_{0}}^{\text{\rm loc}}(\mathbb{R}^{n})$ and from the conditions $r \le \frac{p}{p_{0}}$, we have
$\gamma \in A^{\text{\rm loc}}_{{p}/{r}}(\mathbb{R}^{n})$. Hence, using H\"older's inequality and properties of the sequence $\{s_{k}\}$,
$$
\int\limits_{\mathbb{R}^{n}}\gamma^{p}(x)s^{p}_{k}(x)[\delta^{l}_{r}(x+\frac{I^{n}}{2^{k}})\varphi]^{p}\,dx
\le C_{1} \sum\limits_{m \in \mathbb{Z}^{n}}\int\limits_{Q^{n}_{k,m}}\gamma^{p}(x)s^{p}_{k}(x)[\delta^{l}_{r}(Q^{n}_{k,m})\varphi]^{p}\,dx \le
$$
\begin{equation}
\label{4.40}
\le C_{2} \sum\limits_{m \in \mathbb{Z}^{n}}\int\limits_{Q^{n}_{k,m}}\gamma^{p}(x)s^{p}_{k}(x)\biggl(\int\limits_{Q^{n}_{k,m}}2^{2kn}\frac{\gamma^{r}(y)}{\gamma^{r}(y)}\int\limits_{\frac{I^{n}}{2^{k}}}|\Delta^{l}(h)\varphi(y)|^{r}\,dhdy\biggr)^{\frac{p}{r}}\,dx \le C_{3}\int\limits_{\mathbb{R}^{n}}t^{p}_{k}(y)[\overline{\Delta}^{l}_{r}(2^{-k})\varphi(y)]^{p}\,dy.
\end{equation}

The required embedding follows from estimate \eqref{4.40}.

\textit{Step 2}. Using the arguments employed in Example  \ref{Ex2.1}, we conclude that the sequence $\{t_{k}\}$ satisfies the hypotheses of
Corollary~\ref{Ca4.4} with $0 < r \le \frac{p}{p_{0}}$. Hence,
$\widetilde{B}^{l}_{p,q,r_{1}}(\mathbb{R}^{n},\{t_{k}\})=\widetilde{B}^{l}_{p,q,r_{2}}(\mathbb{R}^{n},\{t_{k}\})$, the norms being equivalent
for $0 < r_{1} \le r_{2} \le \frac{p}{p_{0}}$. So the theorem will be proved if we check the embedding $\widetilde{B}^{l}_{p,q,r}(\mathbb{R}^{n},\{t_{k}\}) \subset \overline{B}^{l}_{p,q,r}(\mathbb{R}^{n},\{t_{k}\})$ for any $0<r \le \frac{p}{p_{0}}$.

Let $\overline{t}^{p}_{k}(x):=\sum\limits_{m \in \mathbb{Z}^{n}}\chi_{Q^{n}_{k,m}}(x)2^{kn}\|t_{k}|L_{p}(Q^{n}_{k,m})\|^{p}:=\sum\limits_{m \in \mathbb{Z}^{n}}\chi_{Q^{n}_{k,m}}(x)2^{kn}t_{k,m}$ for $k \in \mathbb{N}_{0}$, $x \in \mathbb{R}^{n}$.

Let $\varphi \in \widetilde{B}^{l}_{p,q,r}(\mathbb{R}^{n},\{t_{k}\})$, then by Corollary \ref{Ca4.4} $\varphi=\sum\limits_{j=0}^{\infty}v^{l}_{j}(\varphi)$ in the sense of
$L_{r}^{\text{\rm loc}}({R}^{n})$, where $v^{l}_{j}(x)=\sum\limits_{m \in \mathbb{Z}^{n}}\beta_{j,m}N^{l}_{j,m}(x)$ for $x \in  \mathbb{R}^{n}$, $j \in \mathbb{N}_{0}$.

We set $\varphi_{1,k}:=\sum\limits_{i=0}^{k}v^{l}_{i}$, $\varphi_{2,k}:=\varphi - \varphi_{1}$. Clearly, for $k \in \mathbb{N}_{0}$,
\begin{equation}
\label{4.41}
\frac{1}{C}\|t_{k}\overline{\Delta}^{l}_{r}(2^{-k})\varphi|L_{p}(\mathbb{R}^{n})\| \le \|t_{k}\overline{\Delta}^{l}_{r}(2^{-k})\varphi_{1,k}|L_{p}(\mathbb{R}^{n})\|+\|t_{k}\overline{\Delta}^{l}_{r}(2^{-k})\varphi_{2,k}|L_{p}(\mathbb{R}^{n})\|=:S_{1,k}+S_{2,k}.
\end{equation}

Arguing as in the proof of Theorem~\ref{Th4.1} and using Example 2.1 and estimate \eqref{4.39}, we obtain, for $\mu \le \min\{1,r,q\}$,
\begin{equation}
\begin{split}
\label{4.42}
&S_{1,k} \le C \biggl(\sum\limits_{j=0}^{k}2^{(l-\alpha_{2})p(j-k)}\biggl(\sum\limits_{m \in \mathbb{Z}^{n}}t^{p}_{j,m}2^{\frac{jnp}{r}}\|v^{l}_{j}|L_{r}(Q^{n}_{j,m})\|^{p}\biggr)^{\frac{\mu}{p}}\biggr)^{\frac{1}{\mu}} \le \\
&\le C \biggl(\sum\limits_{j=0}^{k}2^{(l-\alpha_{2})p(j-k)}\biggl(\sum\limits_{m \in \mathbb{Z}^{n}}t^{p}_{j,m}|\beta_{j,m}|^{p}\biggr)^{\frac{\mu}{p}}\biggr)^{\frac{1}{\mu}}.
\end{split}
\end{equation}

To estimate $S_{2,k}$ we first note that
\begin{equation}
\label{4.43}
\begin{split}
&\biggl[\overline{\Delta}^{l}_{r}(2^{-k})v^{l}_{j}(x)\biggr]^{p} \le C \biggl(|v^{l}_{j}(x)|^{r}+\sum\limits_{i=1}^{l}\biggl(2^{kn}\int\limits_{\frac{I^{n}}{2^{k}}}|v^{l}_{j}(x+ih)|^{r}\,dh\biggr)\biggr)^{\frac{p}{r}}  \le \\
&\le \biggl(|v^{l}_{j}(x)|^{r}+\int\limits_{2(l+1)Q^{n}_{k,m}}|v^{l}_{j}(y)|^{r}\,dy\biggr)^{\frac{p}{r}}.
\end{split}
\end{equation}

Next, employing property 2) from \eqref{1.1}, (1.4),
 the (local) doubling property of the weight $\gamma^{p}$,
and using properties of the functions $N^{l}_{k,m}$, we have, for $\mu \le \min\{1,r,q\}$, 
$$
\int\limits_{\mathbb{R}^{n}}\gamma^{p}(x)s^{p}_{k}(x)|v^{l}_{j}(x)|^{p}\,dx \le C 2^{(k-j)p\alpha_{1}}\int\limits_{\mathbb{R}^{n}}\gamma^{p}(x)s^{p}_{j}(x)|v^{l}_{j}(x)|^{p}\,dx  \le \\
$$
$$
 \le C 2^{(k-j)p\alpha_{1}} \sum\limits_{m \in \mathbb{Z}^{n}}\int\limits_{Q^{n}_{k,m}}\gamma^{p}(x)s^{p}_{j}(x)\biggl(\sum\limits_{\substack{\widetilde{m} \in \mathbb{Z}^{n}\\
x \in \operatorname{supp}N^{l}_{j,\widetilde{m}}}}\beta^{p}_{j,\widetilde{m}}\biggr) \le C 2^{(k-j)p\alpha_{1}} \sum\limits_{\widetilde{m} \in \mathbb{Z}^{n}}2^{kn}\widehat{\gamma}^{p}_{j,\widetilde{m}}s^{p}_{j,\widetilde{m}}|\beta_{j,\widetilde{m}}|^{p} \le
$$
\begin{equation}
\label{4.44}
\le C 2^{(k-j)p\alpha_{1}}\sum\limits_{m \in \mathbb{Z}^{n}}t^{p}_{j,m}|\beta_{j,m}|^{p}.
\end{equation}

Using Example \ref{Ex2.1}, estimate \eqref{2.16}, property 1) of \eqref{1.1}, Lemma \ref{Lm4.2}, we have, for $j \geq k$,
\begin{equation}
\begin{split}
\label{4.45}
&\sum\limits_{m \in \mathbb{Z}^{n}}\int\limits_{Q^{n}_{k,m}}t_{k}^{p}(x)2^{\frac{knp}{r}}\|v^{l}_{j}|L_{r}(2(1+l)Q^{n}_{k,m})\|^{p}\,dx \le C \sum\limits_{m \in \mathbb{Z}^{n}}t^{p}_{k,m}2^{\frac{knp}{r}}\|v^{l}_{j}|L_{r}(Q^{n}_{k,m})\|^{p} \le \\
&\le C 2^{(k-j)p\alpha_{1}}\sum\limits_{m \in \mathbb{Z}^{n}}t^{p}_{j,m}2^{\frac{jnp}{r}}\|v^{l}_{j}|L_{r}(Q^{n}_{j,m})\|^{p} \le C 2^{(k-j)p\alpha_{1}}\sum\limits_{m \in \mathbb{Z}^{n}}t^{p}_{j,m}|\beta_{j,m}|^{p}.
\end{split}
\end{equation}

Combining estimates \eqref{4.42}, \eqref{4.43}, \eqref{4.44}, \eqref{4.45}, and arguing as in the proof of \eqref{4.32},
\begin{equation}
\begin{split}
\label{4.46}
S_{2,k} \le C 2^{k \alpha_{1}}\left(\sum\limits_{j=k}^{\infty}2^{-j\mu \alpha_{1}}\left(\sum\limits_{m \in \mathbb{Z}^{n}}t^{p}_{j,m}|\beta_{j,m}|^{p}\right)^{\frac{\mu}{p}}\right).
\end{split}
\end{equation}

Substituting estimates \eqref{4.42}, \eqref{4.46} in \eqref{4.41} and taking into account Hardy's inequality, we obtain
\begin{equation}
\begin{split}
\label{4.47}
\biggl(\sum\limits_{k=1}^{\infty}\|t_{k}\overline{\Delta}^{l}_{r}(2^{-k})\varphi|L_{p}(\mathbb{R}^{n})\|^{q}\biggr)^{\frac{1}{q}} \le  C \biggl(\biggl(\sum\limits_{m \in \mathbb{Z}^{n}}t^{p}_{j,m}|\beta_{j,m}|^{p}\biggr)^{\frac{q}{\frac{p}{}}}\biggr)^{\frac{1}{q}}.
\end{split}
\end{equation}

Using estimate \eqref{4.47} in combination with Corollary~\ref{Ca4.4} completes the proof of Theorem~\ref{Th2.5}.

\section{Embedding theorems for the spaces $\widetilde{B}^{l}_{p,q,r}(\mathbb{R}^{n},\{t_{k}\})$}

Let $\beta=\{\beta_{j}\}_{j=1}^{\infty}$ be a~sequence of nonnegative numbers, $w=\{w_{j,m}\}_{j \in \mathbb{N}, m \in \mathbb{Z}^{n}}$ be
a~multiple sequence of nonnegative numbers.

For $0<p,q\le \infty$, we set (with corresponding modifications in the case $p,q=\infty$)
\begin{equation}
\begin{split}
\label{5.1}
&l_{q}(\beta l_{p}(w)):=\{a=a_{j,m}: a_{j,m} \in \mathbb{R}, \|a|l_{q}(\beta l_{p}(w))\| < \infty\}, \mbox{ where } \\
&\|a|l_{q}(\beta l_{p}(w))\|=\left(\beta^{q}_{j}\left(\sum\limits_{m \in \mathbb{Z}^{n}}w^{p}_{j,m}|a_{j,m}|^{p}\right)^{\frac{q}{p}}\right)^{\frac{1}{q}}.
\end{split}
\end{equation}

\begin{Th}(\cite{KLS})
\label{Th5.1}
Let $0<p_{i},q_{i}\le \infty$ for $i=1,2$.

{\rm 1)} The space $l_{q}(\beta^{1} l_{p}(w^{1}))$ is continuously embedded into $l_{q}(\beta^{2} l_{p}(w^{2}))$ if and only if
\begin{equation}
\begin{split}
\label{5.2}
&\sum\limits_{j=1}^{\infty}\biggl(\frac{\beta^{2}_{j}}{\beta^{1}_{j}}\biggr)^{q^{*}}\biggl(\sum\limits_{m \in \mathbb{Z}^{n}}\biggl(\frac{w^{2}_{j,m}}{w^{1}_{j,m}}\biggr)^{p^{*}}\biggr)^{\frac{q^{*}}{p^{*}}} < \infty, \mbox{ where }\\
&\frac{1}{p^{*}}:=\max\{0,\frac{1}{p_{2}}-\frac{1}{p_{1}}\}, \frac{1}{q^{*}}:=\max\{0,\frac{1}{q_{2}}-\frac{1}{q_{1}}\}
\end{split}
\end{equation}

{\rm 2)} The space $l_{q}(\beta^{1} l_{p}(w^{1}))$ is compactly embedded into $l_{q}(\beta^{2} l_{p}(w^{2}))$ if and only if
condition \eqref{5.2} is satisfied, and moreover,
\begin{equation}
\label{5.3}
\lim\limits_{j \to \infty}\frac{\beta^{2}_{j}}{\beta^{1}_{j}}\biggl(\sum\limits_{m \in \mathbb{Z}^{n}}\biggl(\frac{w^{2}_{j,m}}{w^{1}_{j,m}}\biggr)^{p^{*}}\biggr)^{\frac{1}{p^{*}}}=0 \ \ \mbox{ if } q^{*}=\infty
\end{equation}
and
\begin{equation}
\label{5.4}
\lim\limits_{|m| \to \infty}\frac{w^{1}_{j,m}}{w^{2}_{j,m}}=\infty \ \ \mbox{ for all } \ j \in \mathbb{N} \ \ \mbox{ if } p^{*}=\infty.
\end{equation}
\end{Th}

As a direct corollary to Theorem~\ref{Th5.1} and Corollary~\ref{Ca4.4} we obtain

\begin{Ca}
Let  $i=1,2$ and let $0 < p^{i},q^{i},r^{i} \le \infty$, $\theta^{i} \in (0,\min\{p^{i},r^{i}\}]$, $p^{i}_{\theta^{i}}=\frac{p^{i}}{\theta^{i}}$.
Next, for $i=1,2$, let $\{t^{i}_{k}\} \in X^{\alpha^{i}_{3}}_{\alpha^{i},\sigma^{i},p^{i}}$ be a $p$-admissible  weight sequence, $\alpha^{i}_{1} > n(\frac{1}{\theta^{i}}-\frac{1}{r^{i}})$, $\sigma^{i}_{1}=r^{i}(p^{i}_{\theta^{i}})'$, $\sigma^{i}_{2}=p^{i}$, $l > \alpha^{i}_{2}$. Then

{\rm 1)} the space $\widetilde{B}^{l}_{p_{1},q_{1},r_{1}}(\mathbb{R}^{n},\{t^{1}_{k}\})$ is continuously embedded into
$\widetilde{B}^{l}_{p_{2},q_{2},r_{2}}(\mathbb{R}^{n},\{t^{2}_{k}\})$ if
\begin{equation}
\begin{split}
\label{5.6}
&\sum\limits_{j=0}^{\infty}\biggl(\sum\limits_{m \in \mathbb{Z}^{n}}\biggl(\frac{t^{2}_{j,m}}{t^{1}_{j,m}}\biggr)^{p^{*}}\biggr)^{\frac{q^{*}}{p^{*}}} < \infty, \mbox{ where }\\
&\frac{1}{p^{*}}:=\max\{0,\frac{1}{p_{2}}-\frac{1}{p_{1}}\}, \frac{1}{q^{*}}:=\max\{0,\frac{1}{q_{2}}-\frac{1}{q_{1}}\};
\end{split}
\end{equation}

{\rm 2)} the  space $\widetilde{B}^{l}_{p_{1},q_{1},r_{1}}(\mathbb{R}^{n},\{t^{1}_{k}\})$  is compactly embedded into
$\widetilde{B}^{l}_{p_{2},q_{2},r_{2}}(\mathbb{R}^{n},\{t^{2}_{k}\})$  if condition \eqref{4.5} is satisfied and, moreover,
\begin{equation}
\label{5.7}
\lim\limits_{j \to \infty}\biggl(\sum\limits_{m \in \mathbb{Z}^{n}}\biggl(\frac{t^{2}_{j,m}}{t^{1}_{j,m}}\biggr)^{p^{*}}\biggr)^{\frac{1}{p^{*}}}=0 \mbox{ if } q^{*}=\infty
\end{equation}
and
\begin{equation}
\label{5.8}
\lim\limits_{|m| \to \infty}\frac{t^{1}_{j,m}}{t^{2}_{j,m}}=\infty \mbox{ for all } j \in \mathbb{N}_{0} \mbox{ if } p^{*}=\infty.
\end{equation}
\end{Ca}

\section{Traces of the spaces $\widetilde{B}^{l}_{p,q,r}(\mathbb{R}^{n},\{t_{k}\})$  on planes}

In this section we assume that $n \geq 2$. Throughout the section we fix a~natural number $n' < n$ and define $n'':=n-n'$.
The point $x=(x_{1},\dots,x_{n}) \in \mathbb{R}^{n}$ will be denoted by $(x',x'')=(x'_{1},\dots,x'_{n'},x''_{n'+1},\dots,x''_{n})$
(similarly, we put $m:=(m',m'')$ for $m \in \mathbb{Z}^{n}$).
We identify the space $\mathbb{R}^{n'}$ with the plane given in the space $\mathbb{R}^{n}$ by the equation $x''=0$.

Let $\{t_{k}\} \in X^{\alpha_{3}}_{\alpha,\sigma,p}$ be a $p$-admissible  weight sequence and let $\{t_{k,m}\}$ be the $p$-associated multiple sequence.
Given $k \in \mathbb{N}_{0}$, $m' \in \mathbb{Z}^{n'}$, we set  $t'_{k,m'}=t_{k,(m',0)}$. Next, let $t'_{k}(x'):=2^{\frac{kn'}{p}}\sum\limits_{m \in \mathbb{Z}^{n'}}\chi_{\widetilde{Q}^{n'}_{k,m'}}(x')t'_{k,m'}$ for $k \in \mathbb{N}_{0}$, $x' \in \mathbb{R}^{n'}$.
So, we have $t'_{k}(x')=2^{\frac{-kn''}{p}}t_{k}(x',0)$ for $k \in \mathbb{N}_{0}$, $x' \in \mathbb{R}^{n'}$.

For $p,r \in (0,\infty]$, $p \neq \infty$, $\theta \in (0,\min\{p,r\}]$ we define $p_{\theta}:=\frac{p}{\theta}$ (as in \S\,4).
In this section it will be convenient to denote by  $\overline{p_{\theta}}$
the dual exponent to~$p_{\theta}$. In other words, $\frac{1}{p_{\theta}}+\frac{1}{\overline{p_{\theta}}}=1$.

In defining the trace of the space $\widetilde{B}^{l}_{p,q,r}(\mathbb{R}^{n},\{t_{k}\})$ we shall follow the idea of~\cite{Moura}.
(where the trace of the space $B^{\{s_{k}\}}_{p(\cdot),q}(\mathbb{R}^{n})$ was considered).

Given $l \in \mathbb{N}$ we set $\Sigma^{l}:=\bigcup\limits_{k=0}^{\infty}\Sigma^{l}_{k}$ (for the definition of $\Sigma^{l}_{k}$ see \S\,4).
Clearly, $\Sigma^{l} \subset C(\mathbb{R}^{n})$. Hence, it makes sense to talk about the pointwise trace
for a~function $f \in  \Sigma^{l} \bigcap \widetilde{B}^{l}_{p,q,r}(\mathbb{R}^{n},\{t_{k}\})$.

In other words, the function $\operatorname{tr}\left|_{x''=0}f\right.:=f(x',0)$ is well-defined.

In order to define the trace of an arbitrary function $\varphi \in \widetilde{B}^{l}_{p,q,r}(\mathbb{R}^{n},\{t_{k}\})$ we shall require the following simple result.

\begin{Lm}
\label{Lm6.1}
Let $p,q \in (0,\infty)$, $r \in (0,\infty]$, $\theta \in (0,\min\{r,p\}]$, $\alpha_{3} \geq 0$, $\alpha_{1} > n(\frac{1}{\theta}-\frac{1}{r})$, $l > \alpha_{2}$, $\sigma_{1}=\theta\overline{p_{\theta}}$, $\sigma_{2}=p$ and let $\{t_{k}\} \in X^{\alpha_{3}}_{\alpha,\sigma,p}$ be a $p$-admissible  weight sequence.
Next, assume that, for some $l' \geq l$, $r' \in (0,\infty]$, $\theta' \in (0,\min\{r',p\}]$, $\alpha'_{3} \geq 0$, $\alpha'_{1} > n(\frac{1}{\theta'}-\frac{1}{r})$, $\alpha'_{2} < l'$ $\sigma'_{1} = \theta'\overline{p_{\theta'}}$, $\sigma'_{2}=p$ and $\{t'_{k}\} \in X^{\alpha'_{3}}_{\alpha',\sigma',p}$ and any function $f \in \Sigma^{l} \bigcap \widetilde{B}^{l}_{p,q,r}(\mathbb{R}^{n},\{t_{k}\})$,
the following estimate holds $$
\|f(\cdot, 0)|\widetilde{B}^{l'}_{p,q,r'}(\mathbb{R}^{n'},\{t'_{k}\})\| \le C \|f|\widetilde{B}^{l}_{p,q,r}(\mathbb{R}^{n},\{t_{k}\})\|.
$$
in which the constant $C > 0$ is independent of the function $f$.

Then, for any function $\varphi \in \widetilde{B}^{l}_{p,q,r}(\mathbb{R}^{n},\{t_{k}\})$,
here exists a~unique
{\rm (}up to a~nullset with respect to the $n'$-dimensional Lebesgue measure{\rm )}
function $\varphi' \in \widetilde{B}^{l'}_{p,q,r'}(\mathbb{R}^{n'},\{t'_{k}\})$  such that if
$\|\varphi-\varphi_{j}|\widetilde{B}^{l}_{p,q,r}(\mathbb{R}^{n},\{t_{k}\})\| \to 0$ as $j \to \infty$  for some sequence
$\{\varphi_{j}\} \in \Sigma^{l} \bigcap \widetilde{B}^{l}_{p,q,r}(\mathbb{R}^{n},\{t_{k}\})$, then
$\|\varphi'-\varphi_{j}(\cdot,0)|\widetilde{B}^{l'}_{p,q,r'}(\mathbb{R}^{n'},\{t'_{k}\})\| \to 0$ as $j \to \infty$, and moreover,
$$
\|\varphi'|\widetilde{B}^{l'}_{p,q,r'}(\mathbb{R}^{n'},\{t'_{k}\})\| \le C \|\varphi|\widetilde{B}^{l}_{p,q,r}(\mathbb{R}^{n},\{t_{k}\})\|.
$$
\end{Lm}

The proof of this lemma repeats the corresponding arguments in~\cite{Moura} with due account of Theorem \ref{Th2.4} and Remark \ref{R4.3}.

\begin{Def}
\label{Def6.1}
Under the hypotheses of Lemma~\ref{Lm6.1} let $\varphi \in \widetilde{B}^{l}_{p,q,r}(\mathbb{R}^{n},\{t_{k}\})$.
The
function~$\varphi'$ constructed in Lemma~\ref{Lm6.1} will be called the trace of the function $\varphi$ and denoted by $
\operatorname{tr}|_{x''=0}\varphi$. By the trace of the space $\widetilde{B}^{l}_{p,q,r}(\mathbb{R}^{n},\{t_{k}\})$ on the plane, as given
in the space $\mathbb{R}^{n}$ by the equation $x''=0$, we shall mean the set of classes of equivalent functions
$\varphi' \in \widetilde{B}^{l'}_{p,q,r'}(\mathbb{R}^{n'},\{t'_{k}\})$ of which each is the trace of some function
 $\varphi \in \widetilde{B}^{l}_{p,q,r}(\mathbb{R}^{n},\{t_{k}\})$. The corresponding linear space will be denoted by $\operatorname{Tr}|_{x'' = 0}\widetilde{B}^{l}_{p,q,r}(\mathbb{R}^{n},\{t_{k}\})$;
 the norm on this space is defined as
$$
\|\varphi'\mid \operatorname{Tr}|_{x'' = 0}\widetilde{B}^{l}_{p,q,r}(\mathbb{R}^{n},\{t_{k}\})\|:=\inf\limits_{\varphi'=\operatorname{tr}|_{x''=0}\varphi}\|\varphi|\widetilde{B}^{l}_{p,q,r}(\mathbb{R}^{n},\{t_{k}\})\|.
$$
\end{Def}

In what follows under the conditions of Lemma \ref{Lm6.1} we shall also denote by $\operatorname{Tr}$ the linear operator $\operatorname{Tr}:\widetilde{B}^{l}_{p,q,r}
(\mathbb{R}^{n},\{t_{k}\}) \to \widetilde{B}^{l'}_{p,q,r'}(\mathbb{R}^{n'},\{t'_{k}\})$ defined by $\operatorname{Tr}[\varphi](x')=\operatorname{tr}|_{x''=0}\varphi(x')$ for $x' \in \mathbb{R}^{n'}$.

Recall that in \S\,4 we defined, for $k \in \mathbb{N}_{0}$, $m=(m',m'') \in \mathbb{Z}^{n}$,
$$
N^{l}_{k,m}(x):=\prod\limits_{i=1}^{n}N^{l}(2^{k}(x_{i}-\frac{m_{i}}{2^{k}})) \mbox{ for } x \in \mathbb{R}^{n},
$$
and hence
$$
N^{l}_{k,m}(x):=N^{l}_{k,m'}(x')N^{l}_{k,m''}(x'')\mbox{ for } x=(x',x'') \in \mathbb{R}^{n}.
$$

We note that for any $k \in \mathbb{N}_{0}$, $m' \in \mathbb{Z}^{n'}$,
\begin{equation}
\label{6.1}
N^{l}_{k,m'}(x')=\sum\limits_{m'' \in Z^{n''}}N^{l}_{k,(m',m'')}(x',0) \mbox{ for } x' \in \mathbb{R}^{n'}.
\end{equation}

The number of terms on the right of \eqref{6.1} is in fact finite and is bounded by some number independent of $m'$ and~$x'$.
This follows from the fact that the splines $N^{l}_{k,m}$ form a~partition of unity and that the multiplicity of intersections of
the supports of splines $N^{l}_{k,m}$ is finite (and is independent of~$m$).

Corollary \ref{Ca4.4} enables one to obtain necessary and sufficient conditions for the trace of the space $\widetilde{B}^{l}_{p,q,r}(\mathbb{R}^{n},\{t_{k}\})$.

\begin{Th}
\label{Th6.1}
Let $p,q \in (0,\infty)$, $r \in (0,\infty]$, $\theta \in (0,\min\{r,p\}]$, $\alpha_{3} \geq 0$, $\alpha_{1} > n(\frac{1}{\theta}-\frac{1}{r})$, $l > \alpha_{2}$,
$\sigma_{1}=\theta\overline{p_{\theta}}$, $\sigma_{2}=p$ and let $\{t_{k}\} \in X^{\alpha_{3}}_{\alpha,\sigma,p}$
be a~$p$-admissible  weight sequence
such that the weight sequence $\{t'_{k}\} \in \widetilde{X}^{\alpha_{3}}_{\alpha',\sigma',p}$ with  $l' \geq l$, $r' \in (0,\infty]$, $\theta' \in (0,\min\{r',p\}]$, $\alpha'_{3} \geq 0$, $\alpha'_{1} > n(\frac{1}{\theta'}-\frac{1}{r})$, $\alpha'_{2} < l'$ $\sigma'_{1} = \theta'\overline{p_{\theta'}}$, $\sigma'_{2}=p$.
 Then the operator $\operatorname{Tr}:\widetilde{B}^{l}_{p,q,r}(\mathbb{R}^{n},\{t_{k}\}) \to \widetilde{B}^{l'}_{p,q,r'}(\mathbb{R}^{n'},\{t'_{k}\})$
is bounded and there exists a~{\rm (}nonlinear{\rm )} bounded operator $\operatorname{Ext}:\widetilde{B}^{l'}_{p,q,r'}(\mathbb{R}^{n'},\{t'_{k}\})
\to \widetilde{B}^{l}_{p,q,r}(\mathbb{R}^{n},\{t_{k}\})$ such that $\operatorname{Tr} \circ \operatorname{Ext}=Id$ on the space
$\widetilde{B}^{l'}_{p,q,r'}(\mathbb{R}^{n'},\{t'_{k}\})$. In particular,
 $$
 \operatorname{Tr}|_{x'' = 0}\widetilde{B}^{l}_{p,q,r}(\mathbb{R}^{n},\{t_{k}\})= \widetilde{B}^{l'}_{p,q,r'}(\mathbb{R}^{n'},\{t'_{k}\}),
 $$
the corresponding norms being equivalent.
\end{Th}

\textbf{Proof}. Note that under the hypotheses of Theorem~\ref{Th6.1} one may apply Corollaries \ref{Ca4.3} and~\ref{Ca4.4}
to the spaces $\widetilde{B}^{l}_{p,q,r}(\mathbb{R}^{n},\{t_{k}\})$ and $\widetilde{B}^{l'}_{p,q,r'}(\mathbb{R}^{n'},\{t'_{k}\})$.

The proof is naturally split into two parts.

\textit{1}.
Let $\varphi \in \widetilde{B}^{l}_{p,q,r}(\mathbb{R}^{n},\{t_{k}\}) \bigcap \Sigma^{l}$. Then $\varphi \in \widetilde{B}^{l'}_{p,q,r}(\mathbb{R}^{n},\{t_{k}\})\bigcap \Sigma^{l}$ by
Corollary \ref{Ca4.3}, the corresponding norms being equivalent. Using Corollary~\ref{Ca4.4},
\begin{equation}
\begin{split}
\label{6.2}
\varphi=\sum\limits_{k=0}^{\infty}v^{l'}_{k}(\varphi) \mbox{ in the sense of } L_{r}^{\text{\rm loc}}(\mathbb{R}^{n}), \mbox{ where } v^{l'}_{k}(\varphi)(x)=\sum\limits_{m \in \mathbb{Z}^{n}}\alpha_{k,m}(\varphi)N^{l'}_{k,m}(x) \mbox{ for } x \in \mathbb{R}^{n}.
\end{split}
\end{equation}

Moreover,
\begin{equation}
\label{6.3}
 \biggl(\sum\limits_{k=0}^{\infty}\biggl(\sum\limits_{m \in \mathbb{Z}^{n}}t^{p}_{k,m}|\alpha_{k,m}(\varphi)|^{p}\biggr)^{\frac{q}{p}}\biggr)^{\frac{1}{q}} \le C \|\varphi|\widetilde{B}^{l'}_{p,q,r}(\mathbb{R}^{n},\{t_{k}\})\| \le C\|\varphi|\widetilde{B}^{l}_{p,q,r}(\mathbb{R}^{n},\{t_{k}\})\|.
\end{equation}

We set
$$
\alpha'_{k,m'}:=\sum_{m'' \in \mathbb{Z}^{n''}}\alpha_{k,(m',m'')}(\varphi)N^{l}_{k,m''}(0) \mbox{ for } k \in \mathbb{N}_{0}, m' \in \mathbb{Z}^{n'}.
$$

Hence,
\begin{equation}
\begin{split}
\label{6.4}
v'^{l'}_{k}(x'):=\operatorname{tr}|_{x''=0}v^{l'}_{k}(x).=\sum\limits_{m \in Z^{n}}\alpha_{k,m}(\varphi)N^{l'}_{k,m}(x',0)= \sum\limits_{m' \in Z^{n'}}\alpha'_{k,m'}N^{l'}_{k,m'}(x')    \mbox{ for } x' \in \mathbb{R}^{n'}.
\end{split}
\end{equation}

In view of \eqref{2.16} and \eqref{6.2}
\begin{equation}
\label{6.5}
|\alpha'_{k,m'}|t'_{k,m'} \le C\sum_{\substack{m'' \in \mathbb{Z}^{n''}\\
\operatorname{supp}N^{l}_{k,m}\bigcap \mathbb{R}^{n'} \neq
\emptyset}}|\alpha_{k,(m',m'')}(\varphi)|t_{k,(m',m'')}, \qquad  k \in \mathbb{N}_{0}, m' \in \mathbb{Z}^{n'}.
\end{equation}

Next, we easily obtain
\begin{equation}
\label{6.6}
 \biggl(\sum\limits_{k=0}^{\infty}\biggl(\sum\limits_{m' \in \mathbb{Z}^{n'}}t'^{p}_{k,m'}|\alpha'_{k,m'}|^{p}\biggr)^{\frac{q}{p}}\biggr)^{\frac{1}{q}} \le C  \biggl(\sum\limits_{k=0}^{\infty}\biggl(\sum\limits_{m \in \mathbb{Z}^{n}}t^{p}_{k,m}|\alpha_{k,m}(\varphi)|^{p}\biggr)^{\frac{q}{p}}\biggr)^{\frac{1}{q}}.
\end{equation}
where the constant $C>0$ is independent of of the function~$\varphi$.

Using \eqref{6.3}, \eqref{6.6} it is found by Corollary \ref{4.4} that
\begin{equation}
\begin{split}
\label{6.7}
&\|\sum\limits_{k=0}^{N}v^{l'}_{k}(\cdot,0)|\widetilde{B}^{l'}_{p,q,r'}(\mathbb{R}^{n'},\{t'_{k,m'}\})\| \le  C\biggl(\sum\limits_{k=0}^{\infty}\biggl(\sum\limits_{m' \in \mathbb{Z}^{n'}}t'^{p}_{k,m'}|\alpha'_{k,m'}|^{p}\biggr)^{\frac{q}{p}}\biggr)^{\frac{1}{q}} \le  C\|\varphi|\widetilde{B}^{l}_{p,q,r}(\mathbb{R}^{n},\{t_{k,m}\})\|.
\end{split}
\end{equation}

From \eqref{6.7} and Lemma \ref{Lm6.1} it follows that there exists the trance $\varphi'$ of the function $\varphi$ on the plane $x''=0$. Besides,
\begin{equation}
\label{6.8}
\|\varphi' |\widetilde{B}^{l'}_{p,q,r'}(\mathbb{R}^{n'},\{t'_{k}\})\| \le C \|\varphi| \widetilde{B}^{l}_{p,q,r}(\mathbb{R}^{n},\{t_{k}\})\|
\end{equation}
the constant $C>0$ being independent of the function~$\varphi$.
This proves proves the boundedness of the trace operator  $\operatorname{Tr}:\widetilde{B}^{l}_{p,q,r}(\mathbb{R}^{n},\{t_{k}\}) \to \widetilde{B}^{l'}_{p,q,r'}(\mathbb{R}^{n'},\{t'_{k}\})$.

\textit{2}.
Let $\varphi' \in \widetilde{B}^{l'}_{p,q,r'}(\mathbb{R}^{n'},\{t'_{k}\})$. By the hypotheses of the theorem and using Corollary \ref{Ca4.4}
\begin{equation}
\begin{split}
\label{6.9}
&\varphi'=\sum\limits_{k=0}^{\infty}v'^{l'}_{k}(\varphi') \mbox{ in the sense of } L_{r}^{\text{\rm loc}}(\mathbb{R}^{n'}), \mbox{ where }v'^{l'}_{k}(\varphi)(x')=\sum\limits_{m' \in \mathbb{Z}^{n'}}\alpha'_{k,m'}(\varphi)N^{l'}_{k,m'}(x'), \ \ x' \in \mathbb{R}^{n'}.
\end{split}
\end{equation}

We set
\begin{equation}\label{6.10}
\begin{gathered}
\alpha_{k,m}=\alpha'_{k,m'} \mbox{ for } m' \in \mathbb{Z}^{n'}, m'' \in \mathbb{Z}^{n''} \hbox{ and } N^{l'}_{k,m''}(0) \neq 0 \\
\alpha_{k,(m',m'')}=0 \mbox{ for } m' \in \mathbb{Z}^{n'} \hbox{ and } N^{l'}_{k,m''}(0)=0,\\
v^{l'}_{k}(x):=\sum\limits_{m \in \mathbb{Z}^{n'}}\alpha_{k,m}(\varphi)N^{l'}_{k,m}(x) \mbox{ for } x \in \mathbb{R}^{n}.
\end{gathered}
\end{equation}

Hence, using Corollaries \ref{Ca4.3}, \ref{Ca4.4}, one easily shows that the series $\sum\limits_{k=0}^{\infty}v^{l'}_{k}$ converges in $L^{\text{\rm loc}}_{r}(\mathbb{R}^{n})$
to some function $\varphi \in \widetilde{B}^{l}_{p,q,r}(\mathbb{R}^{n},\{t_{k}\})$, and moreover,
\begin{equation}
\label{6.11}
\|\varphi\mid|\widetilde{B}^{l}_{p,q,r}(\mathbb{R}^{n},\{t_{k}\})\| \le C_{1} \|\varphi \mid|\widetilde{B}^{l'}_{p,q,r}(\mathbb{R}^{n},\{t_{k}\})\| \le C_{2} \|\varphi' \mid \widetilde{B}^{l'}_{p,q,r'}(\mathbb{R}^{n'},\{t'_{k}\})\|
\end{equation}
where the constants $C_{1},C_{2}>0$ are independent of the function $\varphi'$.

We set $\operatorname{Ext}[\varphi']:=\varphi$ for $\varphi' \in \widetilde{B}^{l'}_{p,q,r'}(\mathbb{R}^{n'},\{t'_{k}\})$. Then by \eqref{6.11} the operator $\operatorname{Ext}:\widetilde{B}^{l'}_{p,q,r'}(\mathbb{R}^{n'},\{t'_{k}\}) \to \widetilde{B}^{l}_{p,q,r}(\mathbb{R}^{n},\{t_{k}\})$ is bounded.
As an immediate consequence of the construction of the function $\varphi$ we see that $\varphi'=\operatorname{tr}|_{x''=0}\varphi$, and hence
$  \operatorname{Tr} \circ \operatorname{Ext}=\hbox{Id}$ on the space $\widetilde{B}^{l'}_{p,q,r'}(\mathbb{R}^{n'},\{t'_{k}\})$.
This proves the theorem.

We illustrate Theorem \ref{Th6.1} on several examples.

\begin{Example}
\label{Ex6.1}
   Let $p,q \in [1,\infty)$, $r \in [1,p]$, a~weight sequence $\{s_{k}\} \in ^{\text{\rm loc}}Y^{\alpha_{3}}_{\alpha_{1},\alpha_{2}}=\widetilde{X}^{\alpha_{3}}_{\alpha,\infty,p}$ for $\alpha_{1}>\frac{n''}{p}$, $l>\alpha_{2}$. Then

   \begin{equation}
   \label{6.12}
    \operatorname{Tr}|_{x''=0}\widetilde{B}^{l}_{p,q,r}(\mathbb{R}^{n},\{s_{k,m}\})=\widetilde{B}^{l}_{p,q,r}(\mathbb{R}^{n'},\{s'_{k,m'}\}).
   \end{equation}

If now $\{s_{k}\} \in Y^{\alpha_{3}}_{\alpha_{1},\alpha_{2}}$, then for $\alpha_{1}>\frac{n''}{p}$, $l>\alpha_{2}$ it follows by Corollary \ref{Th2.7} that
\begin{equation}
\label{6.13}
    \operatorname{Tr}|_{x''=0}B^{\{s_{k}\}}_{p,q}(\mathbb{R}^{n})=B^{\{s'_{k}\}}_{p,q}(\mathbb{R}^{n'}).
   \end{equation}
For constant exponents $p,q \in [1,\infty)$ equality \eqref{6.13} coincides with that from~\cite{Moura}.

For $p,q \in [1,\infty)$, $r=p$, $s_{k}=2^{ks}$, $\alpha > \frac{n''}{p}$, $l > s$ we obtain  the classical result of O.\,V.~Besov (a~characterization of the trace of the
classical Besov space on the plane; see~\cite{Be1}, Theorems 1.1, 2.1,~2.2).
\end{Example}

\begin{Example}
\label{Ex6.2}
Let $p \in (1,\infty)$, $q \in [1,\infty)$, $r \in [1,p)$ and a weight $\gamma \in A^{\text{\rm loc}}_{\frac{p}{r}}(\mathbb{R}^{n})$.
We set $t_{k}(x')=\gamma_{k}(x'):=2^{k(s+\frac{n'}{p})}\sum\limits_{m' \in \mathbb{Z}^{n'}}\chi_{\widetilde{Q}^{n'}_{k,m'}}(x')\|\gamma|L_{p}(\Sigma^{n'',n'}_{k,m})\|$
for $k \in \mathbb{N}_{0}$, $x' \in \mathbb{R}^{n'}$. As a~particular case of Theorem~\ref{Th6.1} we obtain a~characterization of the trace of the weighted Besov space $B^{s}_{p,q}(\mathbb{R}^{n},\gamma)$ on the hyperplane.

Indeed, the arguments used in Example \ref{Ex2.1} we obtain
$\{\gamma_{k}\} \in X^{\alpha_{3}}_{\alpha,\sigma,p}$ for $\sigma_{1}=r\overline{p_{r}}$, $\sigma_{2}=p$, $\alpha_{1}=\alpha_{2}=s$.
Hence, using Remark \ref{R2.9} with       $s > \frac{1}{r}$, $l > s$
\begin{equation}
\label{6.14}
\operatorname{Tr}|_{x_{n}=0} B^{s}_{p,q}(\mathbb{R}^{n},\gamma)=\widetilde{B}^{l}_{p,q,r}(\mathbb{R}^{n-1},\{\gamma_{k}\}).
\end{equation}

It is worth noting that this assertion is new and may not be derived using the available atomic decomposition machinery. Indeed, the number of zero moments
for the atoms from the trace decomposition is governed by the exponent $\alpha_{1}$ for  $\{\gamma_{k}\} \in Y^{\alpha_{3}}_{\alpha_{1},\alpha_{2}}$.
The moment condition needs not be tested for  $\alpha_{1} > 0$. In a~much lesser generality an analogue of \eqref{6.14} was obtained in~\cite{HaSch},
where a~model weight depending only on the distance to the origin was examined. More precisely, $\gamma^{p}(x)=|x|^{\alpha}$ in a~small neighbourhood of
the origin with $-n+1 < \alpha < (n-1)(p-1)$. Such a~choice of the weight has enabled the authors to skip testing the zero moment condition for the
corresponding atoms from the trace decomposition.
\end{Example}

\textbf{Concluding remarks.} Consideration of the principal results obtained in this paper shows that the differences
$\delta^{l}_{r}$ may be looked upon as the most natural replacements of the differences $\overline{\Delta}^{l}_{r}$ and $\Delta^{l}$ in the definition of
weighted Besov spaces (with fairly complicated weight) and Besov spaces of variable smoothness. The exponent~$r$ proves closely related with the
exponents $\alpha_{1}$, $\sigma_{1}$.

Speaking informally, it may be stated that the worth is the behaviour in the integral sense of the variable smoothness $\{t_{k}\}$ (the exponents $\alpha_{1}$, $\sigma_{1}$ are small)
the smaller exponent~$r$ should be taken in the differences $\delta^{l}_{r}$ in order to reveal the meaningful properties of the corresponding Besov spaces of variable smoothness.
It is worth noting that in essence this idea is contained in the book~\cite{HN}.

However, the methods of \cite{HN} are capable of dealing with weighted Besov and Lizorkin--Triebel spaces with Muckenhoupt weights, but
they do not apply in the case of spaces of variable smoothness.

\textbf{Acknowledgements.} The author is indebted to Prof.\ M.\,L.~Gol'dman for valuable comments and useful discussions of the results obtained in this paper.
The author is also grateful to all the participants in the seminar ``Spaces of differentiable real multivariate functions'' and to its supervisor O.\,V.~Besov.

\end{document}